\documentclass[a4paper,12pt]{article}
\usepackage{amsmath}
\usepackage{amscd}
\usepackage[xdvi]{graphicx}
\usepackage[dvips]{color}
\usepackage{times}
\usepackage{txfonts}

\makeatletter

\@addtoreset{equation}{section}

\newcounter{def}[section]
\renewcommand{\thedef}{\stepcounter{def}\thesection.\@arabic\c@def }

\begin{document}
\setlength{\baselineskip}{24pt}
\begin{center}
\textbf{\LARGE{A spectral theory of linear operators on rigged Hilbert spaces under analyticity conditions}}
\end{center}

\setlength{\baselineskip}{14pt}

\begin{center}
Institute of Mathematics for Industry, Kyushu University, Fukuoka,
819-0395, Japan

\large{Hayato CHIBA} \footnote{E mail address : chiba@imi.kyushu-u.ac.jp}
\end{center}
\begin{center}
Jul 29, 2011; last modified Sep 12, 2014
\end{center}

\begin{center}
\textbf{Abstract}
\end{center}
A spectral theory of linear operators on rigged Hilbert spaces
(Gelfand triplets) is developed under the assumptions that 
a linear operator $T$ on a Hilbert space $\mathcal{H}$ is a perturbation of a selfadjoint operator,
and the spectral measure of the selfadjoint operator has an 
analytic continuation near the real axis in some sense.
It is shown that there exists a dense subspace $X$ of $\mathcal{H}$ such that 
the resolvent $(\lambda -T)^{-1}\phi$ of the operator $T$ has an analytic 
continuation from the lower half plane to the upper half plane as an $X'$-valued holomorphic function
for any $\phi \in X$, even when $T$ has a continuous spectrum on $\mathbf{R}$, where $X'$ is a dual space of $X$.
The rigged Hilbert space consists of three spaces $X \subset \mathcal{H} \subset X'$.
A generalized eigenvalue and a generalized eigenfunction in $X'$ are defined by using the analytic continuation
of the resolvent as an operator from $X$ into $X'$.
Other basic tools of the usual spectral theory, such as a spectrum, resolvent, Riesz projection and 
semigroup are also studied in terms of a rigged Hilbert space.
They prove to have the same properties as those of the usual spectral theory.
The results are applied to estimate asymptotic behavior of solutions of evolution equations.
\\[0.2cm]



\section{Introduction}

A spectral theory of linear operators on topological vector spaces is one of the central issues in
functional analysis.
Spectra of linear operators provide us with much information about the operators.
However, there are phenomena that are not explained by spectra.
Consider  a linear evolution equation $dx/dt = Tx$ defined by some linear operator $T$.
It is known that if the spectrum of $T$ is included in the left half plane, any solutions $x(t)$
decay to zero as $t\to \infty$ with an exponential rate,
while if there is a point of the spectrum on the right half plane, there are solutions that diverge as $t\to \infty$
(this is true at least for a sectorial operator \cite{Hen}).
On the other hand, if the spectrum set is included in the imaginary axis,
the asymptotic behavior of solutions is far from trivial;
for a finite dimensional problem, a solution $x(t)$ is a polynomial in $t$, however,
for an infinite dimensional case, a solution can decay exponentially even if the spectrum does not lie on the left half plane.
In this sense, the spectrum set does not determine the asymptotic behavior of solutions.
Such an exponential decay of a solution is known as Landau damping
in plasma physics \cite{Cra}, and is often observed for Schr\"{o}dinger operators \cite{His, Reed}.
Now it is known that such an exponential decay can be induced by resonance poles or generalized eigenvalues.

Eigenvalues of a linear operator $T$ are singularities of the resolvent $(\lambda -T)^{-1}$.
Resonance poles are obtained as singularities of a continuation of the resolvent in some sense.
In the literature, resonance poles are defined in several ways:
Let $T$ be a selfadjoint operator (for simplicity) on a Hilbert space $\mathcal{H}$ 
with the inner product $(\, \cdot \, ,\, \cdot \,)$.
Suppose that $T$ has the continuous spectrum $\sigma _c(T)$ on the real axis.
For Schr\"{o}dinger operators, spectral deformation (complex distortion) technique is often employed 
to define resonance poles \cite{His}.
A given operator $T$ is deformed by some transformation so that the continuous spectrum $\sigma _c(T)$ moves to
the upper (or lower) half plane.
Then, resonance poles are defined as eigenvalues of the deformed operator.
One of the advantages of the method is that studies of resonance poles are reduced to the usual spectral 
theory of the deformed operator on a Hilbert space.
Another way to define resonance poles is to use analytic continuations of matrix elements of the resolvent.
By the definition of the spectrum, the resolvent $(\lambda -T)^{-1}$ diverges in norm when $\lambda \in \sigma_c (T)$.
However, the matrix element $((\lambda -T)^{-1}\phi, \phi)$ for some ``good" function $\phi\in \mathcal{H}$
may exist for $\lambda \in \sigma_c (T)$, and the function $f(\lambda ) = ((\lambda -T)^{-1}\phi, \phi)$ may have an 
analytic continuation from the lower half plane to the upper half plane through an interval on $\sigma _c(T)$.
Then, the analytic continuation may have poles on the upper half plane, which is called a resonance pole or 
a generalized eigenvalue.
In the study of reaction-diffusion equations, the Evans function is often used, whose zeros give eigenvalues of a given differential operator.
Resonance poles can be defined as zeros of an analytic continuation of the Evans function \cite{ZumHow}.
See \cite{His, Rau, Reed} for other definitions of resonance poles.

Although these methods work well for some special classes of Schr\"{o}dinger operators,
an abstract spectral theory of resonance poles has not been developed well.
In particular, a precise definition of an eigenfunction associated with a resonance pole is not obvious in general.
Clearly a pole of a matrix element or the Evans function does not provide an eigenfunction.
In Chiba \cite{Chi}, a definition of the eigenfunction associated with a resonance pole is suggested for a certain
operator obtained from the Kuramoto model (see Sec.4).
It is shown that the eigenfunction is a distribution, not a usual function.
This suggests that an abstract theory of topological vector spaces should be employed for the study of 
a resonance pole and its eigenfunction of an abstract linear operator.

The purpose in this paper is to give a correct formulation of resonance poles and eigenfunctions in terms of operator theory
on rigged Hilbert spaces (Gelfand triplets).
Our approach based on rigged Hilbert spaces allows one to develop a spectral theory of resonance poles
in a parallel way to ``standard course of functional analysis".
To explain our idea based on rigged Hilbert spaces, let us consider the multiplication operator 
$\mathcal{M} : \phi (\omega ) \mapsto \omega \phi (\omega )$ on the Lebesgue space $L^2(\mathbf{R})$.
The resolvent is given as
\begin{eqnarray*}
((\lambda -\mathcal{M})^{-1}\phi, \psi^*) = \int_{\mathbf{R}}\! \frac{1}{\lambda -\omega }\phi (\omega ) \psi (\omega ) d\omega, 
\end{eqnarray*}
where $\psi^* = \overline{\psi (\omega )}$, which is employed to avoid the complex conjugate of $\psi (\omega )$
in the right hand side.
This function of $\lambda $ is holomorphic on the lower half plane, and it does not exist for $\lambda \in \mathbf{R}$;
the continuous spectrum of $\mathcal{M}$ is the whole real axis.
However, if $\phi$ and $\psi$ have analytic continuations near the real axis, the right hand side has 
an analytic continuation from the lower half plane to the upper half plane, which is given by
\begin{eqnarray*}
\int_{\mathbf{R}}\! \frac{1}{\lambda -\omega }\phi (\omega ) \psi (\omega ) d\omega  + 2\pi \mathrm{i} \phi (\lambda )\psi (\lambda ),
\end{eqnarray*}
where $\mathrm{i} := \sqrt{-1}$.
Let $X$ be a dense subspace of $L^2(\mathbf{R})$ consisting of functions having analytic continuations near the real axis.
A mapping, which maps $\phi \in X$ to the above value, defines a continuous linear functional on $X$, that is,
an element of the dual space $X'$, if $X$ is equipped with a suitable topology.
Motivated by this idea, we define the linear operator $A(\lambda ) : X\to X'$ to be
\begin{equation}
\langle A(\lambda )\psi \,|\, \phi \rangle = \left\{ \begin{array}{ll}
\displaystyle 
     \int_{\mathbf{R}}\! \frac{1}{\lambda -\omega  } \psi (\omega ) \phi (\omega ) d\omega  
           + 2\pi \mathrm{i} \psi (\lambda ) \phi (\lambda ) & (\mathrm{Im}(\lambda ) > 0), \\[0.4cm]
\displaystyle  \lim_{y\to -0} 
  \int_{\mathbf{R}}\! \frac{1}{x + \mathrm{i} y -\omega } \psi (\omega ) \phi (\omega ) d\omega  & (x = \lambda \in \mathbf{R}), \\[0.4cm]
\displaystyle  \int_{\mathbf{R}}\! \frac{1}{\lambda -\omega } \psi (\omega ) \phi (\omega ) d\omega
 & (\mathrm{Im}(\lambda ) < 0),
\end{array} \right.
\label{1-1}
\end{equation}
for $\psi, \phi\in X$, where $\langle \, \cdot  \,|\, \cdot \, \rangle$ is a paring for $(X', X)$.
When $\mathrm{Im} (\lambda ) < 0$, $A(\lambda ) = (\lambda -\mathcal{M})^{-1}$, while when $\mathrm{Im} (\lambda ) \geq 0$,
$A(\lambda )\psi$ is not included in $L^2(\mathbf{R})$ but an element of $X'$.
In this sense, $A(\lambda )$ is called the analytic continuation of the resolvent of $\mathcal{M}$
in the generalized sense.
In this manner, the triplet $X\subset L^2(\mathbf{R}) \subset X'$, which is called the rigged Hilbert space 
or the Gelfand triplet \cite{Gel1, Mau}, is introduced.

In this paper, a spectral theory on a rigged Hilbert space is proposed for an operator of the form $T = H+K$,
where $H$ is a selfadjoint operator on a Hilbert space $\mathcal{H}$,
whose spectral measure has an analytic continuation near the real axis,
when the domain is restricted to some dense subspace $X$ of $\mathcal{H}$, as above.
$K$ is an operator densely defined on $X$ satisfying certain boundedness conditions.
Our purpose is to investigate spectral properties of the operator $T = H+K$.
At first, the analytic continuation $A(\lambda )$ of the resolvent $(\lambda -H)^{-1}$
is defined as an operator from $X$ into $X'$ in the same way as Eq.(\ref{1-1}).
In general, $A(\lambda ) : X \to X'$ is defined on a nontrivial Riemann surface of $\lambda $ so that 
when $\lambda $ lies on the original complex plane, it coincides with the usual resolvent $(\lambda -H)^{-1}$.
The usual eigen-equation $(\lambda -T)v = 0$ is rewritten as
\begin{eqnarray*}
(\lambda -H) \circ (id - (\lambda -H)^{-1}K)v = 0.
\end{eqnarray*}
By neglecting the first factor and replacing $(\lambda -H)^{-1}$ by its analytic continuation $A(\lambda )$,
we arrive at the following definition:
If the equation
\begin{equation}
(id - A(\lambda )K^\times) \mu = 0
\end{equation}
has a nonzero solution $\mu$ in $X'$, such a $\lambda $ is called a generalized eigenvalue (resonance pole)
and $\mu$ is called a generalized eigenfunction,
where $K^\times : X' \to X'$ is a dual operator of $K$.
When $\lambda $ lies on the original complex plane, the above equation is reduced to the usual eigen-equation.
In this manner, resonance poles and corresponding eigenfunctions are naturally obtained
without using spectral deformation technique or poles of matrix elements.

Similarly, the resolvent in the usual sense is given by
\begin{eqnarray*}
(\lambda -T)^{-1} =  (\lambda -H)^{-1}\circ (id - K(\lambda -H)^{-1})^{-1}.
\end{eqnarray*}
Motivated by this, an analytic continuation of the resolvent of $T$ in the generalized sense is defined to be
\begin{equation}
\mathcal{R}_\lambda = A(\lambda ) \circ (id - K^\times A(\lambda ))^{-1} : X\to X',
\end{equation}
(the operator $K^\times A(\lambda )$ is well defined because of the assumption (X8) below).
When $\lambda $ lies on the original complex plane, this is reduced to the usual resolvent $(\lambda -T)^{-1}$.
With the aid of the generalized resolvent $\mathcal{R}_\lambda $,
basic concepts in the usual spectral theory, such as eigenspaces, algebraic multiplicities,
point/continuous/residual spectra, Riesz projections are extended to those defined on a rigged Hilbert space.
It is shown that they have the same properties as the usual theory.
For example, the generalized Riesz projection $\Pi_0$ for an isolated resonance pole $\lambda _0$ is defined by
the contour integral of the generalized resolvent.
\begin{equation}
\Pi_0 = \frac{1}{2\pi \mathrm{i} } \int_{\gamma }\! \mathcal{R}_\lambda d\lambda  : X\to X'.
\end{equation}
Properties of the generalized Riesz projection $\Pi_0$ is investigated in detail.
Note that in the most literature, the eigenspace associated with a resonance pole is \textit{defined} to be 
the range of the Riesz projection.
In this paper, the eigenspace of a resonance pole is defined as the set of solutions of the eigen-equation,
and it is \textit{proved} that it coincides with the range of the Riesz projection as the standard functional analysis.
Any function $\phi \in X$ proves to be uniquely decomposed as $\phi = \mu_1 + \mu_2$, where
$\mu_1 \in \Pi_0  X$ and $\mu_2 = (id - \Pi_0)X$, both of which are elements of $X'$.
These results play an important role when applying the theory to dynamical systems \cite{Chi}.
The generalized Riesz projection around a resonance pole $\lambda _0$ on the left half plane (resp.
on the imaginary axis) defines a stable subspace (resp. a center subspace) in the generalized sense,
both of which are subspaces of $X'$.
Then, the standard idea of the dynamical systems theory may be applied to investigate the asymptotic behavior and bifurcations of an 
infinite dimensional dynamical system.
Such a dynamics induced by a resonance pole is not captured by the usual eigenvalues.

Many properties of the generalized spectrum (the set of singularities of $\mathcal{R}_\lambda $) will be shown.
In general, the generalized spectrum consists of the generalized point spectrum (the set of resonance poles),
the generalized continuous spectrum and the generalized residual spectrum 
(they are not distinguished in the literature).
If the operator $K$ satisfies a certain compactness condition, the Riesz-Schauder theory on a rigged Hilbert space
applies to conclude that the generalized spectrum consists only of a countable number of resonance poles 
having finite multiplicities.
It is remarkable that even if the operator $T$ has the continuous spectrum (in the usual sense),
the generalized spectrum consists only of a countable number of
resonance poles when $K$ satisfies the compactness condition.
Since the topology on the dual space $X'$ is weaker than that on the Hilbert space $\mathcal{H}$,
the continuous spectrum of $T$ disappears, while eigenvalues remain to exist as the generalized spectrum.
This fact is useful to estimate embedded eigenvalues.
Eigenvalues embedded in the continuous spectrum is no longer embedded in our spectral theory.
Thus, the Riesz projection is applicable to obtain eigenspaces of them.
Our theory is also used to estimate an exponential decay of the semigroup $e^{\mathrm{i} Tt}$ generated by $\mathrm{i} T$.
It is shown that resonance poles induce an exponential decay of the semigroup 
even if the operator $\mathrm{i} T$ has no spectrum on the left half plane.

Although resonance poles have been well studied for Schr\"{o}dinger operators,
a spectral theory in this paper is motivated by establishing bifurcation theory of 
infinite dimensional dynamical systems, for which spectral deformation technique is not applied.
In Chiba \cite{Chi}, a bifurcation structure of an infinite dimensional coupled oscillators (Kuramoto model) is investigated
by means of rigged Hilbert spaces.
It is shown that when a resonance pole of a certain linear operator, which is obtained by the linearization of the system around a steady state,
gets across the imaginary axis as a parameter of the system varies, then a bifurcation occurs.
For this purpose, properties of generalized eigenfunctions developed in this paper play an important role.
In Section 4 of the present article, the linear stability analysis of the Kuramoto model will be given
to demonstrate how our new theory is applied to the study of dynamical systems.
In particular, a spectral decomposition theorem of a certain non-selfadjoint non-compact operator
will be proved, which seems not to be obtained by the classical theory of resonance poles.

Throughout this paper, $\mathsf{D}(\cdot)$ and $\mathsf{R}(\cdot)$ denote the domain and range of an operator,
respectively.


\section{Spectral theory on a Hilbert space}

This section is devoted to a review of the spectral theory of a perturbed selfadjoint operator on a Hilbert space
to compare the spectral theory on a rigged Hilbert space developed after Sec.3.
Let $\mathcal{H}$ be a Hilbert space over $\mathbf{C}$.
The inner product is defined so that
\begin{equation}
(a \phi, \psi) = (\phi, \overline{a} \psi) = a (\phi, \psi), 
\label{2-1}
\end{equation}
where $\overline{a}$ is the complex conjugate of $a\in \mathbf{C}$.
Let us consider an operator $T := H  + K$ defined on a dense subspace of $\mathcal{H}$, 
where $H$ is a selfadjoint operator, and $K$ is a compact operator on $\mathcal{H}$ which need not be selfadjoint.
Let $\lambda $ and $v = v_\lambda $ be an eigenvalue and an eigenfunction, respectively, of the operator $T$
defined by the equation $\lambda v = H v + Kv$.
This is rearranged as 
\begin{equation}
(\lambda -H) (id - (\lambda -H)^{-1}K)v = 0,
\label{2-4}
\end{equation}
where $id$ denotes the identity on $\mathcal{H}$.
In particular, when $\lambda $ is not an eigenvalue of $H$,
it is an eigenvalue of $T$ if and only if $id - (\lambda -H )^{-1}K$ is not injective in $\mathcal{H}$.
Since the essential spectrum is stable under compact perturbations (see Kato \cite{Kato}, Theorem IV-5.35), 
the essential spectrum $\sigma_e(T)$ of $T$ is the same as that of $H$, which lies on the real axis.
Since $K$ is a compact perturbation, the Riesz-Schauder theory shows that the spectrum outside the real axis consists of the discrete spectrum; 
for any $\delta >0$, the number of eigenvalues satisfying $|\mathrm{Im}(\lambda )| \geq \delta $
is finite, and their algebraic multiplicities are finite.
Eigenvalues may accumulate only on the real axis.
To find eigenvalues embedded in the essential spectrum $\sigma_e(T)$ is a difficult and important problem.
In this paper, a new spectral theory on rigged Hilbert spaces will be developed to obtain such embedded 
eigenvalues and corresponding eigenspaces.
\\[-0.2cm]

Let $R_\lambda = (\lambda -T)^{-1}$ be the resolvent.
Let $\lambda _j$ be an eigenvalue of $T$ outside the real axis, and $\gamma _j$ be a simple
closed curve enclosing $\lambda _j$ separated from the rest of the spectrum.
The projection to the generalized eigenspace $V_j := \bigcup_{n\geq 1} \mathrm{Ker} (\lambda _j - T)^n$
is given by
\begin{equation}
\Pi_j = \frac{1}{2\pi \mathrm{i} } \int_{\gamma _j}\! R_\lambda d\lambda . 
\end{equation}
\\[-0.2cm]

Let us consider the semigroup $e^{\mathrm{i} Tt}$ generated by $\mathrm{i} T$.
Since $\mathrm{i} H$ generates the $C^0$-semigroup $e^{\mathrm{i} H t}$ and 
$K$ is compact, $\mathrm{i} T$ also generates the $C^0$-semigroup (see Kato \cite{Kato}, Chap.IX).
It is known that $e^{\mathrm{i} Tt}$ is obtained by the Laplace inversion formula (Hille and Phillips \cite{Hil}, Theorem 11.6.1)
\begin{equation}
e^{\mathrm{i} Tt}\phi = \frac{1}{2\pi \mathrm{i} } \lim_{x\to \infty} \int^{x-\mathrm{i} y}_{-x-\mathrm{i} y}\!
e^{\mathrm{i} \lambda t} (\lambda -T)^{-1}\phi d\lambda , \quad x,y \in \mathbf{R}, 
\label{2-16}
\end{equation}
for $t>0 $ and $\phi \in \mathsf{D}(T)$, where $y > 0$ is chosen so that all eigenvalues $\lambda $ of $T$ satisfy $\mathrm{Im}(\lambda ) > -y$,
and the limit $x\to \infty$ exists with respect to the topology of $\mathcal{H}$.
Thus the contour is the horizontal line on the lower half plane.
Let $\varepsilon >0$ be a small number and $\lambda _0, \cdots ,\lambda _N$ eigenvalues of $T$
satisfying $\mathrm{Im}(\lambda _j) \leq -\varepsilon ,\, j=0, \cdots ,N$.
The residue theorem provides
\begin{eqnarray*}
e^{\mathrm{i} Tt}\phi &=& \frac{1}{2\pi \mathrm{i} } \int_{\mathbf{R}}\!
e^{\mathrm{i} xt + \varepsilon t} (x - \mathrm{i} \varepsilon  -T)^{-1}\phi dx \\
& & \quad + \frac{1}{2\pi \mathrm{i} } \sum^N_{j=0} \int_{\gamma _j}\! e^{\mathrm{i} \lambda t} (\lambda -T)^{-1}\phi d\lambda,
\end{eqnarray*}
where $\gamma _j$ is a sufficiently small closed curve enclosing $\lambda _j$.
Let $M_j$ be the smallest integer such that $(\lambda _j - T)^{M_j} \Pi_j = 0$.
This is less or equal to the algebraic multiplicity of $\lambda _j$.
Then, $e^{\mathrm{i} Tt}$ is calculated as
\begin{eqnarray*}
e^{\mathrm{i} Tt}\phi &=& \frac{1}{2\pi \mathrm{i} } \int_{\mathbf{R}}\!
e^{\mathrm{i} xt + \varepsilon t} (x - \mathrm{i} \varepsilon  -T)^{-1}\phi dx \\
& &  \quad + 
\sum^N_{j=0} e^{\mathrm{i} \lambda _j t} \sum^{M_j-1}_{k=0} \frac{(-\mathrm{i} t)^k}{k!}(\lambda _j - T)^k \Pi_j\phi .
\end{eqnarray*}
The second term above diverges as $t\to \infty$ because $\mathrm{Re}(\mathrm{i} \lambda _j) \geq \varepsilon $.
On the other hand, if there are no eigenvalues on the lower half plane, 
we obtain
\begin{eqnarray*}
e^{\mathrm{i} Tt}\phi &=& \frac{1}{2\pi \mathrm{i} } \int_{\mathbf{R}}\!
e^{\mathrm{i} xt + \varepsilon t} (x - \mathrm{i} \varepsilon  -T)^{-1}\phi dx,
\end{eqnarray*}
for any small $\varepsilon >0$.
In such a case, the asymptotic behavior of $e^{\mathrm{i} Tt}$ is quite nontrivial.
One of the purposes in this paper is to give a further decomposition of the first term above 
under certain analyticity conditions to determine the dynamics of $e^{\mathrm{i} Tt}$.


\section{Spectral theory on a Gelfand triplet}

In the previous section, we give the review of the spectral theory of the operator $T = H+K$ on $\mathcal{H}$.
In this section, the notion of spectra, eigenfunctions, resolvents and projections are extended by means of a rigged Hilbert space.
It will be shown that they have similar properties to those on $\mathcal{H}$.
They are used to estimate the asymptotic behavior of the semigroup $e^{\mathrm{i} Tt}$
and to find embedded eigenvalues.

\subsection{Rigged Hilbert spaces}

Let $X$ be a locally convex Hausdorff topological vector space over $\mathbf{C}$ and $X'$ its dual space.
$X'$ is a set of continuous anti-linear functionals on $X$.
For $\mu \in X'$ and $\phi \in X$, $\mu (\phi )$ is denoted by
$\langle \mu \,|\, \phi \rangle$.
For any $a,b \in \mathbf{C},\, \phi, \psi \in X$ and $\mu, \xi \in X'$, the equalities
\begin{eqnarray}
& & \langle \mu \,|\,  a \phi + b\psi\rangle 
   = \overline{a} \langle \mu \,|\,  \phi \rangle + \overline{b} \langle \mu \,|\, \psi \rangle, \\
& & \langle a\mu + b\xi \,|\, \phi \rangle
   = a \langle \mu \,|\, \phi \rangle + b \langle \xi \,|\, \phi \rangle,
\end{eqnarray}
hold. 
In this paper, an element of $X'$ is called a generalized function \cite{Gel0, Gel1}.
Several topologies can be defined on the dual space $X'$.
Two of the most usual topologies are the weak dual topology (weak * topology) 
and the strong dual topology (strong * topology).
A sequence $\{ \mu_j\} \subset X'$ is said to be weakly convergent to $\mu \in X'$
if $\langle \mu_j \,|\, \phi \rangle \to \langle \mu  \,|\, \phi \rangle$ for each $\phi \in X$;
a sequence $\{ \mu_j\} \subset X'$ is said to be strongly convergent to $\mu \in X'$
if $\langle \mu_j \,|\,\phi \rangle \to \langle \mu \,|\, \phi \rangle$ uniformly on any bounded subset of $X$.

Let $\mathcal{H}$ be a Hilbert space with the inner product $(\cdot\, , \, \cdot)$ such that $X$ is a dense subspace of
$\mathcal{H}$.
Since a Hilbert space is isomorphic to its dual space, we obtain $\mathcal{H} \subset X'$ through $\mathcal{H} \simeq \mathcal{H}'$.
\\[0.2cm]
\textbf{Definition \thedef.} If a locally convex Hausdorff topological vector space $X$ is a dense subspace of 
a Hilbert space $\mathcal{H}$ and a topology of $X$ is stronger than that of $\mathcal{H}$, the triplet
\begin{equation}
X \subset \mathcal{H} \subset X'
\end{equation}
is called the \textit{rigged Hilbert space} or the \textit{Gelfand triplet}.
The \textit{canonical inclusion} $i: X \to X'$ is defined as follows; for $\psi\in X$,
we denote $i(\psi)$ by $\langle \psi |$, which is defined to be
\begin{equation}
i(\psi)(\phi) = \langle \psi \,|\, \phi \rangle = (\psi, \phi),
\label{3-4}
\end{equation}
for any $\phi \in X$ (note that we also use $\mathrm{i} = \sqrt{-1}$). 
The inclusion from $\mathcal{H}$ into $X'$ is also defined as above.
It is easy to show that the canonical inclusion is injective if and only if 
$X$ is a dense subspace of $\mathcal{H}$, and the canonical inclusion is continuous 
(for both of the weak dual topology and the strong dual topology) if and only if
a topology of $X$ is stronger than that of $\mathcal{H}$ (see Tr\'{e}ves \cite{Tre}).

A topological vector space $X$ is called Montel if it is barreled and every bounded set of $X$ is relatively compact.
A Montel space has a convenient property that on a bounded set $A$ of a dual space of a Montel space,
the weak dual topology coincides with the strong dual topology.
In particular, a weakly convergent series in a dual of a Montel space also converges with respect to 
the strong dual topology (see Tr\'{e}ves \cite{Tre}).
Furthermore, a linear map from a topological vector space to a Montel space is a compact operator
if and only if it is a bounded operator.
It is known that the theory of rigged Hilbert spaces works best when the space $X$ is a Montel or a nuclear space
\cite{Gel1}.
See Grothendieck \cite{Gro} and Komatsu \cite{Kom} for sufficient conditions for a topological
vector space to be a Montel space or a nuclear space.


\subsection{Generalized eigenvalues and eigenfunctions}

\begin{figure}
\begin{center}
\includegraphics{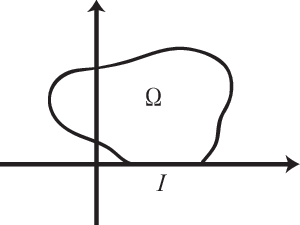}
\caption{A domain on which $E[\psi, \phi](\omega )$ is holomorphic. \label{fig1}}
\end{center}
\end{figure}

Let $\mathcal{H}$ be a Hilbert space over $\mathbf{C}$ and $H$ a selfadjoint operator densely defined on $\mathcal{H}$
with the spectral measure $\{ E(B)\}_{B\in \mathcal{B}}$; that is, $H$ is expressed as
$H = \int_{\mathbf{R}}\! \omega dE(\omega )$.
Let $K$ be some linear operator densely defined on $\mathcal{H}$.
Our purpose is to investigate spectral properties of the operator $T:=H+K$.
Let $\Omega \subset \mathbf{C}$ be a simply connected open domain in the upper half plane such that the intersection of the
real axis and the closure of $\Omega $ is a connected interval $\tilde{I}$.
Let $I = \tilde{I}\backslash \partial \tilde{I}$ be an open interval (see Fig.\ref{fig1}).
For a given $T = H+K$, we suppose that there exists a locally convex Hausdorff vector space $X(\Omega )$ 
over $\mathbf{C}$ satisfying following conditions.
\\[0.2cm]
\textbf{(X1)} $X(\Omega )$ is a dense subspace of $\mathcal{H}$.
\\
\textbf{(X2)} A topology on $X(\Omega )$ is stronger than that on $\mathcal{H}$.
\\
\textbf{(X3)} $X(\Omega )$ is a quasi-complete barreled space.
\\
\textbf{(X4)} For any $\phi \in X(\Omega )$, the spectral measure $(E(B)\phi, \phi)$ is absolutely continuous on the interval $I$.
Its density function, denoted by $E[\phi, \phi](\omega )$, has an analytic continuation to $\Omega \cup I$.
\\
\textbf{(X5)} For each $\lambda \in I \cup \Omega $, the bilinear form
$E[\, \cdot\,, \, \cdot \,](\lambda ): X(\Omega ) \times X(\Omega ) \to \mathbf{C}$ is separately continuous
(i.e. $E[\, \cdot\,, \, \phi \,](\lambda ): X(\Omega ) \to \mathbf{C}$ and $E[\, \phi \,, \, \cdot \,](\lambda ): X(\Omega ) \to \mathbf{C}$
are continuous for fixed $\phi \in X(\Omega )$).
\\[0.2cm]
Because of (X1) and (X2), the rigged Hilbert space $X(\Omega ) \subset \mathcal{H} \subset X(\Omega )'$ is well defined,
where $X(\Omega )'$ is a space of continuous \textit{anti}-linear functionals and
the canonical inclusion $i$ is defined by Eq.(\ref{3-4}).
Sometimes we denote $i(\psi)$ by $\psi$ for simplicity by identifying $iX(\Omega )$ with $X(\Omega )$.
The assumption (X3) is used to define Pettis integrals and Taylor expansions of $X(\Omega )'$-valued
holomorphic functions in Sec.3.5
(refer to Tr\'{e}ves \cite{Tre} for basic terminology of topological vector spaces such as quasi-complete and barreled space.
In this paper, to understand precise definitions of them is not so important;
it is sufficient to know that an integral and holomorphy of $X(\Omega )'$-valued functions are well-defined if 
$X(\Omega )$ is quasi-complete barreled. See Appendix for more detail).
For example, Montel spaces, Fr\'{e}chet spaces, Banach spaces and Hilbert spaces are barreled.
Due to the assumption (X4) with the aid of the polarization identity, 
we can show that $(E(B)\phi, \psi)$ is absolutely continuous on $I$ for any $\phi, \psi \in X(\Omega )$.
Let $E[\phi, \psi] (\omega )$ be the density function;
\begin{equation}
d(E(\omega )\phi, \psi) = E[\phi, \psi] (\omega ) d\omega , \quad \omega \in I.
\label{3-8}
\end{equation}
Then, $E[\phi, \psi] (\omega )$ is holomorphic in $\omega \in I\cup \Omega $.
We will use the above notation for any $\omega \in \mathbf{R}$ for simplicity, although the absolute continuity is
assumed only on $I$.
Since $E[\phi, \psi] (\omega )$ is absolutely continuous on $I$, $H$ is assumed not to have eigenvalues on $I$.
(X5) is used to prove the continuity of a certain operator (Prop.3.7).

Let $A$ be a linear operator densely defined on $X(\Omega )$.
Then, the dual operator $A'$ is defined as follows:
the domain $\mathsf{D}(A')$ is the set of elements $\mu \in X(\Omega )'$ such that the mapping
$\phi \mapsto \langle \mu \,|\, A\phi \rangle$ from $\mathsf{D}(A) \subset X(\Omega )$ into $\mathbf{C}$ is continuous.
Then, $A' : \mathsf{D}(A')\to X(\Omega )'$ is defined by
\begin{equation}
\langle A' \mu \,|\, \phi \rangle = \langle \mu \,|\,A\phi \rangle,\quad
\phi \in \mathsf{D}(A),\, \mu\in \mathsf{D}(A').
\end{equation}
If $A$ is continuous on $X(\Omega )$, then $A'$ is continuous on $X(\Omega )'$
for both of the weak dual topology and the strong dual topology.
The (Hilbert) adjoint $A^*$ of $A$ is defined through $(A\phi, \psi) = (\phi , A^* \psi)$ as usual
when $A$ is densely defined on $\mathcal{H}$.
\\[0.2cm]
\textbf{Lemma \thedef.} Let $A$ be a linear operator densely defined on $\mathcal{H}$.
Suppose that there exists a dense subspace $Y$ of $X(\Omega )$ such that $A^*Y \subset X(\Omega )$
so that the dual $({A^*})'$ is defined. 
Then, $(A^*)' $ is an extension of $A$ and $i\circ A=(A^*)'\circ i\,|_{\mathsf{D}(A)}$.
In particular, $\mathsf{D}((A^*)') \supset i\mathsf{D}(A)$.
\\[0.2cm]
\textbf{Proof.} By the definition of the canonical inclusion $i$, we have
\begin{eqnarray}
i(A \psi) (\phi) = (A\psi, \phi) = (\psi , A^* \phi)
 = \langle \psi  \,|\, A^*\phi \rangle  = \langle (A^*)' \psi \,|\, \phi \rangle,
\label{3-10}
\end{eqnarray}
for any $\psi \in \mathsf{D}(A)$ and $\phi \in Y$. \hfill $\blacksquare$
\\[-0.2cm]

In what follows, we denote $(A^*)'$ by $A^\times$.
Thus Eq.(\ref{3-10}) means $i\circ A = A^\times \circ i\,|_{\mathsf{D}(A)}$.
Note that $A^\times = A'$ when $A$ is selfadjoint.
For the operators $H$ and $K$, we suppose that
\\[0.2cm]
\textbf{(X6)} there exists a dense subspace $Y$ of $X(\Omega )$ such that $HY \subset X(\Omega )$.
\\
\textbf{(X7)} $K$ is $H$-bounded and $K^*Y \subset X(\Omega )$.
\\
\textbf{(X8)} $K^\times A(\lambda ) iX(\Omega ) \subset iX(\Omega )$ for any $\lambda \in \{ \mathrm{Im}(\lambda )<0\}
\cup I \cup \Omega $.
\\[0.2cm]
The operator $A(\lambda ) : iX(\Omega ) \to X(\Omega )'$ will be defined later.
Recall that when $K$ is $H$-bounded (relatively bounded with respect to $H$),
$\mathsf{D}(T) = \mathsf{D}(H)$ and $K(\lambda -H)^{-1}$ is bounded on $\mathcal{H}$ for $\lambda \notin \mathbf{R}$.
In some sense, (X8) is a ``dual version" of this condition because $A(\lambda )$ proves to be an extension of $(\lambda -H)^{-1}$.
In particular, we will show that $K^\times A(\lambda ) i = i (K(\lambda -H)^{-1})$ when $\mathrm{Im} (\lambda ) < 0$.
Our purpose is to investigate the operator $T = H  + K$ with these conditions.
Due to (X6) and (X7), the dual operator $T^\times$ of $T^* = H + K^*$ is well defined.
It follows that $\mathsf{D}(T^\times) = \mathsf{D}(H^\times) \cap \mathsf{D}(K^\times)$ and 
\begin{eqnarray*}
\mathsf{D}(T^\times) \supset i\mathsf{D}(T) =  i\mathsf{D}(H) \supset iY.
\end{eqnarray*}
In particular, the domain of $T^\times$ is dense in $X(\Omega )'$.

To define the operator $A(\lambda )$, we need the next lemma.
\\[0.2cm]
\textbf{Lemma \thedef.} Suppose that a function $q(\omega )$ is integrable on $\mathbf{R}$
and holomorphic on $\Omega \cup I$. Then, the function
\begin{equation}
Q(\lambda ) = \left\{ \begin{array}{ll}
\displaystyle \int_{\mathbf{R}}\! \frac{1}{\lambda -\omega } q(\omega )d\omega & (\mathrm{Im}(\lambda ) < 0), \\[0.4cm]
\displaystyle \int_{\mathbf{R}}\! \frac{1}{\lambda -\omega } q(\omega )d\omega + 2\pi \mathrm{i} q(\lambda ) & 
(\lambda \in \Omega ), \\
\end{array} \right.
\label{3-11}
\end{equation}
is holomorphic on $\{ \lambda \, | \, \mathrm{Im}(\lambda ) < 0\} \cup \Omega \cup I$.
\\[0.2cm]
\textbf{Proof.} Putting $\lambda  = x + \mathrm{i} y$ with $x, y \in \mathbf{R}$ yields
\begin{eqnarray*}
\int_{\mathbf{R}}\! \frac{1}{\lambda -\omega } q(\omega )d\omega
 = \int_{\mathbf{R}}\! \frac{x-\omega }{(x-\omega )^2 + y^2} q(\omega ) d\omega 
 - \mathrm{i}  \int_{\mathbf{R}}\! \frac{y}{(x-\omega )^2 + y^2} q(\omega ) d\omega .  
\end{eqnarray*}
Due to the formula of the Poisson kernel, the equalities
\begin{eqnarray*}
\lim_{y\to +0} \int_{\mathbf{R}}\! \frac{y}{(x-\omega )^2 + y^2} q(\omega ) d\omega = \pi q(x),
\quad \lim_{y\to -0} \int_{\mathbf{R}}\! \frac{y}{(x-\omega )^2 + y^2} q(\omega ) d\omega = -\pi q(x),
\end{eqnarray*}
hold when $q$ is continuous at $x\in I$ (Ahlfors \cite{Ahl}). Thus we obtain
\begin{eqnarray*}
\lim_{y\to -0} \int_{\mathbf{R}}\! \frac{1}{\lambda -\omega } q(\omega )d\omega 
 = \lim_{y\to +0} \left( 
        \int_{\mathbf{R}}\! \frac{1}{\lambda -\omega } q(\omega )d\omega + 2\pi \mathrm{i} q(\lambda )\right)
 = \pi V(x) + \pi \mathrm{i} q(x),
\end{eqnarray*}
where
\begin{eqnarray*}
V(x) := \lim_{y\to 0} \frac{1}{\pi}\int_{\mathbf{R}}\! \frac{x-\omega }{(x-\omega )^2 + y^2} q(\omega ) d\omega
\end{eqnarray*}
is the Hilbert transform of $q$.
It is known that $V(x)$ is Lipschitz continuous on $I$ if $q(x)$ is (see Titchmarsh \cite{Tit}).
Therefore, two holomorphic functions in Eq.(\ref{3-11}) coincide with one another on $I$ and they are continuous on $I$.
This proves that $Q(\lambda )$ is holomorphic on $\{ \lambda \, | \, \mathrm{Im}(\lambda ) < 0\} \cup \Omega \cup I$.
\hfill $\blacksquare$
\\[-0.2cm]

Put $u_\lambda  = (\lambda - H)^{-1}\psi$ for $\psi \in \mathcal{H}$.
In general, $u_\lambda $ is not included in $\mathcal{H}$ when $\lambda \in I$ because of the continuous spectrum of $H$.
Thus $u_\lambda $ does not have an analytic continuation from the lower half plane to $\Omega $
with respect to $\lambda $ as an $\mathcal{H}$-valued function.
To define an analytic continuation of $u_\lambda $, we regard it as a generalized function in $X(\Omega )'$ by the 
canonical inclusion. 
Then, the action of $i((\lambda -H )^{-1}\psi)$ is given by
\begin{eqnarray*}
i((\lambda -H )^{-1}\psi) (\phi) = ((\lambda  -H )^{-1} \psi, \phi)
 = \int_{\mathbf{R}}\! \frac{1}{\lambda -\omega } E [\psi, \phi](\omega ) d\omega , \quad \mathrm{Im}(\lambda ) < 0. 
\end{eqnarray*}
Because of the assumption (X4), this quantity has an analytic continuation to $\Omega \cup I$ as
\begin{eqnarray*}
\int_{\mathbf{R}}\! \frac{1}{\lambda -\omega } E [\psi, \phi](\omega )  d\omega 
+ 2\pi \mathrm{i}  E [\psi, \phi](\lambda ) ,
\quad \lambda \in \Omega .
\end{eqnarray*}
Motivated by this observation, define the operator $A(\lambda ) : iX(\Omega ) \to X(\Omega )'$ to be
\begin{equation}
\langle A(\lambda )\psi \,|\, \phi \rangle = \left\{ \begin{array}{ll}
\displaystyle 
     \int_{\mathbf{R}}\! \frac{1}{\lambda -\omega } E [\psi, \phi](\omega ) d\omega 
           + 2\pi \mathrm{i}  E [\psi, \phi](\lambda )  & (\lambda \in \Omega ), \\[0.4cm]
\displaystyle  \lim_{y\to -0} 
  \int_{\mathbf{R}}\! \frac{1}{x + \mathrm{i} y -\omega } E [\psi, \phi](\omega ) d\omega  & (\lambda =x\in I), \\[0.4cm]
\displaystyle  \int_{\mathbf{R}}\! \frac{1}{\lambda -\omega } E [\psi, \phi](\omega ) d\omega
 & (\mathrm{Im}(\lambda ) < 0),
\end{array} \right.
\label{3-12}
\end{equation}
for any $\psi \in iX(\Omega ),\,\, \phi \in X(\Omega )$.
Indeed, we can prove by using (X5) that $A(\lambda )\psi $ is a continuous functional.
Due to Lemma 3.3, $\langle A(\lambda )\psi \,|\, \phi \rangle$ is holomorphic on 
$\{ \mathrm{Im}(\lambda ) < 0\} \cup \Omega \cup I$.
When $\mathrm{Im}(\lambda ) < 0$, we have $\langle A(\lambda )\psi \,|\, \phi \rangle = ( (\lambda - H)^{-1}\psi, \phi )$.
In this sense, the operator $A(\lambda )$ is called the analytic continuation of the resolvent
$(\lambda -H )^{-1}$ as a generalized function.
By using it, we extend the notion of eigenvalues and eigenfunctions.

Recall that the equation for eigenfunctions of $T$ is given by $(id - (\lambda -H)^{-1}K)v = 0$.
Since the analytic continuation of $(\lambda -H )^{-1}$ in $X(\Omega )'$ is $A(\lambda )$,
we make the following definition.
\\[0.2cm]
\textbf{Definition \thedef.} Let $\mathsf{R}(A(\lambda ))$ be the range of $A(\lambda )$.
If the equation
\begin{equation}
(id - A(\lambda )K^\times) \mu = 0
\label{3-15}
\end{equation}
has a nonzero solution $\mu $ in $\mathsf{R}(A(\lambda ))$ for some 
$\lambda \in \Omega \cup I \cup \{ \lambda \, | \, \mathrm{Im} (\lambda ) < 0\}$,
$\lambda $ is called a \textit{generalized eigenvalue} of $T$ and $\mu$ is called a \textit{generalized eigenfunction}
associated with $\lambda $.
A generalized eigenvalue on $\Omega$ is called a \textit{resonance pole}
(the word ``resonance" originates from quantum mechanics \cite{Reed}).

Note that the assumption (X8) is used to define $A(\lambda )K^\times \mu$ for $\mu \in \mathsf{R}(A(\lambda ))$
because the domain of $A(\lambda )$ is $iX(\Omega )$.
Applied by $K^\times$, Eq.(\ref{3-15}) is rewritten as
\begin{equation}
(id - K^\times A(\lambda )) K^\times \mu = 0.
\label{3-16}
\end{equation}
If $K^\times \mu = 0$, Eq.(\ref{3-15}) shows $\mu = 0$.
This means that if $\mu \neq 0$ is a generalized eigenfunction, $K^\times \mu \neq 0$ and $id - K^\times A(\lambda )$
is not injective on $iX(\Omega )$.
Conversely, if $id - K^\times A(\lambda )$ is not injective on $iX(\Omega )$, there is a function $\phi \in iX(\Omega )$
such that $(id - K^\times A(\lambda ))\phi = 0$.
Applying $A(\lambda )$ from the left, we see that $A(\lambda )\phi$ is a generalized eigenfunction.
Hence, $\lambda $ is a generalized eigenvalue if and only if $id - K^\times A(\lambda )$ is not injective on $iX(\Omega )$.
\\[0.2cm]
\textbf{Theorem \thedef.} Let $\lambda $ be a generalized eigenvalue of $T$ and $\mu$ a generalized eigenfunction
associated with $\lambda $. Then the equality
\begin{equation}
T^\times  \mu = \lambda \mu
\end{equation}
holds.
\\[0.2cm]
\textbf{Proof.} At first, let us show $\mathsf{D}(\lambda -H^\times) \supset \mathsf{R}(A(\lambda ))$.
By the operational calculus, we have 
$E[\psi, (\overline{\lambda} -H)\phi ](\omega ) = (\lambda -\omega ) E[\psi, \phi](\omega )$.
When $\lambda \in \Omega $, this gives
\begin{eqnarray*}
\langle A(\lambda )\psi \,|\, (\overline{\lambda } - H) \phi \rangle 
&=& \int_{\mathbf{R}}\! \frac{1}{\lambda -\omega } E [\psi, (\overline{\lambda } - H) \phi](\omega ) d\omega 
           + 2\pi \mathrm{i}  E [\psi, (\overline{\lambda } - H) \phi](\lambda ) \\
&=& \int_{\mathbf{R}}\! E [\psi, \phi](\omega ) d\omega 
           + 2\pi \mathrm{i}  (\lambda -\omega )|_{\omega =\lambda } E [\psi, \phi](\lambda ) \\
&=& \langle \psi \,|\, \phi \rangle,
\end{eqnarray*}
for any $\psi \in X(\Omega )$ and $\phi \in Y$.
It is obvious that $\langle \psi  \,|\, \phi \rangle$ is continuous in $\phi$ with respect to the topology of $X(\Omega )$.
This proves that $\mathsf{D}(\lambda -H^\times) \supset \mathsf{R}(A(\lambda ))$ and 
$(\lambda - H^\times) A(\lambda ) = id : iX(\Omega ) \to iX(\Omega )$.
When $\mu$ is a generalized eigenfunction, $\mu \in \mathsf{D}(\lambda -H^\times)$ because $\mu = A(\lambda )K^\times \mu$.
Then,  Eq.(\ref{3-15}) provides
\begin{eqnarray*}
(\lambda - H^\times)  (id - A(\lambda )K^\times) \mu
=(\lambda - H^\times - K^\times ) \mu = (\lambda -T^\times )\mu =0.
\end{eqnarray*}
The proofs for the cases $\lambda \in I$ and $\mathrm{Im} (\lambda ) < 0$ are done in the same way.
\hfill $\blacksquare$
\\[-0.2cm]

This theorem means that $\lambda $ is indeed an eigenvalue of the dual operator $T^\times$.
In general, the set of generalized eigenvalues is a proper subset of the set of eigenvalues of $T^\times$.
Since the dual space $X(\Omega )'$ is ``too large", 
typically every point on $\Omega $ is an eigenvalue of $T^\times$
(for example, consider the triplet $X \subset L^2(\mathbf{R}) \subset X'$ 
and the multiplication operator $\mathcal{M}$ on $L^2(\mathbf{R})$, where $X$ is the set of entire functions.
Every point on $\mathbf{C}$ is an eigenvalue of the dual operator $\mathcal{M}^\times : X' \to X'$, 
while there are no generalized eigenvalues).
In this sense, generalized eigenvalues are wider concept than eigenvalues of $T$, while narrower concept than
eigenvalues of $T^\times$ (see Prop.3.17 for more details).
In the literature, resonance poles are defined as poles of 
an analytic continuation of a matrix element of the resolvent \cite{Reed}.
Our definition is based on a straightforward extension of the usual eigen-equation 
and it is suitable for systematic studies of resonance poles.

\subsection{Properties of the operator $A(\lambda )$}

Before defining a multiplicity of a generalized eigenvalue, it is convenient to 
investigate properties of the operator $A(\lambda )$.
For $n=1,2,\cdots $ let us define the linear operator $A^{(n)}(\lambda ) : iX(\Omega ) \to X(\Omega )'$ to be
\begin{equation}
\langle A^{(n)}(\lambda )\psi \,|\, \phi \rangle = \left\{ \begin{array}{l}
\displaystyle 
     \int_{\mathbf{R}}\! \frac{1}{(\lambda -\omega )^n} E [\psi, \phi](\omega ) d\omega 
    + 2\pi \mathrm{i}  \frac{(-1)^{n-1}}{(n-1)!} \frac{d^{n-1}}{dz^{n-1}}\Bigl|_{z=\lambda }
               E [\psi, \phi](z), \,\, (\lambda \in \Omega ), \\[0.4cm]
\displaystyle  \lim_{y\to -0} 
  \int_{\mathbf{R}}\! \frac{1}{(x + \mathrm{i} y -\omega )^n} E [\psi, \phi](\omega ) d\omega, \quad (\lambda =x\in I), \\[0.4cm]
\displaystyle  \int_{\mathbf{R}}\! \frac{1}{(\lambda -\omega )^n} E [\psi, \phi](\omega ) d\omega,
\quad (\mathrm{Im}(\lambda ) < 0).
\end{array} \right.
\end{equation}
It is easy to show by integration by parts that $\langle A^{(n)}(\lambda )\psi \,|\, \phi \rangle$ is an
analytic continuation of $( (\lambda -H)^{-n} \psi , \phi )$ from the lower half plane to $\Omega $.
$A^{(1)}(\lambda )$ is also denoted by $A(\lambda )$ as before.
The next proposition will be often used to calculate the generalized resolvent and projections.
\\[0.2cm]
\textbf{Proposition \thedef.} For any integers $j\geq n \geq 0$. the operator $A^{(j)}(\lambda )$ satisfies
\\[0.2cm]
(i) $(\lambda -H^\times)^nA^{(j)}(\lambda ) = A^{(j-n)}(\lambda )$, where $A^{(0)}(\lambda ) := id$.
\\[0.2cm]
(ii) $A^{(j)}(\lambda )(\lambda -H^\times)^n|_{iX(\Omega )\cap \mathsf{D}(A^{(j)}(\lambda )(\lambda -H^\times)^n)}
 = A^{(j-n)}(\lambda )|_{iX(\Omega )\cap \mathsf{D}(A^{(j)}(\lambda )(\lambda -H^\times)^n)}$.
\\
In particular, $A(\lambda )(\lambda -H^\times)\mu = \mu$ when $(\lambda -H^\times)\mu \in iX(\Omega )$.
\\[0.2cm]
(iii) $\displaystyle \frac{d^j}{d\lambda ^j} \langle A(\lambda )\psi \,|\,  \phi \rangle  
 = (-1)^j j! \langle A^{(j+1)}(\lambda )\psi \,|\, \phi \rangle,\,\, j=0,1,\cdots $.
\\[0.2cm]
(iv) For each $\psi \in X(\Omega )$, $A(\lambda ) \psi $ is expanded as
\begin{equation}
 A(\lambda ) \psi 
 = \sum^\infty_{j=0} (\lambda _0 - \lambda )^j A^{(j+1)}(\lambda_0 ) \psi ,
\label{3-21}
\end{equation}
where the right hand side converges with respect to the strong dual topology.
\\[0.2cm]
\textbf{Proof.} (i) Let us show $(\lambda -H^\times) A^{(j)}(\lambda ) = A^{(j-1)}(\lambda )$.
We have to prove that $\mathsf{D}(\lambda -H^\times) \supset \mathsf{R}(A^{(j)}(\lambda ))$.
For this purpose, put $\mu_\lambda (y) = \langle A^{(j)}(\lambda )\psi \,|\, (\overline{\lambda } - H) y \rangle$ for 
$\psi\in X(\Omega )$ and $y\in Y$.
It is sufficient to show that the mapping $y\mapsto \mu_\lambda (y)$ from $Y$ into $\mathbf{C}$ is continuous
with respect to the topology on $X(\Omega )$.
Suppose that $\mathrm{Im}(\lambda ) > 0$.
By the operational calculus, we obtain
\begin{eqnarray}
\mu_\lambda (y) &=& \int_{\mathbf{R}}\! \frac{1}{(\lambda -\omega )^j} E[\psi, (\overline{\lambda }-H)y](\omega )d\omega 
+2\pi\mathrm{i} \frac{(-1)^{j-1}}{(j-1)!}\frac{d^{j-1}}{dz^{j-1}}\Bigl|_{z=\lambda }E[\psi,(\overline{\lambda }-H)y](z)\nonumber \\
&=&  \int_{\mathbf{R}}\! \frac{\lambda -\omega }{(\lambda -\omega )^j} E[\psi, y](\omega )d\omega 
+2\pi\mathrm{i} \frac{(-1)^{j-1}}{(j-1)!}\frac{d^{j-1}}{dz^{j-1}}\Bigl|_{z=\lambda }(\lambda -z)E[\psi, y](z)\nonumber \\
&=& ((\lambda -H)^{1-j}\psi, y) 
+ 2\pi\mathrm{i} \frac{(-1)^{j-2}}{(j-2)!}\frac{d^{j-2}}{dz^{j-2}}\Bigl|_{z=\lambda }E[\psi, y](z).
\label{3-1}
\end{eqnarray}
Since $E[\psi, y](z)$ is continuous in $y\in X(\Omega )$ (the assumption (X5)) and $E[\psi, y](z)$ 
is holomorphic in $z$, for any $\varepsilon >0$, there exists a neighborhood $U_1$ of zero in $X(\Omega )$ such that 
$|(d^{j-2}/dz^{j-2})E[\psi, y] (z)| < \varepsilon $ at $z=\lambda $ for $y \in U_1 \cap Y$.
Let $U_2$ be a neighborhood of zero in $\mathcal{H}$ such that $|| y ||_{\mathcal{H}} < \varepsilon $ for $y \in U_2$.
Since the topology on $X(\Omega )$ is stronger than that on $\mathcal{H}$, 
$U_2 \cap X(\Omega )$ is a neighborhood of zero in $X(\Omega )$.
If $y\in  U_1 \cap U_2 \cap Y$, we obtain
\begin{eqnarray*}
| \mu_\lambda (y)| \leq || (\lambda -H)^{1-j}\psi || \varepsilon  + 2\pi \mathrm{i} \frac{(-1)^{j-2}}{(j-2)!} \varepsilon .
\end{eqnarray*}
Note that $(\lambda -H)^{1-j}$ is bounded when $\lambda \notin \mathbf{R}$ and $1-j\leq 0$ because $H$ is selfadjoint.
This proves that $\mu_\lambda$ is continuous, so that $\mu_\lambda = (\lambda -H^\times)A^{(j)}(\lambda )\psi \in X(\Omega )'$.
The proof of the continuity for the case $\mathrm{Im}(\lambda ) < 0$ is done in the same way.
When $\lambda \in I$, there exists a sequence $\{ \lambda _j\}^\infty_{j=1}$ in the lower half plane such that
$\mu_\lambda (y)  = \lim_{j\to \infty}\mu_{\lambda_j} (y)$.
Since $X(\Omega )$ is barreled, Banach-Steinhaus theorem is applicable to conclude that
the limit $\mu_\lambda$ of continuous linear mappings is also continuous. 
This proves $\mathsf{D}(\lambda -H^\times) \supset \mathsf{R}(A^{(j)}(\lambda ))$
and $(\lambda -H^\times) A^{(j)}(\lambda )$ is well defined 
for any $\lambda \in \{ \mathrm{Im} (\lambda ) < 0\} \cup I \cup \Omega $.
Then, the above calculation immediately shows that $(\lambda -H^\times) A^{(j)}(\lambda ) = A^{(j-1)}(\lambda )$.
By the induction, we obtain (i).

(ii) is also proved by the operational calculus as above, and
(iii) is easily obtained by induction.

For (iv), since $\langle A(\lambda )\psi \,|\, \phi \rangle$ is holomorphic,
it is expanded in a Taylor series as
\begin{eqnarray}
\langle A(\lambda )\psi \,|\, \phi \rangle &=& \sum^\infty_{j=0} \frac{1}{j!}\frac{d^j}{d\lambda ^j}
        \Bigl|_{\lambda =\lambda _0} \langle A(\lambda )\psi \,|\, \phi \rangle (\lambda -\lambda _0)^j \nonumber \\ 
&=& \sum^\infty_{j=0} (\lambda_0 - \lambda )^j \langle A^{(j+1)}(\lambda _0) \psi \,|\, \phi \rangle,
\label{3-22}
\end{eqnarray}
for each $\phi ,\psi \in X(\Omega )$.
This means that the functional $A(\lambda )\psi $ is weakly holomorphic in $\lambda $.
Then, $A(\lambda )\psi $ turns out to be strongly holomorphic and expanded as Eq.(\ref{3-21})
by Thm.A.3(iii) in Appendix, in which basic facts on $X(\Omega )'$-valued holomorphic functions are given.
 \hfill $\blacksquare$
\\[-0.2cm]

Unfortunately, the operator $A(\lambda ) : iX(\Omega ) \to X(\Omega )'$ is not continuous
if $iX(\Omega )$ is equipped with the relative topology from $X(\Omega )'$.
Even if $\langle \psi | \to 0$ in $iX(\Omega )\subset X(\Omega )'$, the value $E[\psi, \phi](\lambda )$ does not tend 
to zero in general because the topology on $X(\Omega )'$ is weaker than that on $X(\Omega )$ .
However, $A(\lambda )$ proves to be continuous if $iX(\Omega )$ is equipped with the topology induced from $X(\Omega )$
by the canonical inclusion.
\\[0.2cm]
\textbf{Proposition \thedef.} 
$A(\lambda ) \circ i  : X(\Omega ) \to X(\Omega )'$ is continuous if $X(\Omega )'$ is equipped with the weak dual topology.
\\[0.2cm]
\textbf{Proof.} 
Suppose $\lambda \in \Omega $ and fix $\phi \in X(\Omega )$.
Because of the assumption (X5), for any $\varepsilon >0$, there exists
a neighborhood $U_1$ of zero in $X(\Omega )$ such that $|E[\psi, \phi] (\lambda )| < \varepsilon $ for $\psi \in U_1$.
Let $U_2$ be a neighborhood of zero in $\mathcal{H}$ such that $|| \psi ||_{\mathcal{H}} < \varepsilon $ for $\psi \in U_2$.
Since the topology on $X(\Omega )$ is stronger than that on $\mathcal{H}$, 
$U_2 \cap X(\Omega )$ is a neighborhood of zero in $X(\Omega )$.
If $\psi\in U: = U_1 \cap U_2$,
\begin{eqnarray*}
| \langle A(\lambda ) \psi \,|\, \phi \rangle |
 &\leq & || (\lambda -H)^{-1} ||_{\mathcal{H}}\cdot || \phi ||_{\mathcal{H}}\cdot || \psi ||_{\mathcal{H}}
 + 2\pi \,|E[\psi, \phi] (\lambda )| \\
&=& \left( || (\lambda -H)^{-1} ||_{\mathcal{H}}\cdot || \phi ||_{\mathcal{H}} + 2\pi \right) \varepsilon .
\end{eqnarray*}
This proves that $A(\lambda ) \circ i$ is continuous in the weak dual topology.
The proof for the case $\mathrm{Im}(\lambda ) < 0$ is done in a similar manner.
When $\lambda \in I$, there exists a sequence $\{ \lambda _j\}^\infty_{j=1}$ in the lower half plane such that
$A(\lambda ) \circ i = \lim_{j\to \infty}A(\lambda _j) \circ i$.
Since $X(\Omega )$ is barreled, Banach-Steinhaus theorem is applicable to conclude that
the limit $A(\lambda ) \circ i $ of continuous linear mappings is also continuous.
\hfill $\blacksquare$
\\[-0.2cm]

Now we are in a position to define an algebraic multiplicity and a generalized eigenspace of generalized eigenvalues.
Usually, an eigenspace is defined as a set of solutions of the equation $(\lambda -T)^nv = 0$.
For example, when $n=2$, we rewrite it as
\begin{eqnarray*}
(\lambda -H -K)(\lambda -H -K)v
 = (\lambda - H )^2(id - (\lambda -H )^{-2} K(\lambda - H)) \circ (id - (\lambda - H)^{-1}K)v = 0.
\end{eqnarray*}
Dividing by $(\lambda -H)^{2}$ yields
\begin{eqnarray*}
(id - (\lambda -H )^{-2} K(\lambda - H)) \circ (id - (\lambda - H)^{-1}K)v = 0.
\end{eqnarray*}
Since the analytic continuation of $(\lambda -H)^{-n}$ in $X(\Omega )'$ is $A^{(n)}(\lambda )$, we consider the equation
\begin{eqnarray*}
(id - A^{(2)}(\lambda )K^\times (\lambda -H^\times ))
\circ (id- A(\lambda )K^\times) \,\mu = 0.
\end{eqnarray*}
Motivated by this observation, we define the operator $B^{(n)}(\lambda ) : \mathsf{D}(B^{(n)}(\lambda ))
\subset X(\Omega )' \to X(\Omega )'$ to be
\begin{equation}
B^{(n)}(\lambda )=id - A^{(n)}(\lambda )K^\times (\lambda -H^\times)^{n-1}.
\label{3-20}
\end{equation}
Then, the above equation is rewritten as $B^{(2)}(\lambda )B^{(1)}(\lambda )\mu = 0$.
The domain of $B^{(n)}(\lambda )$ is the domain of $A^{(n)}(\lambda )K^\times (\lambda -H^\times)^{n-1}$.
The following equality is easily proved.
\begin{eqnarray}
(\lambda -H^\times)^k B^{(j)}(\lambda )
 = B^{(j-k)}(\lambda ) (\lambda -H^\times)^k|_{\mathsf{D}(B^{(j)}(\lambda ))}, \quad j > k.
\label{3-25} 
\end{eqnarray}
\textbf{Definition \thedef.} 
Let $\lambda $ be a generalized eigenvalue of the operator $T$.
The generalized eigenspace of $\lambda $ is defined by
\begin{equation}
V_\lambda = \bigcup_{m\geq 1}\mathrm{Ker}\, B^{(m)}(\lambda ) \circ B^{(m-1)}(\lambda ) \circ \cdots \circ B^{(1)}(\lambda ).
\end{equation}
We call $\mathrm{dim} V_\lambda $ the algebraic multiplicity of the generalized eigenvalue $\lambda $.
\\[0.2cm]
\textbf{Theorem \thedef.} For any $\mu \in V_\lambda $, there exists an integer $M$ such that 
$(\lambda - T^\times) ^M \mu = 0$.
\\[0.2cm]
\textbf{Proof.} Suppose that $B^{(M)}(\lambda ) \circ \cdots \circ B^{(1)}(\lambda ) \mu = 0$.
Put $\xi = B^{(M-1)}(\lambda ) \circ \cdots \circ B^{(1)}(\lambda ) \mu $.
Eq.(\ref{3-25}) shows
\begin{eqnarray*}
0 &=& (\lambda -H^\times)^{M-1} B^{(M)}(\lambda ) \xi \\
&=& B^{(1)}(\lambda ) (\lambda -H^\times)^{M-1}\xi = (id - A(\lambda )K^\times) (\lambda -H^\times)^{M-1}\xi.
\end{eqnarray*}
Since $\mathsf{D}(\lambda -H^\times) \supset \mathsf{R}(A(\lambda ))$,
it turns out that $(\lambda -H^\times)^{M-1}\xi \in \mathsf{D}(\lambda -H^\times)$.
Then, we obtain
\begin{eqnarray*}
0 &=& (\lambda -H^\times)(id - A(\lambda )K^\times) (\lambda -H^\times)^{M-1}\xi \\
 &=& (\lambda -H^\times - K^\times )(\lambda -H^\times)^{M-1}\xi = (\lambda -T^\times )(\lambda -H^\times)^{M-1}\xi.
\end{eqnarray*}
By induction, we obtain $(\lambda - T^\times) ^M \mu = 0$. \hfill $\blacksquare$
\\[-0.2cm]

In general, the space $V_\lambda $ is a proper subspace of the usual eigenspace 
$\bigcup_{m\geq 1} \mathrm{Ker}\, (\lambda - T^\times)^m$ of $T^\times$.
Typically $\bigcup_{m\geq 1} \mathrm{Ker}\, (\lambda - T^\times)^m$ becomes of
infinite dimensional because the dual space $X(\Omega )'$ is ``too large",
however, $V_\lambda $ is a finite dimensional space in many cases.

\subsection{Generalized resolvents}

In this subsection, we define a generalized resolvent.
As the usual theory, it will be used to construct projections and semigroups.
Let  $R_\lambda = (\lambda -T)^{-1}$ be the resolvent of $T$ as an operator on $\mathcal{H}$.
A simple calculation shows 
\begin{equation}
R_\lambda \psi = (\lambda -H)^{-1} \left( id - K(\lambda -H )^{-1} \right)^{-1}\psi.
\label{3-29}
\end{equation}
Since the analytic continuation of $(\lambda -H)^{-1}$ in the dual space is $A(\lambda )$, we make the following definition.
In what follows, put $\hat{\Omega } = \Omega \cup I \cup \{ \lambda \, | \, \mathrm{Im} (\lambda ) < 0 \}$.
\\[0.2cm]
\textbf{Definition \thedef.}
If the inverse $(id - K^\times A(\lambda ))^{-1}$ exists,
define the generalized resolvent $\mathcal{R}_\lambda : iX(\Omega ) \to X(\Omega )'$ to be
\begin{equation}
\mathcal{R}_\lambda =A(\lambda ) \circ (id - K^\times A(\lambda ))^{-1}
      = (id - A(\lambda )K^\times)^{-1} \circ A(\lambda ), \quad \lambda 
 \in \hat{\Omega }.
\label{3-30}
\end{equation}
The second equality follows from $(id - A(\lambda )K^\times) A(\lambda ) = A(\lambda ) (id - K^\times A(\lambda ))$.
Recall that $id - K^\times A(\lambda )$ is injective on $iX(\Omega )$ if and only if 
$id - A(\lambda )K^\times$ is injective on $\mathsf{R}(A(\lambda ))$.
\\[-0.2cm]

Since $A(\lambda )$ is not continuous, $\mathcal{R}_\lambda $ is not a continuous operator in general.
However, it is natural to ask whether $\mathcal{R}_\lambda \circ i : X(\Omega ) \to X(\Omega )'$ is continuous or not
because $A(\lambda ) \circ i$ is continuous.
\\[0.2cm]
\textbf{Definition \thedef.}
The generalized resolvent set $\hat{\varrho}(T)$ is defined to be the set of points $\lambda \in \hat{\Omega }$ satisfying following:
there is a neighborhood $V_\lambda \subset \hat{\Omega }$ of $\lambda $ such that 
for any $\lambda ' \in V_\lambda $, $\mathcal{R}_{\lambda '} \circ i$ is a densely 
defined continuous operator from $X(\Omega )$ into $X(\Omega )'$,
where $X(\Omega )'$ is equipped with the weak dual topology,
and the set $\{ \mathcal{R}_{\lambda '} \circ i(\psi)\}_{\lambda ' \in V_\lambda }$ is bounded in $X(\Omega )'$ for each $\psi \in X(\Omega )$.
The set $\hat{\sigma }(T):= \hat{\Omega }\backslash \hat{\varrho}(T)$ is called the \textit{generalized spectrum} of $T$.
The \textit{generalized point spectrum} $\hat{\sigma }_p(T)$ is the set of points $\lambda \in \hat{\sigma }(T)$ at which
$id - K^\times A(\lambda )$ is not injective (this is the set of generalized eigenvalues).
The \textit{generalized residual spectrum} $\hat{\sigma }_r(T)$ is the set of points $\lambda \in \hat{\sigma }(T)$ 
such that the domain of $\mathcal{R}_\lambda \circ i$ is not dense in $X(\Omega )$.
The \textit{generalized continuous spectrum} is defined to be 
$\hat{\sigma }_c(T) = \hat{\sigma }(T)\backslash (\hat{\sigma }_p(T)\cup \hat{\sigma }_r(T))$.
\\[-0.2cm]

By the definition, $\hat{\varrho}(T)$ is an open set.
To require the existence of the neighborhood $V_\lambda $ in the above definition is introduced by Waelbroeck \cite{Wae} 
(see also Maeda \cite{Mae}) for the spectral theory on locally convex spaces.
If $\hat{\varrho}(T)$ were simply defined to be the set of points such that $\mathcal{R}_\lambda \circ i$ is a densely defined
continuous operator as in the Banach space theory, $\hat{\varrho}(T)$ is not an open set in general.
If $X(\Omega )$ is a Banach space and the operator $i^{-1}K^\times A(\lambda )i$ is continuous on $X(\Omega )$ for each
$\lambda \in \hat{\Omega }$, we can show that $\lambda \in \hat{\varrho}(T)$ if and only if 
$id - i^{-1}K^\times A(\lambda )i$ has a continuous inverse on $X(\Omega )$ (Prop.3.18).
\\[0.2cm]
\textbf{Theorem \thedef.}
\\
(i) For each $\psi \in X(\Omega )$, $\mathcal{R}_\lambda i\psi $ is an $X(\Omega )'$-valued holomorphic function
in $\lambda \in \hat{\varrho} (T)$.
\\
(ii) Suppose $\mathrm{Im} (\lambda ) < 0$ and $\lambda \in \hat{\varrho} (T) \cap \varrho (T)$,
where $\varrho (T)$ is the resolvent set of $T$ in $\mathcal{H}$-sense.
Then, $\langle \mathcal{R}_\lambda \psi \,|\, \phi \rangle 
 = ( (\lambda -T)^{-1} \psi ,\phi)$ for any $\psi, \phi \in X(\Omega )$.
\\[-0.2cm]

This theorem means that $\langle \mathcal{R}_\lambda \psi \,|\, \phi \rangle$
is an analytic continuation of $( (\lambda -T)^{-1} \psi , \phi)$ from 
the lower half plane to $\hat{\varrho}(T)$ through the interval $I$.
We always suppose that the domain of $\mathcal{R}_\lambda \circ i$ is continuously extended to
the whole $X(\Omega )$ when $\lambda \notin \hat{\sigma} (T)$.
The significant point to be emphasized is that to prove the \textit{strong} holomorphy of $\mathcal{R}_\lambda \circ i(\psi )$,
it is sufficient to assume that $\mathcal{R}_\lambda \circ i : X(\Omega ) \to X(\Omega )'$
is continuous in the \textit{weak} dual topology on $X(\Omega )'$.
\\[0.2cm]
\textbf{Proof of Thm.3.12.}
Since $\hat{\varrho}(T)$ is open, 
when $\lambda \in \hat{\varrho}(T)$, $\mathcal{R}_{\lambda +h}$ exists for sufficiently small $h\in \mathbf{C}$.
Put $\psi _\lambda = i^{-1}(id - K^\times A(\lambda ))^{-1}i (\psi)$ for $\psi \in X(\Omega )$.
It is easy to verify the equality
\begin{eqnarray*}
\mathcal{R}_{\lambda +h}i(\psi) - \mathcal{R}_\lambda i(\psi)
= (A(\lambda +h) - A(\lambda ))i(\psi_\lambda ) + \mathcal{R}_{\lambda +h}i \circ 
i^{-1}K^\times (A(\lambda +h) - A(\lambda ))i(\psi_\lambda ).
\end{eqnarray*}
Let us show that $i^{-1}K^\times A(\lambda )i (\psi )\in X(\Omega )$ is holomorphic in $\lambda $.
For any $\psi, \phi \in X(\Omega )$, we obtain
\begin{eqnarray*}
& & \langle \phi \,|\, i^{-1}K^\times A(\lambda )i\psi \rangle
 = ( \phi , i^{-1}K^\times A(\lambda )i\psi )
 = \overline{( i^{-1}K^\times A(\lambda )i\psi, \phi) } \\
&=& \langle \overline{K^\times A(\lambda )i \psi \,|\, \phi } \rangle 
 = \langle \overline{A(\lambda )i \psi \,|\, K^*\phi } \rangle .
\end{eqnarray*}
From the definition of $A(\lambda )$, it follows that $\langle \phi \,|\, i^{-1}K^\times A(\lambda )i\psi \rangle$
is holomorphic in $\overline{\lambda }$.
Since $X(\Omega )$ is dense in $X(\Omega )'$, $\langle \mu \,|\, i^{-1}K^\times A(\lambda )i\psi \rangle$ is 
holomorphic in $\overline{\lambda }$ for any $\mu \in X(\Omega )'$ by Montel's theorem.
This means that $i^{-1}K^\times A(\lambda )i \psi$ is weakly holomorphic.
Since $X(\Omega )$ is a quasi-complete locally convex space, any weakly holomorphic function is 
holomorphic with respect to the original topology (see Rudin \cite{Rud}).
This proves that $i^{-1}K^\times A(\lambda )i \psi$ is holomorphic in $\lambda $
(note that the weak holomorphy in $\overline{\lambda }$ implies the strong holomorphy in $\lambda $
because functionals in $X(\Omega )'$ are \textit{anti}-linear).

Next, the definition of $\hat{\varrho}(T)$ implies that the family $\{ \mathcal{R}_\mu \circ i\}_{\mu \in V_\lambda }$
of continuous operators is bounded in the pointwise convergence topology.
Due to Banach-Steinhaus theorem (Thm.33.1 of \cite{Tre}), the family is equicontinuous.
This fact and the holomorphy of $A(\lambda )$ and $i^{-1}K^\times A(\lambda )i (\psi )$ prove that
$\mathcal{R}_{\lambda +h}i(\psi)$ converges to $\mathcal{R}_{\lambda }i(\psi)$ as $h\to 0$ with respect to the 
weak dual topology.
In particular, we obtain
\begin{equation}
\lim_{h\to 0} \frac{\mathcal{R}_{\lambda +h}i - \mathcal{R}_\lambda i}{h}(\psi)
= \frac{dA}{d\lambda }(\lambda )i (\psi_\lambda ) + \mathcal{R}_{\lambda }i \circ 
\frac{d}{d\lambda } (i^{-1}K^\times A(\lambda ) i )(\psi _\lambda ),
\label{thm3-12}
\end{equation}
which proves that
$\mathcal{R}_\lambda i(\psi)$ is holomorphic in $\lambda $ with respect to the weak dual topology on $X(\Omega )'$.
Since $X(\Omega )$ is barreled, the weak dual holomorphy implies the strong dual holomorphy (Thm.A.3 (iii)). 

Let us prove (ii). Suppose $\mathrm{Im} (\lambda ) < 0$.
Note that $\mathcal{R}_\lambda \circ i$ is written as 
$\mathcal{R}_\lambda  \circ i = A(\lambda ) \circ (id - i^{-1}K^\times A(\lambda ) i)^{-1}$.
We can show the equality
\begin{equation}
(id - i^{-1}K^\times A(\lambda ) i) f = (id -K(\lambda -H)^{-1}) f \in X(\Omega ).
\end{equation}
Indeed, for any $f, \psi \in X(\Omega )$, we obtain
\begin{eqnarray*}
\langle (i - K^\times A(\lambda ) i)f \,|\, \psi \rangle
&=& \langle if \,|\, \psi \rangle - \langle A(\lambda )if \,|\, K^* \psi \rangle \\
&=& \langle if \,|\, \psi \rangle - \langle i \circ (\lambda -H)^{-1}f \,|\, K^* \psi \rangle \\
&=& (f, \psi) - (K(\lambda -H)^{-1} f, \psi) = ((id -K(\lambda -H)^{-1}) f, \psi ).
\end{eqnarray*}
Thus, $\mathcal{R}_\lambda $ satisfies for $\phi = (id - i^{-1}K^\times A(\lambda ) i) f$ that
\begin{eqnarray*}
\mathcal{R}_\lambda i \phi  
&=& A(\lambda )i \circ (id - i^{-1}K^\times A(\lambda ) i)^{-1} \phi \\
&=& i (\lambda -H)^{-1} \circ (id - K(\lambda -H)^{-1})^{-1} \phi = i(\lambda -T)^{-1}\phi.
\end{eqnarray*}
Since $\lambda \in \hat{\rho} (T)$, $(id - i^{-1}K^\times A(\lambda ) i) X(\Omega )$ is dense in $X(\Omega )$
and $\mathcal{R}_\lambda i : X(\Omega ) \to X(\Omega )'$ is continuous.
Since $\lambda \in \rho (T)$, $i(\lambda -T)^{-1} : \mathcal{H} \to X(\Omega )'$ is continuous.
Therefore, taking the limit proves that $\mathcal{R}_\lambda i \phi = i(\lambda -T)^{-1}\phi$ holds for 
any $\phi \in X(\Omega )$.
 \hfill $\blacksquare$
\\[0.2cm]
\textbf{Remark.} Even when $\lambda $ is in the continuous spectrum of $T$,
Thm.3.12 (ii) holds as long as $(\lambda -T)^{-1}$ exists and $i\circ (\lambda -T)^{-1} : \mathcal{H} \to X(\Omega )'$
is continuous.
In general, the continuous spectrum of $T$ is not included in the generalized spectrum because the topology of $X(\Omega )'$
is weaker than that of $\mathcal{H}$.
\\[0.2cm]
\textbf{Proposition \thedef.}
The generalized resolvent satisfies
\\
(i) $(\lambda - T^\times ) \circ \mathcal{R}_\lambda = id|_{iX(\Omega )}$
\\
(ii) If $\mu \in X(\Omega )'$ satisfies $(\lambda -T^\times) \mu \in iX(\Omega )$, then
$\mathcal{R}_\lambda \circ (\lambda  - T^\times )\mu = \mu$.
\\
(iii)  $T^\times \circ \mathcal{R}_\lambda|_{iY} = \mathcal{R}_\lambda \circ T^\times |_{iY}$.
\\[0.2cm]
\textbf{Proof.}
Prop.3.6 (i) gives $id = (\lambda -H^\times) A(\lambda ) = (\lambda -T^\times + K^\times) A(\lambda ).$
This proves
\begin{eqnarray*}
& & (\lambda - T^\times ) \circ A(\lambda ) = id -  K^\times  A(\lambda ) \\
& \Rightarrow & (\lambda - T^\times ) \circ A(\lambda ) \circ 
    (id -  K^\times  A(\lambda ))^{-1} = (\lambda - T^\times ) \circ
\mathcal{R}_\lambda  = id.
\end{eqnarray*}
Next, when $(\lambda -T^\times) \mu \in iX(\Omega )$, $A(\lambda )(\lambda -T^\times) \mu $ is well defined and
Prop.3.6 (ii) gives 
\begin{eqnarray*}
A(\lambda )(\lambda -T^\times) \mu 
 = A(\lambda )(\lambda -H^\times - K^\times)\mu = (id - A(\lambda ) K^\times )\mu .
\end{eqnarray*}
This proves $\mu = (id - A(\lambda ) K^\times )^{-1}A(\lambda )(\lambda -T^\times )\mu
 = \mathcal{R}_\lambda (\lambda -T^\times )\mu$.
Finally, note that $(\lambda -T^\times) iY = i(\lambda -T)Y \subset iX(\Omega )$ because of the assumptions (X6), (X7).
Thus part (iii) of the proposition immediately follows from (i), (ii). \hfill $\blacksquare$

\subsection{Generalized projections}

Let $\Sigma \subset \hat{\sigma}(T)$ be a bounded subset of the generalized spectrum,
which is separated from the rest of the spectrum by a simple closed curve 
$\gamma \subset \Omega \cup I \cup \{ \lambda \, | \, \mathrm{Im} (\lambda ) < 0\}$.
Define the operator $\Pi_\Sigma : iX(\Omega ) \to X(\Omega )'$ to be
\begin{equation}
\Pi_\Sigma \phi= \frac{1}{2\pi \mathrm{i} } \int_{\gamma }\! \mathcal{R}_\lambda \phi \,d\lambda,
\quad \phi \in iX(\Omega ),
\label{3-33}
\end{equation}
where the integral is defined as the Pettis integral.
Since $X(\Omega )$ is assumed to be barreled by (X3), $X(\Omega )'$ is quasi-complete and satisfies the
convex envelope property (see Appendix A).
Since $\mathcal{R}_\lambda \phi $ is strongly holomorphic in $\lambda $ (Thm.3.12),
the Pettis integral of $\mathcal{R}_\lambda \phi$ exists by Thm.A.1.
See Appendix A for the definition and the existence theorem of Pettis integrals.
Since $\mathcal{R}_\lambda \circ i : X(\Omega ) \to X(\Omega )'$ is continuous, Thm.A.1 (ii) proves that 
$\Pi_\Sigma \circ i$ is a continuous operator from $X(\Omega )$ into $X(\Omega )'$ equipped with the weak dual topology.
Note that the equality
\begin{equation}
T^\times \int_{\gamma }\! \mathcal{R}_\lambda \phi \,d\lambda = \int_{\gamma }\! T^\times \mathcal{R}_\lambda \phi \,d\lambda,
\label{3-23b}
\end{equation}
holds.
To see this, it is sufficient to show that the set 
$\{ \langle T^\times \mathcal{R}_\lambda \phi \,|\, \psi \rangle\}_{\lambda \in \gamma }$ is bounded
for each $\psi \in X(\Omega )$ due to Thm.A.1 (iii).
Prop.3.13 (i) yields $T^\times \mathcal{R}_\lambda \phi= \lambda \mathcal{R}_\lambda \phi - \phi$.
Since $\lambda \mathcal{R}_\lambda $ is holomorphic and $\gamma $ is compact, 
$\{ \langle T^\times \mathcal{R}_\lambda \phi \,|\, \psi \rangle\}_{\lambda \in \gamma }$ is bounded so that Eq.(\ref{3-23b})
holds.

Although $\Pi_\Sigma \circ \Pi_\Sigma $ is not defined, we call $\Pi_\Sigma$ the \textit{generalized Riesz projection}
for $\Sigma $ because of the next proposition.
\\[0.2cm]
\textbf{Proposition \thedef.}
$\Pi_\Sigma (iX(\Omega )) \cap (id - \Pi_\Sigma )(iX(\Omega )) = \{ 0\}$ and the direct sum satisfies
\begin{equation}
iX(\Omega ) \subset \Pi_\Sigma (iX(\Omega )) \oplus (id - \Pi_\Sigma )(iX(\Omega )) \subset X(\Omega )'.
\end{equation}
In particular, for any $\phi \in X(\Omega )$, there exist $\mu_1, \mu_2$ such that $\phi$ is uniquely decomposed as
\begin{equation}
i(\phi) = \langle \phi | = \mu_1 + \mu_2, \quad \mu_1 \in \Pi_\Sigma (iX(\Omega )),\,\,
 \mu_2 \in (id - \Pi_\Sigma )(iX(\Omega )).
\label{3-34}
\end{equation}
\textbf{Proof.} We simply denote $\langle \phi |$ as $\phi$.
It is sufficient to show that $\Pi_\Sigma (iX(\Omega )) \cap (id - \Pi_\Sigma )(iX(\Omega )) = \{ 0\}$.
Suppose that there exist $\phi , \psi \in iX(\Omega )$ such that $\Pi_\Sigma  \phi = \psi - \Pi_\Sigma \psi$.
Since $\Pi_\Sigma (\phi + \psi ) = \psi \in iX(\Omega )$,
we can again apply the projection to the both sides as $\Pi_\Sigma \circ \Pi_\Sigma (\phi + \psi) = \Pi_\Sigma \psi $.
Let $\gamma '$ be a closed curve which is slightly larger than $\gamma $. Then,
\begin{eqnarray*}
\Pi_\Sigma \circ \Pi_\Sigma (\phi+ \psi)
&=& \left( \frac{1}{2\pi \mathrm{i} }\right) ^2 
  \int_{\gamma '}\! \mathcal{R}_{\lambda '}\left( \int_{\gamma }\!\mathcal{R}_\lambda (\phi +\psi)d\lambda\right) d\lambda '\\
&=& \left( \frac{1}{2\pi \mathrm{i} }\right) ^2 
       \int_{\gamma '}\! \mathcal{R}_{\lambda '} \left(  \int_{\gamma }\! 
            \frac{(\lambda -\lambda ')+(\lambda '-T^\times)}{\lambda -\lambda '}\mathcal{R}_{\lambda }(\phi + \psi) 
              d\lambda \right) d\lambda ' \\
& & \quad - \left( \frac{1}{2\pi \mathrm{i} }\right) ^2 
       \int_{\gamma '}\! \mathcal{R}_{\lambda '}\left( \int_{\gamma }\! \frac{\lambda '-T^\times}{\lambda -\lambda '}
           \mathcal{R}_\lambda (\phi + \psi) d\lambda \right) d\lambda '. 
\end{eqnarray*}
Eq.(\ref{3-23b}) shows
\begin{eqnarray*}
\Pi_\Sigma \circ \Pi_\Sigma (\phi+ \psi) & = & \left( \frac{1}{2\pi \mathrm{i} }\right) ^2 
       \int_{\gamma '}\! \mathcal{R}_{\lambda '} \left(  \int_{\gamma }\! 
            \frac{\lambda -T^\times}{\lambda -\lambda '}\mathcal{R}_{\lambda }(\phi + \psi) d\lambda \right) d\lambda ' \\
& & \quad - \left( \frac{1}{2\pi \mathrm{i} }\right) ^2 
       \int_{\gamma '}\! \mathcal{R}_{\lambda '}\circ (\lambda '-T^\times)
          \left( \int_{\gamma }\! \frac{ \mathcal{R}_\lambda}{\lambda -\lambda '}(\phi + \psi) d\lambda \right) d\lambda '.
\end{eqnarray*}
Prop.3.13 shows
\begin{eqnarray*}
&=& \left( \frac{1}{2\pi \mathrm{i} }\right) ^2  \int_{\gamma '}\! \mathcal{R}_{\lambda '} \left(  \int_{\gamma }\! 
            \frac{\phi + \psi}{\lambda -\lambda '}d\lambda \right) d\lambda ' 
- \left( \frac{1}{2\pi \mathrm{i} }\right) ^2  \int_{\gamma '}\! 
          \left( \int_{\gamma }\! \frac{ \mathcal{R}_\lambda}{\lambda -\lambda '}(\phi + \psi) d\lambda \right) d\lambda '\\
&=& 0 - \left( \frac{1}{2\pi \mathrm{i} }\right) ^2  \int_{\gamma }\! \mathcal{R}_\lambda (\phi + \psi) \cdot
              \int_{\gamma '}\! \frac{1}{\lambda -\lambda '}d\lambda' \cdot d\lambda\\
&=& \frac{1}{2\pi \mathrm{i} } \int_{\gamma }\! \mathcal{R}_\lambda (\phi + \psi) d\lambda 
=  \Pi_\Sigma (\phi+ \psi).
\end{eqnarray*}
This proves that $\Pi_\Sigma \phi = 0$. \hfill $\blacksquare$
\\[0.2cm]
The above proof also shows that as long as $\Pi_\Sigma \phi \in iX(\Omega )$,
$\Pi_\Sigma \circ \Pi_\Sigma$ is defined and $\Pi_\Sigma  \circ \Pi_\Sigma \phi = \Pi_\Sigma  \phi$.
\\[0.2cm]
\textbf{Proposition \thedef.} $\Pi_\Sigma |_{iY}$ is $T^\times$-invariant: 
$\Pi_\Sigma \circ T^\times |_{iY} = T^\times \circ \Pi_\Sigma |_{iY}$.
\\[0.2cm]
\textbf{Proof.}
This follows from Prop.3.13 (iii) and Eq.(\ref{3-23b}). \hfill $\blacksquare$
\\[-0.2cm]

Let $\lambda _0$ be an isolated generalized eigenvalue, which is separated from the rest of the generalized 
spectrum by a simple closed curve $\gamma _0 \subset \Omega \cup I \cup \{ \lambda \, | \, \mathrm{Im} (\lambda ) < 0\}$.
Let
\begin{equation}
\Pi_0 = \frac{1}{2\pi \mathrm{i} } \int_{\gamma _0}\! \mathcal{R}_\lambda d\lambda ,
\end{equation}
be a projection for $\lambda _0$ and $V_0 = \bigcup_{m\geq 1} \mathrm{Ker}\, B^{(m)}(\lambda _0) \circ
\cdots \circ B^{(1)}(\lambda _0)$ a generalized eigenspace of $\lambda _0$.
The main theorem in this paper is stated as follows:
\\[0.2cm]
\textbf{Theorem \thedef.}
If $\Pi_0 iX(\Omega )$ is finite dimensional, then $\Pi_0 iX(\Omega ) = V_0$.
\\[0.2cm]
In the usual spectral theory, this theorem is easily proved by using the resolvent equation.
In our theory, the composition $\mathcal{R}_{\lambda '} \circ \mathcal{R}_\lambda $ is not defined because
$\mathcal{R}_\lambda $ is an operator from $iX(\Omega )$ into $X(\Omega )'$.
As a result, the resolvent equation does not hold and the proof of the above theorem is rather technical.
\\[0.2cm]
\textbf{Proof.}
Let $\mathcal{R}_\lambda = \sum^\infty_{j=-\infty} (\lambda _0 - \lambda )^j E_j$ be a 
Laurent series of $\mathcal{R}_\lambda $, which converges in the strong dual topology (see Thm.A.3).
Since
\begin{eqnarray*}
id = (\lambda - T^\times) \circ \mathcal{R}_\lambda 
 = (\lambda _0 - T^\times - (\lambda _0 - \lambda ))
\circ  \sum^\infty_{j=-\infty} (\lambda _0 - \lambda )^j E_j,
\end{eqnarray*}
we obtain $E_{-n-1} = (\lambda_0 - T^\times) E_{-n}$ for $n= 1,2,\cdots $. 
Thus the equality
\begin{equation}
E_{-n-1} = (\lambda_0 - T^\times) ^nE_{-1} 
\label{3-43}
\end{equation}
holds. 
Similarly, $id|_{iY} = \mathcal{R}_\lambda \circ (\lambda -T^\times)|_{iY}$ (Prop.3.13 (ii)) provides
$E_{-n-1}|_{iY} = E_{-n} \circ (\lambda _0 - T^\times)|_{iY}$.
Thus we obtain $\mathsf{R}(E_{-n-1}|_{iY}) \subseteq \mathsf{R}(E_{-n})$ for any $n\geq 1$.
Since $Y$ is dense in $X(\Omega )$ and the range of $E_{-1} = -\Pi_0$ is finite dimensional,
it turns out that $\mathsf{R}(E_{-n}|_{iY}) = \mathsf{R}(E_{-n})$ and 
$\mathsf{R}(E_{-n-1}) \subseteq \mathsf{R}(E_{-n})$ for any $n\geq 1$.
This implies that the principle part $\sum^{-1}_{-\infty} (\lambda _0-\lambda )^jE_j$ of
the Laurent series is a finite dimensional operator. 
Hence, there exists an integer $M\geq 1$ such that $E_{-M-1} = 0$.
This means that $\lambda _0$ is a pole of $\mathcal{R}_\lambda $ :
\begin{equation}
\mathcal{R}_\lambda  = \sum^\infty_{j=-M} (\lambda _0 - \lambda )^j E_j.
\label{3-40}
\end{equation}

Next, from the equality $(id - A(\lambda )K^\times) \circ \mathcal{R}_\lambda = A(\lambda )$, we have
\begin{eqnarray*}
\left( id - \sum^\infty_{k=0} (\lambda _0 - \lambda )^k
A^{(k+1)}(\lambda _0) K^\times \right) \circ \sum^\infty_{j=-M}(\lambda _0 - \lambda )^j
E_j = \sum^\infty_{k=0}(\lambda _0 - \lambda )^k A^{(k+1)}(\lambda _0).
\end{eqnarray*}
Comparing the coefficients of $(\lambda _0 - \lambda )^{-1}$ on both sides, we obtain
\begin{equation}
(id - A(\lambda _0) K^\times) E_{-1}
 - \sum^M_{j=2} A^{(j)}(\lambda _0) K^\times E_{-j} = 0.
\end{equation}
Substituting Eq.(\ref{3-43}) and $E_{-1} = -\Pi_0$ provides
\begin{eqnarray}
B^{(1)}(\lambda _0)\Pi_0 - \sum^{M}_{j=2} A^{(j)} (\lambda _0) K^\times (\lambda _0 - T^\times )^{j-1} \Pi_0= 0.
\label{3-42}
\end{eqnarray}
In particular, this implies $\mathsf{R}(\Pi_0) \subset \mathsf{D}(B^{(1)}(\lambda _0))$.
Hence, $(\lambda _0 - T^\times ) \Pi_0$ can be rewritten as
\begin{eqnarray*}
(\lambda _0 - T^\times ) \Pi_0 = (\lambda _0 - H^\times) \circ (id - A(\lambda _0)K^\times )\Pi_0
 = (\lambda _0 - H^\times) B^{(1)}(\lambda _0)\Pi_0.
\end{eqnarray*}
Then, by using the definition of $B^{(2)}(\lambda _0)$, Eq.(\ref{3-42}) is rearranged as
\begin{eqnarray*}
B^{(2)}(\lambda _0) B^{(1)}(\lambda _0)\Pi_0
   - \sum^{M}_{j=3} A^{(j)} (\lambda _0) K^\times (\lambda _0 - T^\times )^{j-1} \Pi_0 = 0.
\end{eqnarray*}
Repeating similar calculations, we obtain
\begin{equation}
B^{(M)}(\lambda _0) \circ \cdots \circ B^{(1)}(\lambda _0)\Pi_0 = 0.
\label{3-43b}
\end{equation}
This proves $\Pi_0iX(\Omega ) \subset V_0$.

Let us show $\Pi_0iX(\Omega ) \supset V_0$.
From the equality $\mathcal{R}_\lambda \circ (id - K^\times A(\lambda )) = A(\lambda )$, 
we have
\begin{equation}
\sum^\infty_{j=-M} (\lambda _0 - \lambda )^j E_j \circ
\left( id - K^\times\sum^\infty_{k=0} (\lambda _0 - \lambda ) ^k A^{(k+1)}(\lambda _0)\right)
 = \sum^\infty_{k=0} (\lambda _0 - \lambda ) ^k A^{(k+1)}(\lambda _0).
\label{3-47}
\end{equation}
Comparing the coefficients of $(\lambda _0 - \lambda )^k$ on both sides for $k=1,2,\cdots $,
we obtain
\begin{eqnarray}
E_k (id - K^\times A(\lambda _0))\phi - \sum^\infty_{j=1} E_{-j+k} K^\times A^{(j+1)}(\lambda _0) \phi
 = A^{(k+1)}(\lambda _0)\phi,
\label{3-47b}
\end{eqnarray}
for any $\phi \in iX(\Omega )$, where the left hand side is a finite sum.
Note that $K^\times A^{(j)}(\lambda _0)iX(\Omega ) \subset iX(\Omega )$ for any $j=1,2,\cdots $
because $K^\times A(\lambda )iX(\Omega ) \subset iX(\Omega )$ for any $\lambda $ (the assumption (X8)).

Now suppose that $\mu \in V_0$ is a generalized eigenfunction satisfying 
$B^{(M)}(\lambda_0 ) \circ \cdots \circ B^{(1)}(\lambda _0)\mu = 0$.
For this $\mu$, we need the following lemma.
\\[0.2cm]
\textbf{Lemma.} For any $k=0,1,\cdots ,M-1$,
\\
(i) $(\lambda _0 - T^\times)^k \mu = (\lambda _0 - H^\times)^kB^{(k)}(\lambda_0 ) \circ \cdots \circ B^{(1)}(\lambda _0)\mu$.
\\
(ii) $K^\times (\lambda _0 - T^\times)^k\mu \in iX(\Omega )$.
\\[0.2cm]
\textbf{Proof.}
Due to Thm.3.9, $\mu $ is included in the domain of $(\lambda _0 - T^\times)^k$.
Thus the left hand side of (i) indeed exists.
Then, we have
\begin{eqnarray*}
(\lambda _0 - H^\times)^kB^{(k)}(\lambda_0 )
&=& (\lambda _0 - H^\times)^k(id - A^{(k)}(\lambda _0)K^\times (\lambda _0 - H^\times)^{k-1}) \\
&=& (\lambda _0 - H^\times - K^\times) (\lambda _0 - H^\times)^{k-1}
 = (\lambda _0 - T^\times ) (\lambda _0 - H^\times)^{k-1}.
\end{eqnarray*}
Repeating this procedure yields (i).
To prove (ii), let us calculate
\begin{eqnarray*}
0 = K^\times (\lambda _0 - H^\times)^k B^{(M)}(\lambda_0 ) \circ \cdots \circ B^{(1)}(\lambda _0)\mu.
\end{eqnarray*}
Eq.(\ref{3-25}) and the part (i) of this lemma give
\begin{eqnarray*}
0 &=& K^\times B^{(M-k)}(\lambda_0 ) \circ \cdots \circ B^{(k+1)}(\lambda _0) \circ
      (\lambda _0 - H^\times)^k \circ B^{(k)}(\lambda_0 ) \circ \cdots \circ B^{(1)}(\lambda _0)\mu \\
&=& K^\times B^{(M-k)}(\lambda_0 ) \circ \cdots \circ B^{(k+1)}(\lambda _0) \circ (\lambda _0 - T^\times)^k \mu.
\end{eqnarray*}
For example, when $k= M-1$, this is reduced to
\begin{eqnarray*}
0 = K^\times (id - A(\lambda _0)K^\times) \circ (\lambda _0 - T^\times)^{M-1}\mu.
\end{eqnarray*}
This proves $K^\times (\lambda _0 - T^\times)^{M-1}\mu = K^\times A(\lambda _0)K^\times (\lambda _0 - T^\times)^{M-1}\mu
\in iX(\Omega )$.
This is true for any $k=0,1,\cdots ,M-1$;
it follows from the definition of $B^{(j)}(\lambda _0)$'s that $K^\times (\lambda _0 - T^\times)^{k}\mu$ is expressed as
a linear combination of elements of the form $K^\times A^{(j)}(\lambda _0) \xi_j,\, \xi_j\in iX(\Omega )$.
Since $K^\times A^{(j)}(\lambda _0)iX(\Omega ) \subset iX(\Omega )$, 
we obtain $K^\times (\lambda _0 - T^\times)^k\mu \in iX(\Omega )$. \hfill $\blacksquare$
\\[-0.2cm]

Since $K^\times (\lambda _0 - T^\times)^k\mu \in iX(\Omega )$, we can substitute $\phi = K^\times (\lambda _0 - T^\times)^k\mu $
into Eq.(\ref{3-47b}).
The resultant equation is rearranged as
\begin{eqnarray*}
& & E_kK^\times (id - A(\lambda _0)K^\times) (\lambda _0 - T^\times)^k\mu
  - \left( id + \sum^k_{j=1} E_{-j+k} K^\times (\lambda _0 - H^\times )^{k-j} \right) A^{(k+1)} (\lambda _0) K^\times 
     (\lambda _0 - T^\times)^k\mu \\
&= & \sum^\infty_{j=k+1} E_{-j+k} K^\times A^{(j+1)}(\lambda _0) K^\times (\lambda _0 - T^\times)^k\mu.
\end{eqnarray*}
Further, $(\lambda _0 - T^\times)^k
 = (\lambda _0 - H^\times ) ^k B^{(k)}(\lambda _0) \circ \cdots \circ B^{(1)}(\lambda _0)$ provides
\begin{eqnarray}
& & E_kK^\times (\lambda _0- H^\times )^k B^{(k+1)}(\lambda_0) \circ \cdots \circ B^{(1)}(\lambda_0)\mu
 \nonumber \\
& & \quad - \left( id + \sum^k_{j=1} E_{-j+k} K^\times (\lambda _0 - H^\times )^{k-j} \right)
  A^{(k+1)} (\lambda _0) K^\times 
       (\lambda _0 - H^\times ) ^k B^{(k)}(\lambda _0) \circ \cdots \circ B^{(1)}(\lambda _0) \mu\nonumber \\
&= & \sum^\infty_{j=k+1} E_{-j+k} K^\times A^{(j+1)}(\lambda _0) K^\times
(\lambda _0 - T^\times)^k\mu.
\label{3-48}
\end{eqnarray}
On the other hand, comparing the coefficients of $(\lambda _0 - \lambda )^0$ of Eq.(\ref{3-47}) provides
\begin{eqnarray*}
E_0 (id - K^\times A(\lambda _0))\phi
    - \sum^\infty_{j=1} E_{-j} K^\times A^{(j+1)}(\lambda _0) \phi= A(\lambda _0)\phi,
\end{eqnarray*}
for any $\phi \in iX(\Omega )$.
Substituting $\phi = K^\times \mu \in iX(\Omega )$ provides
\begin{equation}
(id + E_0 K^\times) B^{(1)}(\lambda _0)\mu
 = \mu + \sum^\infty_{j=1} E_{-j} K^\times A^{(j+1)}(\lambda _0) K^\times \mu.
\label{3-49}
\end{equation}
By adding Eq.(\ref{3-49}) to Eqs.(\ref{3-48}) for $k=1, \cdots ,M-1$, we obtain
\begin{eqnarray}
& & (id + E_0 K^\times) B^{(1)}(\lambda _0) \mu \nonumber \\
& &\quad   - \sum^{M-1}_{k=1} \left( id + \sum^k_{j=1} E_{-j+k} K^\times (\lambda _0 - H^\times)^{k-j} \right)
  A^{(k+1)} (\lambda _0) K^\times 
       (\lambda _0 - H^\times ) ^k B^{(k)}(\lambda _0) \circ \cdots \circ B^{(1)}(\lambda _0) \mu \nonumber \\
& & \quad \quad + \sum^{M-1}_{k=1} 
 E_kK^\times (\lambda _0- H^\times )^k B^{(k+1)}(\lambda_0) \circ \cdots \circ B^{(1)}(\lambda_0)\mu \nonumber \\
 &=& \mu + \sum^{M-1}_{k=0} \sum^\infty_{j=1} E_{-j}K^\times A^{(j+k+1)}(\lambda _0) K^\times
(\lambda _0 - T^\times)^k \mu.
\label{3-49b}
\end{eqnarray}
The left hand side above is rewritten as
\begin{eqnarray*}
& & \left( id + E_0 K^\times + E_1 K^\times (\lambda _0 - H^\times )\right)
B^{(2)}(\lambda _0) B^{(1)}(\lambda _0) \mu\\
& & \quad - \sum^{M-1}_{k=2} \left( id + \sum^k_{j=1} E_{-j+k} K^\times (\lambda _0 - H^\times )^{k-j} \right)
  A^{(k+1)} (\lambda _0) K^\times 
       (\lambda _0 - H^\times ) ^k B^{(k)}(\lambda _0) \circ \cdots \circ B^{(1)}(\lambda _0) \mu \\
& & \quad \quad + \sum^{M-1}_{k=2} 
   E_kK^\times (\lambda _0 -H^\times )^k B^{(k+1)}(\lambda_0) \circ \cdots \circ B^{(1)}(\lambda_0) \mu.
\end{eqnarray*}
Repeating similar calculations, we can verify that Eq.(\ref{3-49b}) is rewritten as
\begin{eqnarray}
& & \left( id + \sum^{M-1}_{j=0} E_j K^\times (\lambda _0 -H^\times )^{j} \right)
B^{(M)}(\lambda _0) \circ \cdots \circ B^{(1)}(\lambda _0) \mu \nonumber \\
 &=& \mu - \sum^{M-1}_{k=0} \sum^\infty_{j=1}E_{-j}K^\times A^{(j+k+1)}(\lambda _0) K^\times
(\lambda _0 - T^\times)^k \mu.
\label{3-50}
\end{eqnarray}
Since $B^{(M)}(\lambda _0) \circ \cdots \circ B^{(1)}(\lambda _0) \mu = 0$, we obtain
\begin{eqnarray*}
\mu = \sum^{M-1}_{k=0} \sum^\infty_{j=1}E_{-j} K^\times A^{(j+k+1)}(\lambda _0) K^\times
(\lambda _0 - T^\times)^k \mu.
\end{eqnarray*}
Since $\mathsf{R}(E_{-j}) \subset \mathsf{R}(E_{-1}) = \mathsf{R}(\Pi_0)$, this proves $\Pi_0 iX(\Omega ) \supset V_0$.
Thus the proof of $\Pi_0iX(\Omega ) = V_0$ is completed. \hfill $\blacksquare$


\subsection{Properties of the generalized spectrum}

We show a few criteria to estimate the generalized spectrum.
Recall that $\hat{\sigma}_p (T) \subset \sigma _p(T^\times)$ because of Thm.3.5.
The relation between $\hat{\sigma} (T)$ and $\sigma (T)$ is given as follows.
\\[0.2cm]
\textbf{Proposition \thedef.}
Let $\mathbf{C}_- = \{ \mathrm{Im}(\lambda )<0\}$ be an open lower half plane.
Let $\sigma _p(T)$ and $\sigma (T)$ be the point spectrum and the spectrum in the usual sense, respectively.
Then, the following relations hold.
\\
(i) $\hat{\sigma} (T) \cap \mathbf{C}_- \subset \sigma (T)\cap \mathbf{C}_-$.
In particular, $\hat{\sigma}_p (T) \cap \mathbf{C}_- \subset \sigma _p(T)\cap \mathbf{C}_-$
\\
(ii) Let $\Sigma \subset \mathbf{C}_-$ be a bounded subset of $\sigma (T)$ which is separated from the rest of 
the spectrum by a simple closed curve $\gamma $.
Then, there exists a point of $\hat{\sigma }(T)$ inside $\gamma $.
In particular, if $\lambda \in \mathbf{C}_-$ is an isolated point of $\sigma (T)$, then $\lambda \in \hat{\sigma }(T)$.
\\[0.2cm]
\textbf{Proof.} 
Note that when $\lambda \in \mathbf{C}_-$, the generalized resolvent satisfies 
$\mathcal{R}_\lambda \circ i =i \circ (\lambda -T)^{-1} $ due to Thm.3.12.

(i) Suppose that $\lambda \in \varrho (T) \cap \mathbf{C}_-$, where $\varrho (T)$ is the resolvent set of $T$ in the usual sense.
Since $\mathcal{H}$ is a Hilbert space, there is a neighborhood $V_\lambda \subset \varrho (T) \cap \mathbf{C}_-$ of $\lambda $
such that $(\lambda ' - T)^{-1}$ is continuous on $\mathcal{H}$ for any $\lambda '\in V_\lambda $ and the set 
$\{ (\lambda ' - T)^{-1}\psi \}_{\lambda ' \in V_\lambda }$ is bounded in $\mathcal{H}$ for each $\psi \in X(\Omega )$.
Since $i:\mathcal{H} \to X(\Omega )'$ is continuous and since the topology of $X(\Omega )$ is stronger than that of $\mathcal{H}$,
$\mathcal{R}_{\lambda '}\circ i = i \circ (\lambda ' -T)^{-1}$ is a continuous operator from $X(\Omega )$ into $X(\Omega )'$
for any $\lambda ' \in V_\lambda $,
and the set $\{ \mathcal{R}_{\lambda '} \circ i \psi \}_{\lambda ' \in V_\lambda }$ is bounded in $X(\Omega )'$.
This proves that $\lambda \in \hat{\varrho} (T) \cap \mathbf{C}_-$.

Next, suppose that $\lambda \in \mathbf{C}_-$ is a generalized eigenvalue satisfying
$(id - K^\times A(\lambda ) ) i (\psi) =0$ for $\psi \in X(\Omega )$.
Since $\lambda -H$ is invertible on $\mathcal{H}$ when $\lambda \in \mathbf{C}_-$, putting $\phi = (\lambda -H)^{-1}\psi$
provides
\begin{eqnarray*}
(id - K^\times A(\lambda ) ) i (\lambda -H)\phi = (i(\lambda -H) - K^\times i) \phi = i(\lambda -T)\phi = 0,
\end{eqnarray*}
and thus $\lambda \in \sigma _p(T)$.

(ii) Let $\mathcal{P}$ be the Riesz projection for $\Sigma \subset \sigma (T)\cap \mathbf{C}_-$, which is
defined as $\mathcal{P} = (2\pi \mathrm{i} )^{-1} \int_{\gamma }\! (\lambda -T)^{-1} d\lambda $.
Since $\gamma $ encloses a point of $\sigma (T)$, $\mathcal{P} \mathcal{H} \neq \emptyset$.
Since $X(\Omega )$ is dense in $\mathcal{H}$, $\mathcal{P}X(\Omega ) \neq \emptyset$.
This fact and $\mathcal{R}_\lambda \circ i =i \circ (\lambda -T)^{-1}$ prove that the range of the generalized Riesz projection
defined by Eq.(\ref{3-33}) is not zero.
Hence, the closed curve $\gamma $ encloses a point of $\hat{\sigma }(T)$.  
 \hfill $\blacksquare$
\\[-0.2cm]

A few remarks are in order.
If the spectrum of $T$ on the lower half plane consists of discrete eigenvalues,
(i) and (ii) show that $\sigma _p(T)\cap \mathbf{C}_- = \sigma (T)\cap \mathbf{C}_- = \hat{\sigma }(T)\cap \mathbf{C}_-$.
However, it is possible that a generalized eigenvalue on $I$ is not an eigenvalue in the usual sense.
See \cite{Chi} for such an example.
In most cases, the continuous spectrum on the lower half plane is not included in the 
generalized spectrum because the topology on $X(\Omega )'$ is weaker than that on $\mathcal{H}$, 
although the point spectrum and the residual spectrum may remain to exist as the generalized spectrum.
Note that the continuous spectrum on the interval $I$ also disappears;
for the resolvent $(\lambda -T)^{-1} = (\lambda -H)^{-1}(id - K(\lambda -H)^{-1})^{-1}$ in the usual sense,
the factor $(\lambda -H)^{-1}$ induces the continuous spectrum on the real axis because $H$ is selfadjoint.
For the generalized resolvent, $(\lambda -H)^{-1}$ is replaced by $A(\lambda )$, which has no singularities.
This suggests that obstructions when calculating the Laplace inversion formula by using the residue theorem
may disappear.

Recall that a linear operator $L$ from a topological vector space $X_1$ to another topological vector space $X_2$
is said to be bounded if there exists a neighborhood $U\subset X_1$ such that $LU \subset X_2$ is a bounded set.
When $L=L(\lambda )$ is parameterized by $\lambda $, it is said to be bounded uniformly in $\lambda $
if such a neighborhood $U$ is independent of $\lambda $.
When the domain $X_1$ is a Banach space, $L(\lambda )$ is bounded uniformly in $\lambda $
if and only if $L(\lambda )$ is continuous for each $\lambda $
($U$ is taken to be the unit sphere).
Similarly, $L$ is called compact if there exists a neighborhood $U\subset X_1$ such that $LU \subset X_2$ is relatively compact.
When $L=L(\lambda )$ is parameterized by $\lambda $, it is said to be compact uniformly in $\lambda $
if such a neighborhood $U$ is independent of $\lambda $.
When the domain $X_1$ is a Banach space, $L(\lambda )$ is compact uniformly in $\lambda $
if and only if $L(\lambda )$ is compact for each $\lambda $.
When the range $X_2$ is a Montel space, a (uniformly) bounded operator is (uniformly) compact 
because every bounded set in a Montel space is relatively compact.
Put $\hat{\Omega } := \{ \mathrm{Im}(\lambda )< 0\} \cup I \cup \Omega $ as before.
In many applications, $i^{-1}K^\times A(\lambda )i$ is a bounded operator.
In such a case, the following proposition is useful to estimate the generalized spectrum.
\\[0.2cm]
\textbf{Proposition \thedef.} 
Suppose that for $\lambda \in \hat{\Omega }$, there exists a neighborhood $U_\lambda \subset \hat{\Omega }$ of 
$\lambda $ such that $i^{-1}K^\times A(\lambda ')i : X(\Omega ) \to X(\Omega )$ is a bounded operator 
uniformly in $\lambda ' \in U_\lambda $.
If $id - i^{-1}K^\times A(\lambda )i$ has a continuous inverse on $X(\Omega )$, then $\lambda \notin \hat{\sigma }(T)$.
\\[0.2cm]
\textbf{Proof.}
Note that $\mathcal{R}_\lambda \circ i$ is rewritten as 
$\mathcal{R}_\lambda \circ i = A(\lambda )\circ i \circ (id - i^{-1}K^\times A(\lambda )i )^{-1}$.
Since $A(\lambda )\circ i$ is continuous, it is sufficient to prove that there exists a neighborhood $V_\lambda $ of $\lambda $
such that the set $\{ (id - i^{-1}K^\times A(\lambda ')i )^{-1}\psi \}_{\lambda '\in V_\lambda }$
is bounded in $X(\Omega )$ for each $\psi \in X(\Omega )$.
For this purpose, it is sufficient to prove that the mapping $\lambda ' \mapsto (id - i^{-1}K^\times A(\lambda ')i )^{-1}\psi$
is continuous in $\lambda ' \in V_\lambda $.
Since $i^{-1}K^\times A(\lambda )i$ is holomorphic (see the proof of Thm.3.12), 
there is an operator $D(\lambda ,h)$ on $X(\Omega )$ such that
\begin{eqnarray*}
id - i^{-1}K^\times A(\lambda +h)i &=& id - i^{-1}K^\times A(\lambda )i - h D(\lambda , h) \\
&=& \left( id - hD(\lambda ,h) (id - i^{-1}K^\times A(\lambda )i)^{-1}\right) \circ
(id - i^{-1}K^\times A(\lambda )i ).
\end{eqnarray*}
Since $ i^{-1}K^\times A(\lambda )i$ is a bounded operator uniformly in $\lambda \in U_\lambda $, $D(\lambda , h)$ is a bounded 
operator when $h$ is sufficiently small.
Since $(id - i^{-1}K^\times A(\lambda )i)^{-1}$ is continuous by the assumption,
$D(\lambda ,h) (id - i^{-1}K^\times A(\lambda )i)^{-1}$ is a bounded operator.
Then, Bruyn's theorem \cite{Bru} shows that
$id - hD(\lambda ,h) (id - i^{-1}K^\times A(\lambda )i)^{-1}$ has a continuous inverse for sufficiently small $h$ 
and the inverse is continuous in $h$
(when $X(\Omega )$ is a Banach space, Bruyn's theorem is reduced to the existence of the Neumann series).
This proves that $(id - i^{-1}K^\times A(\lambda +h)i )^{-1}\psi$ exists and continuous in $h$ for each $\psi$. 
\hfill $\blacksquare$
\\[0.2cm]
As a corollary, if $X(\Omega )$ is a Banach space and $ i^{-1}K^\times A(\lambda )i$ is a continuous operator on $X(\Omega )$
for each $\lambda $, then $\lambda \in \hat{\varrho}(T)$ if and only if $id - i^{-1}K^\times A(\lambda )i$
has a continuous inverse on $X(\Omega )$.
Because of this proposition, we can apply the spectral theory on locally convex spaces 
(for example, \cite{All,DeV,Moo,Ola,Rin,Sch}) to the operator $id - i^{-1}K^\times A(\lambda )i$
to estimate the generalized spectrum.
In particular, like as Riesz-Schauder theory in Banach spaces, we can prove the next theorem.
\\[0.2cm]
\textbf{Theorem \thedef.}
In addition to (X1) to (X8), suppose that $i^{-1} K^\times A(\lambda )i : X(\Omega ) \to X(\Omega )$ is a compact operator
uniformly in $\lambda \in \hat{\Omega } := \{ \mathrm{Im}(\lambda )< 0\} \cup I \cup \Omega $.
Then, the following statements are true.
\\ 
(i) For any compact set $D \subset \hat{\Omega }$, the number of generalized eigenvalues in $D$ is finite
(thus $\hat{\sigma }_p(T)$ consists of a countable number of generalized eigenvalues and they may accumulate
only on the boundary of $\hat{\Omega }$ or infinity).
\\
(ii) For each $\lambda _0\in \hat{\sigma }_p(T)$, the generalized eigenspace $V_0$ is of finite dimensional and 
$\Pi_0iX(\Omega ) = V_0$.
\\
(iii) $\hat{\sigma }_c(T) = \hat{\sigma }_r(T) = \emptyset$.
\\[-0.2cm]

If $X(\Omega )$ is a Banach space, the above theorem follows from well known Riesz-Schauder theory.
Even if $X(\Omega )$ is not a Banach space, we can prove the same result (see below).
Thm.3.19 is useful to find embedded eigenvalues of $T$:
\\[0.2cm]
\textbf{Corollary \thedef.}
Suppose that $T$ is selfadjoint.
Under the assumptions in Thm.3.19, 
the number of eigenvalues of $T = H+K$ (in $\mathcal{H}$-sense) in any compact set $D \subset I$ is finite.
Their algebraic multiplicities $\mathrm{dim}\, \mathrm{Ker}\, (\lambda -T)$ are finite.
\\[0.2cm]
\textbf{Proof.}
Let $\lambda_0 \in I$ be an eigenvalue of $T$.
It is known that the projection $\mathcal{P}_0$ to the corresponding eigenspace is given by
\begin{equation}
\mathcal{P}_0\phi = \lim_{\varepsilon \to -0} \mathrm{i} \varepsilon \cdot (\lambda _0 + \mathrm{i} \varepsilon -T)^{-1} \phi,
\quad \phi \in \mathcal{H},
\end{equation}  
where the limit is taken with respect to the topology on $\mathcal{H}$.
When $\mathrm{Im} (\lambda ) < 0$, we have $\mathcal{R}_\lambda i (\phi)= i(\lambda -T)^{-1} \phi$
for $\phi \in X(\Omega )$.
This shows
\begin{eqnarray*}
i \circ \mathcal{P}_0 \phi = \lim_{\varepsilon \to -0}
  \mathrm{i} \varepsilon \cdot \mathcal{R}_{\lambda _0 + \mathrm{i} \varepsilon } \circ i (\phi),
\quad \phi\in X(\Omega ).
\end{eqnarray*}
Let $\mathcal{R}_\lambda = \sum^\infty_{j=-\infty} (\lambda _0- \lambda )^j E_{j}$ be the Laurent expansion of 
$\mathcal{R}_\lambda $, which converges around $\lambda _0$.
This provides
\begin{eqnarray*}
i \circ \mathcal{P}_0 = \lim_{\varepsilon \to -0} \mathrm{i} \varepsilon 
\sum^\infty_{j=-\infty} (-\mathrm{i} \varepsilon )^j E_{j} \circ i.
\end{eqnarray*}
Since the right hand side converges with respect to the topology on $X(\Omega )'$, we obtain
\begin{equation}
i \circ \mathcal{P}_0 = -E_{-1} \circ i = \Pi_0 \circ i, \quad E_{-2} = E_{-3} = \cdots =0,
\end{equation}
where $\Pi_0$ is the generalized Riesz projection for $\lambda _0$.
Since $\lambda _0$ is an eigenvalue, $\mathcal{P}_0\mathcal{H} \neq \emptyset$.
Since $X(\Omega )$ is a dense subspace of $\mathcal{H}$, $\mathcal{P}_0X(\Omega ) \neq \emptyset$.
Hence, we obtain $\Pi_0 iX(\Omega ) \neq \emptyset$, which implies that 
$\lambda _0$ is a generalized eigenvalue; $\sigma _p(T) \subset \hat{\sigma }_p(T)$.
Since $\hat{\sigma }_p(T)$ is countable, so is $\sigma _p(T)$. 
Since $\Pi_0 iX(\Omega )$ is a finite dimensional space, so is $\mathcal{P}_0X(\Omega )$.
Then, $\mathcal{P}_0\mathcal{H} = \mathcal{P}_0X(\Omega )$ proves to be finite dimensional
because $\mathcal{P}_0\mathcal{H}$ is the closure of $\mathcal{P}_0X(\Omega )$. 
\hfill $\blacksquare$
\\[0.2cm]
Our results are also useful to calculate eigenvectors for embedded eigenvalues.
In the usual Hilbert space theory, if an eigenvalue $\lambda $ is embedded in the continuous spectrum of $T$,
we can not apply the Riesz projection for $\lambda $ because there are no closed curves in $\mathbf{C}$ which
separate $\lambda $ from the rest of the spectrum.
In our theory, $\hat{\sigma }_c(T) = \hat{\sigma }_r(T) = \emptyset$.
Hence, the generalized eigenvalues are indeed isolated and the Riesz projection $\Pi_0$ is applied to yield
$\Pi_0 iX(\Omega ) = V_0$.
Then, the eigenspace in $\mathcal{H}$-sense is obtained as $V_0\cap \mathsf{D}(T)$.
\\[0.2cm]
\textbf{Proof of Thm.3.19.}
The theorem follows from Riesz-Schauder theory on locally convex spaces developed in Ringrose \cite{Rin}.
Here, we give a simple review of the argument in \cite{Rin}. We denote $X(\Omega ) = X$ and 
$i^{-1}K^\times A(\lambda ) i = C(\lambda )$ for simplicity.
A pairing for $(X', X)$ is denoted by $\langle \, \cdot\, \,|\, \, \cdot \, \rangle_X$.

Since $C(\lambda ): X\to X$ is compact uniformly in $\lambda $, there exists a neighborhood $V$ of zero in $X$,
which is independent of $\lambda $, such that 
$C(\lambda )V \subset X$ is relatively compact.
Put $p(x) = \inf \{ |\lambda | ; x\in \lambda V\}$.
Then, $p$ is a continuous semi-norm on $X$ and $V = \{ x\, | \, p(x) < 1\}$.
Define a closed subspace $M$ in $X$ to be
\begin{equation}
M = \{ x\in X \, | \, p(x) = 0 \} \subset V.
\end{equation}
Let us consider the quotient space $X/M$, whose elements are denoted by $[x]$.
The semi-norm $p$ induces a norm $P$ on $X/M$ by $P([x]) = p(x)$.
If $X/M$ is equipped with the norm topology induced by $P$, we denote the space as $\mathcal{B}$.
The completion of $\mathcal{B}$, which is a Banach space, is denoted by $\mathcal{B}_0$.
The dual space $\mathcal{B}_0'$ of $\mathcal{B}_0$ is a Banach space with the norm
\begin{equation}
|| \mu ||_{\mathcal{B}_0'} := \sup_{P([x]) < 1} |\langle \mu \,|\, [x] \rangle _{\mathcal{B}_0}|,
\end{equation}
where $\langle \, \cdot\, \,|\, \,\cdot \, \rangle_{\mathcal{B}_0}$ is a pairing for $(\mathcal{B}_0', \mathcal{B}_0)$.
Define a subspace $S\subset X'$ to be
\begin{equation}
S = \{ \mu \in X'\, | \, \sup_{x\in V} |\langle \mu \,|\, x \rangle _X | < \infty \}.
\end{equation}
The linear mapping $\hat{\,\,} : S \to \mathcal{B}_0'\,\, (\mu \mapsto \hat{\mu})$ defined through
$\langle \hat{\mu} \,|\, [x] \rangle_{\mathcal{B}_0} = \langle \mu \,|\, x \rangle_X$ is bijective.
Define the operator $Q(\lambda ) : \mathcal{B} \to \mathcal{B}$ to be $Q(\lambda )[x] = [C(\lambda )x]$.
Then, the equality
\begin{equation}
\langle \hat{\mu} \,|\, Q(\lambda )[x] \rangle_{\mathcal{B}_0} = \langle \mu \,|\, C(\lambda )x \rangle_X
\label{3-50b}
\end{equation}
holds for $\mu \in S$ and $x\in X$.
Let $Q_0(\lambda ) : \mathcal{B}_0 \to \mathcal{B}_0$ be a continuous extension of $Q(\lambda )$.
Then, $Q_0(\lambda )$ is a compact operator on a Banach space, and thus the usual Riesz-Schauder theory is applied.
By using Eq.(\ref{3-50b}), it is proved that $z\in \mathbf{C}$ is an eigenvalue of $C(\lambda )$ if and only if 
it is an eigenvalue of $Q_0(\lambda )$.
In this manner, we can prove that 
\\[0.2cm]
\textbf{Theorem \thedef \, \cite{Rin}.}
The number of eigenvalues of the operator $C(\lambda ) : X \to X$ 
is at most countable, which can accumulate only at the origin.
The eigenspaces $\bigcup_{m\geq 1} \mathrm{Ker}\, (z -C(\lambda ))^m$ of nonzero eigenvalues $z$ are finite dimensional.
If $z \neq 0$ is not an eigenvalue, $z -C(\lambda )$ has a continuous inverse on $X$.
See \cite{Rin} for the complete proof.
\\[0.2cm]
\indent 
Now we are in a position to prove Thm.3.19.
Suppose that $\lambda $ is not a generalized eigenvalue.
Then, $1$ is not an eigenvalue of $C(\lambda ) = i^{-1}K^\times A(\lambda )i$.
The above theorem concludes that $id - C(\lambda )$ has a continuous inverse on $X(\Omega )$.
Since $C(\lambda )$ is compact uniformly in $\lambda $, Prop.3.18 implies $\lambda \notin \hat{\sigma }(T)$.
This proves the part (iii) of Thm.3.19.

Let us show the part (i) of the theorem.
Let $z = z(\lambda )$ be an eigenvalue of $C(\lambda )$.
We suppose that $z(\lambda _0) = 1$ so that $\lambda _0$ is a generalized eigenvalue.
As was proved in the proof of Thm.3.12, 
$\langle \mu \,|\, C(\lambda )x \rangle_X$ is holomorphic in $\lambda$.
Eq.(\ref{3-50b}) shows that
$\langle \hat{\mu} \,|\, Q(\lambda )[x] \rangle_{\mathcal{B}_0}$ is holomorphic for any $\hat{\mu} \in \mathcal{B}_0'$
and $[x]\in \mathcal{B}$.
Since $\mathcal{B}_0$ is a Banach space and $\mathcal{B}$ is dense in $\mathcal{B}_0$,
$Q_0(\lambda )$ is a holomorphic family of operators.
Recall that the eigenvalue $z(\lambda )$ of $C(\lambda )$ is also an eigenvalue of $Q_0(\lambda )$ 
satisfying $z(\lambda _0) = 1$.
Then, the analytic perturbation theory of operators (see Chapter VII of Kato \cite{Kato}) shows that
there exists a natural number $p$ such that $z(\lambda )$ is holomorphic as a function of $(\lambda -\lambda _0)^{1/p}$.
Let us show that $z(\lambda )$ is not a constant function.
If $z(\lambda ) \equiv 1$, every point in $\hat{\Omega }$ is a generalized eigenvalue.
Due to Prop.3.17, the open lower half plane is included in the point spectrum of $T$.
Hence, there exists $f = f_\lambda $ in $\mathcal{H}$ such that $f = K(\lambda -H)^{-1}f$ for any $\lambda \in \mathbf{C}_-$.
However, since $K$ is $H$-bounded, there exist nonnegative numbers $a$ and $b$ such that 
\begin{eqnarray*}
|| K(\lambda -H)^{-1} || \leq a || (\lambda -H)^{-1} || + b|| H(\lambda -H)^{-1} ||
 = a|| (\lambda -H)^{-1} || + b|| \lambda (\lambda -H)^{-1} - id||,
\end{eqnarray*}
which tends to zero as $|\lambda | \to \infty$ outside the real axis.
Therefore, $|| f || \leq || K(\lambda -H)^{-1} || \cdot || f || \to 0$, which contradicts with the assumption.
Since $z(\lambda )$ is not a constant, there exists a neighborhood $U\subset \mathbf{C}$ of 
$\lambda _0$ such that $z(\lambda ) \neq 1$ when $\lambda \in U$ and $\lambda \neq \lambda _0$.
This implies that $\lambda \in U \backslash \{ \lambda _0 \}$ is not a generalized eigenvalue and the part (i)
of Thm.3.19 is proved.

Finally, let us prove the part (ii) of Thm.3.19.
Put $\tilde{C}(z) = (z-1)\cdot id + C(z)$ and $\tilde{Q}(z) = (z-1)\cdot id + Q(z)$.
They satisfy $\langle \hat{\mu} \,|\, \tilde{Q}(\lambda )[x] \rangle_{\mathcal{B}_0}
 = \langle \mu \,|\, \tilde{C}(z)x \rangle_{X}$ and 
\begin{eqnarray*}
\langle \hat{\mu} \,|\, (\lambda -\tilde{Q}(z))^{-1}[x] \rangle_{\mathcal{B}_0}
 = \langle \mu \,|\, (\lambda -\tilde{C}(z))^{-1}x \rangle_{X}.
\end{eqnarray*}
Since an eigenspace of $Q(z)$ is finite dimensional, an eigenspace of $\tilde{Q}(z)$ is also finite dimensional.
Thus the resolvent $(\lambda -\tilde{Q}(z))^{-1}$ is meromorphic in $\lambda \in \hat{\Omega }$.
Since $\tilde{Q}(z)$ is holomorphic, $(\lambda -\tilde{Q}(\lambda ))^{-1}$ is also meromorphic.
The above equality shows that $\langle \mu \,|\, (\lambda -\tilde{C}(\lambda ))^{-1}x \rangle_{X}$ is meromorphic for any $\mu \in S$.
Since $S$ is dense in $X'$, it turns out that $(\lambda -\tilde{C}(\lambda ))^{-1}x$ is meromorphic
with respect to the topology on $X$. 
Therefore, the generalized resolvent
\begin{equation}
\mathcal{R}_\lambda \circ i = A(\lambda ) \circ i \circ  (id - i^{-1}K^\times A(\lambda )i)^{-1} 
  = A(\lambda ) \circ i \circ (\lambda -\tilde{C}(\lambda ))^{-1}
\end{equation}
is meromorphic on $\hat{\Omega }$.
Now we have shown that the Laurent expansion of $\mathcal{R}_\lambda $ is of the form (\ref{3-40}) for some $M\geq 0$.
Then, we can prove Eq.(\ref{3-43b}) by the same way as the proof of Thm.3.16.
To prove that $\Pi_0 iX(\Omega )$ is of finite dimensional, we need the next lemma.
\\[0.2cm]
\textbf{Lemma \thedef.}
$\mathrm{dim}\, \mathrm{Ker}\, B^{(n)} (\lambda ) \leq \mathrm{dim}\, \mathrm{Ker}\, (id - K^\times A(\lambda ))$
for any $n\geq 1$.
\\[0.2cm]
\textbf{Proof.} Suppose that $B^{(n)}(\lambda ) \mu = 0$ with $\mu \neq 0$.
Then, we have
\begin{eqnarray*}
K^\times (\lambda - H^\times)^{n-1} B^{(n)}(\lambda ) \mu
&=& K^\times (\lambda - H^\times)^{n-1} (id - A^{(n)}(\lambda )K^\times (\lambda -H^\times)^{n-1}) \mu \\
&=& (id - K^\times A(\lambda )) \circ K^\times (\lambda -H^\times)^{n-1} \mu = 0.
\end{eqnarray*}
If $K^\times (\lambda -H^\times)^{n-1} \mu = 0$, $B^{(n)}(\lambda ) \mu = 0$ yields
$\mu = A^{(n)}(\lambda ) K^\times (\lambda -H^\times)^{n-1} \mu = 0$,
which contradicts with the assumption $\mu \neq 0$.
Thus we obtain $K^\times (\lambda -H^\times)^{n-1} \mu \in \mathrm{Ker}\, (id - K^\times A(\lambda ))$
and the mapping $\mu \mapsto K^\times (\lambda -H^\times)^{n-1} \mu$ is one-to-one. \hfill $\blacksquare$
\\[-0.2cm]

Due to Thm.3.21, $\mathrm{Ker}\, (id - K^\times A(\lambda ))$ is of finite dimensional.
Hence, $\mathrm{Ker}\, B^{(n)} (\lambda ) $ is also finite dimensional for any $n\geq 1$.
This and Eq.(\ref{3-43b}) prove that $\Pi_0iX(\Omega )$ is a finite dimensional space.
By Thm.3.16, $\Pi_0iX(\Omega ) = V_0$, which completes the proof of Thm.3.19 (ii).
\hfill $\blacksquare$

\subsection{Semigroups}

In this subsection, we suppose that
\\
\textbf{(S1)} The operator $\mathrm{i} T = \mathrm{i} (H +K)$ generates a $C^0$-semigroup 
$e^{\mathrm{i} Tt}$ on $\mathcal{H}$ (recall $\mathrm{i}= \sqrt{-1}$).
\\
For example, this is true when  $K$ is bounded on $\mathcal{H}$ or $T$ is selfadjoint.
By the Laplace inversion formula (\ref{2-16}), the semigroup is given as
\begin{equation}
(e^{\mathrm{i} Tt} \psi, \phi)
= \frac{1}{2\pi \mathrm{i} } \lim_{x\to \infty} \int^{x-\mathrm{i} y}_{-x-\mathrm{i} y}\!
e^{\mathrm{i} \lambda t} ((\lambda -T)^{-1} \psi, \phi) d\lambda , \quad x,y \in \mathbf{R},
\label{3-53}
\end{equation}
where the contour is a horizontal line in the lower half plane below the spectrum of $T$.
In Sec.2, we have shown that if there is an eigenvalue of $T$ on the lower half plane,
$e^{\mathrm{i} Tt}$ diverges as $t\to \infty$, while if there are no eigenvalues,
to investigate the asymptotic behavior of $e^{\mathrm{i} Tt}$ is difficult in general.
Let us show that resonance poles induce an exponential decay of the semigroup.

\begin{figure}
\begin{center}
\includegraphics{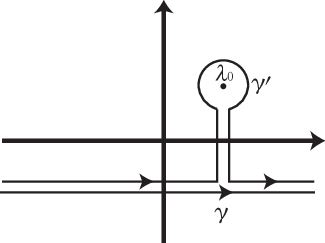}
\caption{ Deformation of the contour. \label{fig2}}
\end{center}
\end{figure}

We use the residue theorem to calculate Eq.(\ref{3-53}).
Let $\lambda _0 \in \Omega $ be an isolated resonance pole of finite multiplicity.
Suppose that the contour $\gamma $ is deformed to the contour $\gamma '$, which lies above $\lambda _0$,
without passing the generalized spectrum $\hat{\sigma}(T)$ except for $\lambda _0$, see Fig.\ref{fig2}.
For example, it is possible under the assumptions of Thm.3.19.
Recall that if $\psi, \phi \in X(\Omega )$, 
$((\lambda -T)^{-1} \psi, \phi)$ defined on the lower half plane
has an analytic continuation $\langle \mathcal{R}_\lambda \psi \,|\, \phi \rangle$ defined on 
$\Omega \cup I \cup \{ \lambda \, | \, \mathrm{Im} (\lambda ) < 0\}$ (Thm.3.12).
Thus we obtain
\begin{equation}
(e^{\mathrm{i} Tt}\psi, \phi)
 = \frac{1}{2\pi \mathrm{i} } \int_{\gamma '}\! 
       e^{\mathrm{i} \lambda t} \langle \mathcal{R}_\lambda \psi \,|\,  \phi \rangle d\lambda 
 - \frac{1}{2\pi \mathrm{i} } \int_{\gamma _0}\! e^{\mathrm{i} \lambda t}
       \langle \mathcal{R}_\lambda \psi \,|\,  \phi \rangle d\lambda ,
\label{3-55}
\end{equation}
where $\gamma _0$ is a sufficiently small simple closed curve enclosing $\lambda _0$.
Let $\mathcal{R}_\lambda = \sum^\infty_{j=-M} (\lambda _0 - \lambda )^jE_j$ be a Laurent
series of $\mathcal{R}_\lambda $ as the proof of Thm.3.16.
Due to Eq.(\ref{3-43}) and $E_{-1} = -\Pi_0$, we obtain
\begin{eqnarray*}
\frac{1}{2\pi \mathrm{i} } \int_{\gamma _0}\! e^{\mathrm{i} \lambda t}
\langle \mathcal{R}_\lambda \psi \,|\,  \phi \rangle d\lambda
 = \sum^{M-1}_{k=0} e^{\mathrm{i} \lambda _0t}\frac{(-\mathrm{i} t)^k}{k!} 
\langle  (\lambda _0 - T^\times)^ k \Pi_0 \psi \,|\,  \phi \rangle,
\label{3-56}
\end{eqnarray*}
where $\Pi_0$ is the generalized projection to the generalized eigenspace of $\lambda _0$.
Since $\mathrm{Im} (\lambda _0) > 0$, 
this proves that the second term in the right hand side of Eq.(\ref{3-55}) decays to zero as $t\to \infty$.
Such an exponential decay (of a part of) the semigroup induced by resonance poles is known as Landau damping
in plasma physics \cite{Cra}, and is often observed for Schr\"{o}dinger operators \cite{Reed}.
A similar calculation is possible without defining the generalized resolvent and the generalized spectrum
as long as the quantity $((\lambda -T)^{-1}\psi, \phi)$ has an analytic continuation for some $\psi$ and $\phi$.
Indeed, this has been done in the literature.

Let us reformulate it by using the dual space to find a decaying state corresponding to $\lambda _0$.
For this purpose, we suppose that
\\[0.2cm]
\textbf{(S2)} the semigroup $\{ (e^{\mathrm{i} Tt})^* \}_{t\geq 0}$ is an equicontinuous $C_0$ semigroup on $X(\Omega )$.
\\[0.2cm]
Then, by the theorem in IX-13 of Yosida \cite{Yos}, the dual semigroup 
$(e^{\mathrm{i} Tt})^\times = ((e^{\mathrm{i} \,T t})^{*})'$ is also an equicontinuous $C_0$ semigroup generated by $\mathrm{i}  T^\times$.
A convenient sufficient condition for (S2) is that:
\\[0.2cm]
\textbf{(S2)'} $K^*|_{X(\Omega )}$ is bounded and 
$\{ e^{\mathrm{i} H t}\}_{t\geq 0}$ is an equicontinuous $C_0$ semigroup on $X(\Omega )$.
\\[0.2cm]
Indeed, the perturbation theory of equicontinuous $C_0$ semigroups \cite{Sin} shows that (S2)' implies (S2).
By using the dual semigroup, Eq.(\ref{3-53}) is rewritten as
\begin{equation}
(e^{\mathrm{i} Tt})^\times  \psi
= \frac{1}{2\pi \mathrm{i} } \lim_{x\to \infty} \int^{x-\mathrm{i} y}_{-x-\mathrm{i} y}\!
 e^{\mathrm{i} \, \lambda t} \mathcal{R}_\lambda \psi d\lambda .
\label{3-58}
\end{equation}
for any $\psi \in iX(\Omega )$.
Similarly, Eq.(\ref{3-55}) yields
\begin{eqnarray}
(e^{\mathrm{i} Tt})^\times \psi 
 = \frac{1}{2\pi \mathrm{i} } \int_{\gamma '}\! e^{\mathrm{i} \, \lambda t} \mathcal{R}_\lambda  \psi d\lambda 
 -  \sum^{M-1}_{k=0} e^{\mathrm{i} \, \lambda_0t}\frac{(\mathrm{i} t)^k}{k!} (\lambda _0 - T^\times) ^ k \Pi_0 \psi,
\label{3-59}
\end{eqnarray}
when $\lambda _0$ is a generalized eigenvalue of finite multiplicity.
For the dual semigroup, the following statements hold.
\\[0.2cm]
\textbf{Proposition \thedef.}
Suppose (S1) and (S2).
\\
(i) A solution of the initial value problem
\begin{equation}
\frac{d}{dt}\xi = \mathrm{i} T^\times \xi, \quad \xi (0) = \mu\in \mathsf{D}(T^\times),
\label{3-60}
\end{equation}
in $X(\Omega )'$ is uniquely given by $\xi (t) = (e^{\mathrm{i} Tt})^\times \mu$.
\\
(ii) Let $\lambda _0$ be a generalized eigenvalue and $\mu_0$ a corresponding generalized eigenfunction.
Then, $(e^{\mathrm{i} Tt})^\times \mu_0 = e^{\mathrm{i} \, \lambda _0 t} \mu_0$.
\\
(iii) Let $\Pi_0$ be a generalized projection for $\lambda _0$.
The space $\Pi_0iX(\Omega )$ is $(e^{\mathrm{i} Tt})^\times$-invariant: 
$(e^{\mathrm{i} Tt})^\times \Pi_0 = \Pi_0 (e^{\mathrm{i} Tt})^\times |_{iX(\Omega )}$.
\\[0.2cm]
\textbf{Proof.} 
Since $\{ (e^{\mathrm{i} Tt})^\times \}_{t\geq 0}$ is an equicontinuous $C_0$ semigroup generated by $\mathrm{i} T^\times$,
(i) follows from the usual semigroup theory \cite{Yos}.
Because of Thm.3.5, we have $\mathrm{i} T^\times \mu_0 = \mathrm{i} \, \lambda _0 \mu_0$.
Then, 
\begin{eqnarray*}
\frac{d}{dt}  e^{\mathrm{i} \, \lambda _0 t}\mu_0
=  \mathrm{i} \, \lambda _0 e^{\mathrm{i} \, \lambda _0 t} \mu_0
= \mathrm{i} T^\times (e^{\mathrm{i} \, \lambda _0 t} \mu_0). 
\end{eqnarray*}
Thus $\xi (t) = e^{\mathrm{i} \, \lambda _0 t}\mu_0$ is a solution of the equation (\ref{3-60}).
By the uniqueness of a solution, we obtain (ii).
Because of Prop.3.13 (iii), we have
\begin{eqnarray*}
& & \frac{d}{dt}(e^{\mathrm{i} Tt})^\times\mathcal{R}_\lambda 
 = \mathrm{i}  T^\times \left( (e^{\mathrm{i} Tt})^\times \mathcal{R}_\lambda  \right), \\
& & \frac{d}{dt}\mathcal{R}_\lambda (e^{\mathrm{i} Tt})^\times|_{iY}
 = \mathcal{R}_\lambda \cdot (e^{\mathrm{i} Tt})^\times \mathrm{i} T^\times |_{iY}
 = \mathrm{i}  T^\times \left(  \mathcal{R}_\lambda (e^{\mathrm{i} Tt})^\times \right)|_{iY}.
\end{eqnarray*}
Hence, both of $(e^{\mathrm{i} Tt})^\times\mathcal{R}_\lambda $ and $\mathcal{R}_\lambda (e^{\mathrm{i} Tt})^\times$
are solutions of the equation (\ref{3-60}).
By the uniqueness, we obtain 
$(e^{\mathrm{i} Tt})^\times\mathcal{R}_\lambda |_{iY}= \mathcal{R}_\lambda (e^{\mathrm{i} Tt})^\times |_{iY}$.
Then, the definition of the projection $\Pi_0$ proves 
$(e^{\mathrm{i} Tt})^\times \Pi_0|_{iY} = \Pi_0 (e^{\mathrm{i} Tt})^\times |_{iY}$
with the aid of Eq.(\ref{3-23b}).
Since $Y$ is dense in $X(\Omega )$ and both operators $(e^{\mathrm{i} Tt})^\times \Pi_0 \circ i$ and
$\Pi_0 (e^{\mathrm{i} Tt})^\times \circ i = \Pi_0 \circ i \circ e^{\mathrm{i} Tt}$ are continuous on $X(\Omega )$,
the equality is true on $iX(\Omega )$.
 \hfill $\blacksquare$
\\[-0.2cm]

By Prop.3.14, any usual function $\phi \in X(\Omega )$ is decomposed as $\langle \phi | = \mu_1 + \mu_2$
with $\mu_1 \in \Pi_0 iX(\Omega )$ and $\mu_2 \in (id - \Pi_0) iX(\Omega )$ in the dual space.
Due to Prop.3.23 (iii) above, this decomposition is $(e^{\mathrm{i} Tt})^\times$-invariant.
When $\lambda _0 \in \Omega $, $ (e^{\mathrm{i} Tt})^\times \mu_1 \in \Pi_0 iX(\Omega )$ decays to zero exponentially as $t\to \infty$.
Eq.(\ref{3-59}) gives the decomposition explicitly.
Such an exponential decay can be well observed if we choose a function, which is sufficiently close to the 
generalized eigenfunction $\mu_0$, as an initial state.
Since $X(\Omega )$ is dense in $X(\Omega )'$ and since $(e^{\mathrm{i} Tt})^\times$ is continuous,
for any $T>0$ and $\varepsilon >0$, there exists a function $\phi_0$ in $X(\Omega )$ such that
\begin{eqnarray*}
| \langle (e^{\mathrm{i} Tt})^\times  \phi_0 \,|\,\psi \rangle  
    - \langle (e^{\mathrm{i} Tt})^\times \mu_0 \,|\, \psi \rangle | < \varepsilon,
\end{eqnarray*}
for $0 \leq t \leq T$ and $\psi \in X(\Omega )$.
This implies that
\begin{equation}
(e^{\mathrm{i} Tt} \phi_0, \psi) \sim \langle (e^{\mathrm{i} Tt})^\times \mu_0 \,|\, \psi \rangle
 = e^{\mathrm{i} \lambda _0 t} \langle \mu_0 \,|\, \psi \rangle,
\end{equation} 
for the interval $0 \leq t \leq T$.
Thus generalized eigenvalues describe the transient behavior of solutions.


\section{An application}

Let us apply the present theory to the dynamics of an infinite dimensional coupled oscillators.
The results in this section are partially obtained in \cite{Chi}.


\subsection{The Kuramoto model}

\begin{figure}
\begin{center}
\includegraphics[]{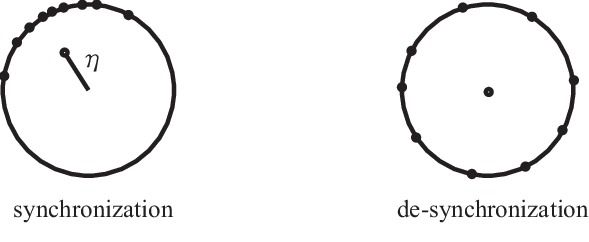}
\caption[]{The order parameter of the Kuramoto model.}
\label{fig3}
\end{center}
\end{figure}

Coupled oscillators are often used as models of collective synchronization phenomena.
One of the important models for synchronization is the Kuramoto model defined by
\begin{equation}
\frac{d\theta _i}{dt} 
= \omega _i + \frac{k}{N} \sum^N_{j=1} \sin (\theta _j - \theta _i),\,\, i= 1, \cdots  ,N,
\label{4-1}
\end{equation}
where $\theta _i = \theta _i(t) \in [ 0, 2\pi )$ denotes the phase of an $i$-th oscillator rotating on a circle,
$\omega _i\in \mathbf{R}$ is a constant called a natural frequency, $k\geq 0$ is a coupling strength,
and where $N$ is the number of oscillators.
When $k>0$, there are interactions between oscillators and collective behavior may appear.
For this system, the order parameter $\eta (t)$, which gives the centroid of oscillators, is defined to be
\begin{equation}
\eta (t) := \frac{1}{N}\sum^N_{j=1} e^{\mathrm{i}  \theta _j(t)}.
\label{4-2}
\end{equation}
If $|\eta (t)|$ takes a positive number, synchronous state is formed,
while if $|\eta (t)|$ is zero on time average, de-synchronization is stable (see Fig.\ref{fig3}).

For many applications, $N$ is too large so that statistical-mechanical description is applied.
In such a case, the continuous limit of the Kuramoto model is often employed:
At first, note that Eq.(\ref{4-1}) can be written as 
\begin{eqnarray*}
\frac{d\theta _i}{dt} 
= \omega _i + \frac{k}{2\mathrm{i} }(\eta (t) e^{-\mathrm{i} \theta_i }-\overline{\eta (t)}
       e^{\mathrm{i} \theta_i }).
\end{eqnarray*}
Keeping it in mind, the continuous model is defined as the equation of continuity of the form
\begin{eqnarray}
\left\{ \begin{array}{ll}
\displaystyle \frac{\partial \rho_t}{\partial t} + 
\frac{\partial }{\partial \theta }
\left(v \rho_t \right) = 0,  \\[0.3cm]
\displaystyle  v:= \omega  + \frac{k}{2\mathrm{i} }(\eta (t) e^{-\mathrm{i} \theta }-\overline{\eta (t)}
       e^{\mathrm{i} \theta }), \\[0.3cm]
\displaystyle \eta (t) := \int_{\mathbf{R}}
 \! \int^{2\pi}_{0} \! e^{\mathrm{i} \theta } \rho_t (\theta , \omega ) g(\omega ) d\theta d\omega .
\end{array} \right.
\label{4-3}
\end{eqnarray}
This is an evolution equation of a probability measure $\rho_t = \rho_t (\theta , \omega )$ 
on $S^1 = [0, 2\pi )$ parameterized by $t \in \mathbf{R}$ and $\omega \in \mathbf{R}$.
Roughly speaking, $\rho_t (\theta , \omega )$ denotes a probability that
an oscillator having a natural frequency $\omega $ is placed at a position $\theta $.
The $\eta$ above is the continuous version of (\ref{4-2}), which is also called the order parameter,
and $g(\omega )$ is a given probability density function for natural frequencies.
This system is regarded as a Fokker-Planck equation of (\ref{4-1}).
Indeed, it is known that the order parameter (\ref{4-2}) for the finite dimensional system converges to that
of the continuous model as $N\to \infty$ in some probabilistic sense \cite{Chi2}.
To investigate the stability and bifurcations of solutions of the system (\ref{4-3}) is a 
famous difficult problem in this field \cite{Chi, Str1}.
It is numerically observed that when $k > 0$ is sufficiently small, then the de-synchronous state
$|\eta | = 0$ is asymptotically stable, while if $k$ exceeds a certain value $k_c$, a nontrivial solution
corresponding to the synchronous state $|\eta | > 0$ bifurcates from the de-synchronous state.
Indeed, Kuramoto conjectured that

\textbf{Kuramoto conjecture \cite{Kura2}.}

\begin{figure}
\begin{center}
\includegraphics[]{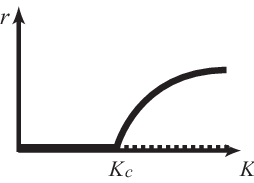}
\caption[]{A bifurcation diagram of the order parameter.
Solid lines denote stable solutions and dotted lines denote unstable solutions.}
\label{fig3-2}
\end{center}
\end{figure}

Suppose that natural frequencies $\omega _i$'s are distributed according to a probability
density function $g(\omega )$.
If $g(\omega )$ is an even and unimodal function such that $g''(0)\neq 0$, then the bifurcation
diagram of $r = |\eta |$ is given as Fig.\ref{fig3-2}; that is, if the coupling strength $k$ is 
smaller than $k_c := 2/(\pi g(0))$, then $r \equiv 0$ is asymptotically stable.
On the other hand, if $k$ is larger than $k_c$, the synchronous state emerges; there exists a positive constant $r_c$ such that 
$r = r_c$ is asymptotically stable.
Near the transition point $k_c$, $r_c$ is of order $O((k- k_c)^{1/2})$.

A function $g(\omega )$ is called unimodal (at $\omega =0$) if $g(\omega _1) > g(\omega _2)$ for $0\leq \omega _1 < \omega _2$
and $g(\omega _1) < g(\omega _2)$ for $\omega _1 < \omega _2 \leq 0$.
See \cite{Kura2} and \cite{Str1} for Kuramoto's discussion.
The purpose here is to prove the linear stability of the de-synchronous state $|\eta | = r= 0$ for $0 < k<k_c$ by
applying our spectral theory when $g(\omega ) = e^{-\omega ^2/2}/\sqrt{2\pi}$ is assumed to be the 
Gaussian distribution as in the most literature.
See Chiba \cite{Chi} for nonlinear analysis and the proof of the bifurcation at $k=k_c$.

At first, let us observe that the difficulty of the conjecture is caused by the continuous spectrum.
Let 
\begin{equation}
Z_j(t, \omega ) := \int^{2\pi}_{0} \! e^{\mathrm{i} j\theta } \rho_t (\theta , \omega ) d\theta 
\label{4-4}
\end{equation}
be the Fourier coefficient of $\rho_t (\theta , \omega )$.
Then, $Z_0(t, \omega ) = 1$ and $Z_j$ satisfy the 
differential equations
\begin{eqnarray}
\frac{dZ_1}{dt} = \mathrm{i} \omega Z_1 + \frac{k}{2} \eta (t) - \frac{k}{2} \overline{\eta (t)} Z_{2},
\label{4-5}
\end{eqnarray}
and
\begin{eqnarray}
\frac{dZ_j}{dt} &=& j\mathrm{i} \omega Z_j + \frac{jk}{2}( \eta (t) Z_{j-1} - \overline{\eta (t)} Z_{j+1}),
\label{4-6}
\end{eqnarray}
for $j=2,3,\cdots $.
Let $L^2(\mathbf{R}, g(\omega )d\omega )$ be the weighted Lebesgue space and put
$P_0(\omega ) : = 1 \in L^2(\mathbf{R}, g(\omega )d\omega )$.
Then, the order parameter is written as $\eta (t) = (Z_1, P_0)$ by using the inner product on
$L^2(\mathbf{R}, g(\omega )d\omega )$.
Since our purpose is to investigate the dynamics of the order parameter, 
let us consider the linearized system of $Z_1$ given by
\begin{equation}
\frac{dZ_1}{dt} = \left( \mathrm{i}  \mathcal{M} + \frac{k}{2} \mathcal{P}\right) Z_1, 
\label{4-7}
\end{equation}
where $\mathcal{M}: \phi(\omega ) \mapsto \omega \phi(\omega )$ is the multiplication operator 
on $L^2 (\mathbf{R}, g(\omega )d\omega )$ and 
$\mathcal{P}$ is the projection on $L^2 (\mathbf{R}, g(\omega )d\omega )$ defined to be
\begin{equation}
\mathcal{P}\phi (\omega ) = \int_{\mathbf{R}} \! \phi (\omega )g(\omega ) d\omega = (\phi, P_0)P_0. 
\end{equation} 
To determine the linear stability of the de-synchronous state $\eta = 0$, we have to investigate
the spectrum and the semigroup of the operator 
$\displaystyle T_1 := \mathrm{i}  \mathcal{M} + \frac{k}{2} \mathcal{P}$.


\subsection{Eigenvalues of the operator $T_1$}

The domain of $T_1 =\mathrm{i}  \mathcal{M} + \frac{k}{2} \mathcal{P}$ is given by 
$\mathsf{D}(\mathcal{M}) \cap \mathsf{D}(\mathcal{P}) = \mathsf{D}(\mathcal{M})$, which is dense in $L^2 (\mathbf{R}, g(\omega )d\omega)$.
Since $\mathcal{M}$ is selfadjoint and since $\mathcal{P}$ is bounded, $T_1$ is a closed operator \cite{Kato}.
Let $\mathfrak{\varrho} (T_1)$ be the resolvent set of $T_1$ and 
$\sigma (T_1) = \mathbf{C}\backslash \mathfrak{\varrho} (T_1)$ the spectrum.
Let $\sigma _p(T_1)$ and $\sigma _c(T_1)$ be the point spectrum (the set of eigenvalues)
and the continuous spectrum of $T_1$, respectively.
\\[0.2cm]
\textbf{Lemma \thedef.}\, (i) Eigenvalues $\lambda $ of $T_1$ are given as roots of
\begin{equation}
\int_{\mathbf{R}} \! \frac{1}{\lambda - \mathrm{i}  \omega }g(\omega )d\omega = \frac{2}{k}.
\label{4-9}
\end{equation} 
(ii) The continuous spectrum of $T_1$ is given by
\begin{equation}
\sigma _c(T_1) = \sigma (\mathrm{i}  \mathcal{M}) = \mathrm{i} \mathbf{R}.
\end{equation}
\textbf{Proof.} \, Part (i) follows from a straightforward calculation of the equation $\lambda v = T_1v$.
Indeed, this equation yields
\begin{eqnarray*}
(\lambda -\mathrm{i} \omega )v = \frac{k}{2} \mathcal{P}v = \frac{k}{2} \cdot (v, P_0)P_0.
\end{eqnarray*}
This is rewritten as $v = k/2 \cdot (v, P_0)(\lambda -\mathrm{i} \omega )^{-1}P_0$.
Taking the inner product with $P_0$, we obtain
\begin{eqnarray*}
1 = \frac{k}{2}((\lambda -\mathrm{i} \omega )^{-1}P_0, P_0),
\end{eqnarray*}
which gives the desired result.
Part (ii) follows from the fact that the essential spectrum is stable under the 
bounded perturbation.
The essential spectrum of $T_1$ is the same as $\sigma (\mathrm{i}  \mathcal{M})$.
Since $\mathcal{M}$ is defined on the weighted Lebesgue space and the weight $g(\omega )$ is the Gaussian,
$\sigma (\mathrm{i}  \mathcal{M}) = \mathrm{i} \cdot \mathrm{supp} (g) = \mathrm{i} \mathbf{R}$.
\hfill $\blacksquare$

Our next task is to calculate roots of Eq.(\ref{4-9}) to obtain eigenvalues of $T_1$.
Put $\displaystyle k_c = \frac{2}{\pi g(0)}$, which is called Kuramoto's transition point.
\\[0.2cm]
\textbf{Lemma \thedef.} When $k$ is larger than $k_c$, there exists a unique eigenvalue $\lambda (k)$
of $T_1$ on the positive real axis.
As $k$ decreases, the eigenvalue $\lambda (k)$ approaches to the imaginary axis, and at $k=k_c$,
it is absorbed into the continuous spectrum and disappears.
When $0 < k< k_c$, there are no eigenvalues.
\\[0.2cm]
\textbf{Proof.} Put $\lambda = x + \mathrm{i} y$ with $ x,y\in \mathbf{R}$,
Eq.(\ref{4-9}) is rewritten as
\begin{equation}
\left\{ \begin{array}{l}
\displaystyle 
\int_{\mathbf{R}} \! \frac{x}{x^2 + (\omega -y)^2}g(\omega )d\omega  = \frac{2}{k},\\[0.4cm]
\displaystyle 
\int_{\mathbf{R}} \! \frac{\omega -y}{x^2 + (\omega -y)^2}g(\omega )d\omega  = 0. \\
\end{array} \right.
\end{equation}
The first equation implies that if there is an eigenvalue $x+\mathrm{i} y$ for $k>0$, then $x>0$.
Next, the second equation is calculated as
\begin{eqnarray*}
0 = \int_{\mathbf{R}} \! \frac{\omega -y}{x^2 + (\omega -y)^2} g(\omega )d\omega 
 = \int^\infty_{0} \! \frac{\omega }{x^2 + \omega ^2} (g(y + \omega ) - g(y - \omega )) d\omega . 
\end{eqnarray*}
Since $g$ is an even function, $y = 0$ is a root of this equation.
Since $g$ is unimodal, $g(y + \omega ) - g(y - \omega ) > 0$ when $y < 0, \omega >0$
and $g(y + \omega ) - g(y - \omega ) <0$ when $y > 0, \omega >0$. Hence, $y=0$ is a unique root.
This proves that an eigenvalue should be on the positive real axis, if it exists.

Let us show the existence.
If $|\lambda |$ is large, Eq.(\ref{4-9}) is expanded as
\begin{eqnarray*}
\frac{1}{\lambda } + O(\frac{1}{\lambda ^2}) = \frac{2}{k}.
\end{eqnarray*}
Thus Rouch\'{e}'s theorem proves that Eq.(\ref{4-9}) has a root $\lambda \sim k/2$ if $k>0$ is sufficiently large.
Its position $\lambda (k)$ is continuous (actually analytic) in $k$ as long as it exists.
The eigenvalue disappears only when $\lambda \to +0$ as $k\to k_c$ for some value $k_c$.
Substituting $y=0$ and taking the limit $x \to +0,\, k\to k_c$, we have 
\begin{eqnarray*}
\lim_{x\to +0} \int_{\mathbf{R}} \! \frac{x}{x^2 + \omega ^2}g(\omega )d\omega = \frac{2}{k_c}
\end{eqnarray*}
The well known formula
\begin{eqnarray*}
\lim_{x\to +0} \int_{\mathbf{R}} \! \frac{x}{x^2 + \omega ^2}g(\omega )d\omega = \pi g(0)
\end{eqnarray*}
provides $\displaystyle k_c =2/\pi g(0)$.
Since $k_c$ is uniquely determined, the eigenvalue $\lambda (k)$ exists for $k>k_c$,
disappears at $k=k_c$ and there are no eigenvalues for $0<k<k_c$.
\hfill $\blacksquare$

This lemma shows that when $k$ is larger than $k_c$,
$Z_1 = 0$ of the equation (\ref{4-7}) is unstable because of the eigenvalue with a positive real part.
However, when $0< k< k_c$, there are no eigenvalues and the spectrum of $T_1$ consists of the 
continuous spectrum on the imaginary axis.
Hence, the usual spectral theory does not provide the stability of solutions.
To handle this difficulty, let us introduce a rigged Hilbert space.


\subsection{A rigged Hilbert space for $T_1$}

To apply our theory, let us define a test function space $X(\Omega )$.
Let $\mathrm{Exp}_+(\beta, n)$ be the set of holomorphic functions on the region 
$\mathbf{C}_n := \{ z\in \mathbf{C} \, | \, \mathrm{Im} (z) \geq -1/n\}$ such that the norm
\begin{equation}
|| \phi ||_{\beta, n} := \sup_{\mathrm{Im} (z) \geq -1/n} e^{-\beta |z|}|\phi (z)|
\label{4-11}
\end{equation}
is finite. With this norm, $\mathrm{Exp}_+(\beta, n)$ is a Banach space.
Let $\mathrm{Exp}_+ (\beta)$ be their inductive limit with respect to $n=1,2,\cdots $
\begin{equation}
\mathrm{Exp}_+ (\beta) = \varinjlim_{n \geq 1} \mathrm{Exp}_+ (\beta, n) = \bigcup_{n \geq 1} \mathrm{Exp}_+ (\beta, n).
\label{4-12}
\end{equation}
Next, define $\mathrm{Exp}_+ $ to be their inductive limit with respect to $\beta = 0,1,2,\cdots $
\begin{equation}
\mathrm{Exp}_+ = \varinjlim_{\beta \geq 0} \mathrm{Exp}_+ (\beta) = \bigcup_{\beta \geq 0} \mathrm{Exp}_+ (\beta).
\label{4-13}
\end{equation}
Thus $\mathrm{Exp}_+$ is the set of holomorphic functions near the upper half plane that can grow at most exponentially.
Then, we can prove the next proposition.
\\[0.2cm]
\textbf{Proposition \thedef.} $\mathrm{Exp}_+$ is a topological vector space satisfying
\\[0.2cm]
(i) $\mathrm{Exp}_+$ is a complete Montel space (see Sec.3.1 for Montel spaces).
\\
(ii) $\mathrm{Exp}_+$ is a dense subspace of $L^2(\mathbf{R}, g(\omega )d\omega )$.
\\
(iii) the topology of $\mathrm{Exp}_+$ is stronger than that of $L^2(\mathbf{R}, g(\omega )d\omega )$.
\\
(iv) the operators $\mathcal{M}$ and $\mathcal{P}$ are continuous on $\mathrm{Exp}_+$.
In particular, $T_1 : \mathrm{Exp}_+ \to \mathrm{Exp}_+$ is continuous (note that it is not continuous on 
$L^2(\mathbf{R}, g(\omega )d\omega )$).
\\[0.2cm]
See \cite{Chi} for the proof.
Thus, $X(\Omega ) := \mathrm{Exp}_+$ satisfies (X1) to (X3) and the rigged Hilbert space
\begin{equation}
\mathrm{Exp}_+ \subset L^2(\mathbf{R}, g(\omega )d\omega )\subset \mathrm{Exp}_+'
\end{equation}
is well-defined.
Furthermore, the operator
\begin{equation}
T := T_1/\mathrm{i}  = \mathcal{M} + \frac{k}{2\mathrm{i} }\mathcal{P}
\end{equation}
satisfies the assumptions (X4) to (X8) with $H = \mathcal{M}$ and $K = \frac{k}{2\mathrm{i} }\mathcal{P}$.
Indeed, the analytic continuation $A(\lambda )$ of the resolvent $(\lambda -\mathcal{M})^{-1}$ is given by
\begin{equation}
\langle A(\lambda )\psi \,|\, \phi \rangle = \left\{ \begin{array}{ll}
\displaystyle 
     \int_{\mathbf{R}}\! \frac{1}{\lambda -\omega  } \psi (\omega ) \phi (\omega ) g(\omega )d\omega  
           + 2\pi \mathrm{i}  \psi (\lambda ) \phi (\lambda )g(\lambda ) & (\mathrm{Im}(\lambda ) > 0), \\[0.4cm]
\displaystyle  \lim_{y\to -0} 
  \int_{\mathbf{R}}\! \frac{1}{x + \mathrm{i} y -\omega } \psi (\omega ) \phi (\omega )g(\omega ) d\omega  & (x = \lambda \in \mathbf{R}), \\[0.4cm]
\displaystyle  \int_{\mathbf{R}}\! \frac{1}{\lambda -\omega } \psi (\omega ) \phi (\omega )g(\omega ) d\omega
 & (\mathrm{Im}(\lambda ) < 0),
\end{array} \right.
\label{4-16}
\end{equation}
for $\psi, \phi\in \mathrm{Exp}_+$.
Since functions in $\mathrm{Exp}_+$ are holomorphic near the upper half plane, (X4) and (X5) are satisfied with
$I = \mathbf{R}$ and $\Omega  = \text{(the upper half plane)}$.
Since $\mathcal{M}$ and $\mathcal{P}$ are continuous on $\mathrm{Exp}_+$, (X6) and (X7) are satisfied 
with $Y = \mathrm{Exp}_+$.
For (X8), note that the dual operator $K^\times$ of $K$ is given as
\begin{equation}
K^\times \mu = \frac{k}{2\mathrm{i} } \langle \mu \,|\, P_0 \rangle \langle P_0 \,| \,\,\in i \mathrm{Exp}_+= iX(\Omega ).
\label{4-17}
\end{equation}
Since the range of $K^\times$ is included in $iX(\Omega )$, (X8) is satisfied.
Therefore, all assumptions in Sec.3 are verified and we can apply our spectral theory to 
the operator $T_1/ \mathrm{i} $.
\\[0.2cm]
\textbf{Remark.}
$T_1$ is not continuous on $\mathrm{Exp}_+(\beta ,n)$ for fixed $\beta > 0$ because of the multiplication 
$\mathcal{M}: \phi \mapsto \omega \phi$.
The inductive limit in $\beta$ is introduced so that it becomes continuous.
The proof of Lemma 4.1 shows that the eigenfunction of $T_1$ associated with $\lambda $ is given by
\begin{eqnarray*}
v_\lambda = \frac{1}{\lambda -\mathrm{i}\omega }, \quad \lambda >0.
\end{eqnarray*}
If $\lambda >0$ is small, $v_\lambda $ is not included in $\mathrm{Exp}_+(\beta ,n)$ for fixed $n$.
The inductive limit in $n$ is introduced so that any eigenfunctions are elements of $\mathrm{Exp}_+$.
Furthermore, the topology of $\mathrm{Exp}_+$ is carefully defined so that
the strong dual $\mathrm{Exp}_+'$ becomes a Fr\'{e}chet Montel space.
It is known that the strong dual of a Montel space is also Montel.
Since $\mathrm{Exp}_+$ is defined as the inductive limit of Banach spaces, its dual is realized as a 
projective limit of Banach spaces $\mathrm{Exp}_+(\beta, n)'$, which is  Fr\'{e}chet by the definition.
Hence, the contraction principle is applicable on $\mathrm{Exp}_+'$, which allows one to prove the 
existence of center manifolds of the system (\ref{4-3}) (see \cite{Chi}), though nonlinear problems are not
treated in this paper.


\subsection{Generalized spectrum of $T_1/\mathrm{i} $}

\begin{figure}
\begin{center}
\includegraphics[]{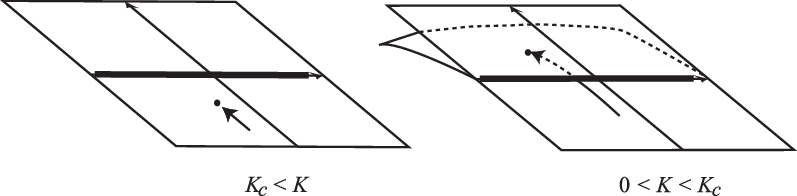}
\caption[]{As $k$ decreases, the eigenvalue of $T_1/\mathrm{i} $ disappears from the original 
complex plane by absorbed into the continuous spectrum on the real axis.
However, it still exists as a resonance pole on the second Riemann sheet of the generalized resolvent.}
\label{fig4}
\end{center}
\end{figure}

For the operator $T_1/\mathrm{i} $, we can prove that (see also Fig.\ref{fig4}) 
\\[0.2cm]
\textbf{Proposition \thedef.} 
\\
(i) The generalized continuous and the generalized residual spectra are empty.
\\
(ii) For any $k>0$, there exist infinitely many generalized eigenvalues on the upper half plane.
\\
(iii) For $k>k_c$, there exists a unique generalized eigenvalue $\lambda (k)$ on the lower half plane,
which is an eigenvalue of $T_1/\mathrm{i} $ in $L^2(\mathbf{R}, g(\omega )d\omega )$-sense.
As $k$ decreases, $\lambda (k)$ goes upward and at $k=k_c$, $\lambda (k)$ gets across the real axis and
it becomes a resonance pole.
When $0<k<k_c$, $\lambda (k)$ lies on the upper half plane and there are no generalized eigenvalues on the 
lower half plane.
\\[0.2cm]
\textbf{Proof.} 
(i) Since $K^\times$ given by (\ref{4-17}) is a one-dimensional operator,
it is easy to verify the assumption of Thm.3.19.
Hence, the generalized continuous and the generalized residual spectra are empty.

(ii) Let $\lambda $ and $\mu$ be a generalized eigenvalue and a generalized eigenfunction.
By Eq.(\ref{3-16}), $\lambda $ and $\mu$ satisfy $(id - K^\times A(\lambda ))K^\times \mu = 0$.
In our case,
\begin{eqnarray*}
\langle K^\times \mu \,|\, \phi \rangle = \frac{k}{2\mathrm{i} } \langle \mu \,|\, P_0 \rangle
\langle P_0 \,|\, \phi \rangle
\end{eqnarray*}
and
\begin{eqnarray*}
\langle K^\times A(\lambda )K^\times \mu \,|\, \phi \rangle = 
\langle A(\lambda )K^\times \mu \,|\, K^*\phi \rangle
= \left( \frac{k}{2 \mathrm{i} } \right)^2 \langle \mu \,|\, P_0 \rangle
\langle P_0 \,|\, \phi \rangle \langle A(\lambda )P_0 \,|\, P_0 \rangle,
\end{eqnarray*}
for any $\phi \in \mathrm{Exp}_+$.
Hence, generalized eigenvalues are given as roots of the equation
\begin{equation}
\frac{2\mathrm{i} }{k} = \langle A(\lambda )P_0 \,|\, P_0 \rangle
 = \left\{ \begin{array}{ll}
\displaystyle 
     \int_{\mathbf{R}}\! \frac{1}{\lambda -\omega  } g(\omega )d\omega  
           + 2\pi \mathrm{i}  g(\lambda ) & (\mathrm{Im}(\lambda ) > 0), \\[0.4cm]
\displaystyle  \int_{\mathbf{R}}\! \frac{1}{\lambda -\omega } g(\omega ) d\omega
 & (\mathrm{Im}(\lambda ) < 0).
\end{array} \right.
\label{4-18}
\end{equation}
Since $g$ is the Gaussian, it is easy to verify that the equation (\ref{4-18}) for $\mathrm{Im}(\lambda ) > 0$ has
infinitely many roots $\{ \lambda _n\}^\infty_{n=0}$ such that $\mathrm{Im}(\lambda _n) \to \infty$
and they approach to the rays $\mathrm{arg} (z) = \pi/4,\, 3\pi /4$ as $n\to \infty$.

(iii) When $\mathrm{Im}(\lambda ) < 0$, the equation (\ref{4-18}) is the same as (\ref{4-9}),
in which $\lambda $ is replaced by $\mathrm{i} \lambda $.
Thus Lemma 4.2 shows that when $k>k_c$, there exists a root $\lambda (k)$ on the lower half plane.
As $k$ decreases, $\lambda (k)$ goes upward and for $0 < k_c < k$, it becomes a root of the first equation of 
(\ref{4-18}) because the right hand side of (\ref{4-18}) is holomorphic in $\lambda $. \hfill $\blacksquare$
\\

Eq.(\ref{3-15}) shows that a generalized eigenfunction associated with $\lambda $ is given by
$\mu = A(\lambda )K^\times \mu = \frac{k}{2\mathrm{i} } \langle \mu \,|\, P_0 \rangle \cdot A(\lambda ) \langle P_0 \,|$.
We can choose a constant $\langle \mu \,|\, P_0 \rangle$ as $\langle \mu \,|\, P_0 \rangle = 2\mathrm{i} /k$.
Then, $\mu = A(\lambda ) \langle P_0 \,| = A(\lambda )i(P_0)$.
When $\mathrm{Im} (\lambda ) < 0$, $\mu$ is a usual function written as $\mu = (\lambda -\omega )^{-1}
\in \mathrm{Exp}_+$, although when $\mathrm{Im} (\lambda ) \geq 0$, $\mu$ is not included in 
$L^2(\mathbf{R}, g(\omega )d\omega) $ but an element of the dual space $\mathrm{Exp}_+'$.
In what follows, we denote generalized eigenvalues by $\{ \lambda _n\}^\infty_{n=0}$ such that 
$|\lambda _n| \leq |\lambda _{n+1}|$ for $n=0,1,\cdots $,
and a corresponding generalized eigenfunction by $\mu_n = A(\lambda_n ) \langle P_0 \,|$.
Thm.3.5 proves that they satisfy $T_1^\times \mu_n = \mathrm{i} \lambda _n \mu_n$.
Note that when $0<k<k_c$, all generalized eigenvalues satisfy $\mathrm{Im} (\lambda _n) > 0$.

Next, let us calculate the generalized resolvent of $T_1/\mathrm{i} $.
Eq.(\ref{3-30}) yields
\begin{eqnarray}
\mathcal{R}_\lambda \phi - A(\lambda ) K^\times \mathcal{R}_\lambda \phi = A(\lambda )\phi \quad 
\Longrightarrow \quad
 \mathcal{R}_\lambda \phi = A(\lambda )\phi 
      + \frac{k}{2\mathrm{i} } \langle \mathcal{R}_\lambda \phi \,|\, P_0 \rangle A(\lambda ) \langle P_0 \,|,
\label{4-19}
\end{eqnarray}
for any $\phi \in \mathrm{Exp}_+$.
Taking the inner product with $P_0$, we obtain
\begin{eqnarray*}
\langle \mathcal{R}_\lambda \phi \,|\, P_0 \rangle
= \frac{\langle A(\lambda ) \phi \,|\, P_0 \rangle}{1-\frac{k}{2\mathrm{i} }\langle A(\lambda ) P_0 \,|\, P_0 \rangle}
 = \frac{\langle A(\lambda )P_0\,|\, \phi \rangle}{1-\frac{k}{2\mathrm{i} }\langle A(\lambda ) P_0 \,|\, P_0 \rangle}.
\end{eqnarray*}
Substituting this into Eq.(\ref{4-19}), we obtain
\begin{eqnarray}
 \mathcal{R}_\lambda \phi = A(\lambda )\phi 
      + \left( 2\mathrm{i} /k -  \langle A(\lambda ) P_0 \,|\, P_0 \rangle \right)^{-1} 
         \langle A(\lambda )P_0\,|\, \phi \rangle \cdot A(\lambda ) \langle P_0\, |.
\label{4-20}
\end{eqnarray}
Then, the generalized Riesz projection for the generalized eigenvalue $\lambda _n$ is given by
\begin{equation}
\Pi_n \phi = \frac{1}{2\pi \mathrm{i} } \int_{\gamma }\! \mathcal{R}_\lambda \phi d\lambda 
 = D_n \langle A(\lambda _n) P_0 \,|\, \phi \rangle \cdot A(\lambda_n ) \langle P_0\, |
 = D_n \langle \mu_n \,|\, \phi \rangle \cdot \mu_n,
\label{4-21}
\end{equation}
or
\begin{equation}
\langle \Pi_n \phi \,|\, \psi \rangle = D_n \langle \mu_n \,|\, \phi \rangle \cdot 
\langle \mu_n \,|\, \psi \rangle,
\end{equation}
where $D_n$ is a constant defined by
\begin{eqnarray*}
D_n = \lim_{\lambda \to \lambda _n} (\lambda -\lambda _n) \cdot
\left( 2\mathrm{i} /k -  \langle A(\lambda ) P_0 \,|\, P_0 \rangle \right)^{-1} .
\end{eqnarray*}
As was proved in Thm.3.16, the range of $\Pi_n$ is spanned by the generalized eigenfunction $\mu_n$.


\subsection{Spectral decomposition of the semigroup}

\begin{figure}
\begin{center}
\includegraphics[]{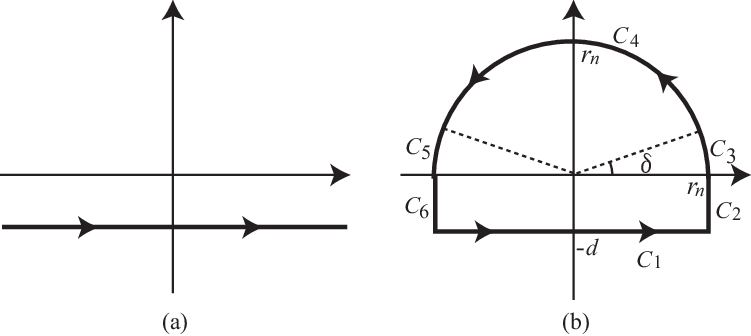}
\caption[]{The contour for the Laplace inversion formula.}
\label{fig5}
\end{center}
\end{figure}

Now we are in a position to give a spectral decomposition theorem of the semigroup generated by 
$T_1 = \mathrm{i} \mathcal{M} + \frac{k}{2} \mathcal{P}$.
Since $\mathrm{i} \mathcal{M}$ generates the $C^0$-semigroup on $L^2(\mathbf{R}, g(\omega )d\omega )$ and
$\mathcal{P}$ is bounded, $T_1$ also generates the $C^0$-semigroup given by
\begin{equation}
e^{T_1t}\phi = \lim_{y\to \infty} \frac{1}{2\pi \mathrm{i} } 
\int^{x+\mathrm{i} y}_{x-\mathrm{i} y}\! e^{\lambda t}(\lambda -T_1)^{-1} \phi d\lambda , 
\end{equation}
for $t>0$, where $x$ is a sufficiently large number.
In $L^2(\mathbf{R}, g(\omega )d\omega )$-theory, we can not deform the contour from the right half plane
to the left half plane because $T_1$ has the continuous spectrum on the imaginary axis.
Let us use the generalized resolvent $\mathcal{R}_\lambda $ of $T_1/\mathrm{i} $.
For this purpose, we rewrite the above as
\begin{equation}
e^{T_1t}\phi = \lim_{y\to \infty} \frac{1}{2\pi \mathrm{i} } 
\int^{y-\mathrm{i} x}_{-y-\mathrm{i} x}\! e^{\mathrm{i} \lambda t}(\lambda -T_1/\mathrm{i} )^{-1} \phi d\lambda , 
\end{equation}
whose contour is the horizontal line on the lower half plane (Fig.\ref{fig5} (a)).
Recall that when $\mathrm{Im} (\lambda ) < 0$, $((\lambda -T_1/\mathrm{i} )^{-1}\phi, \psi)
 = \langle \mathcal{R}_\lambda \phi \,|\, \psi \rangle$ for $\phi, \psi \in \mathrm{Exp}_+$
because of Thm.3.12.
Thus we have
\begin{equation}
\langle e^{T_1t} \phi\,|\, \psi \rangle
 = \lim_{y\to \infty} \frac{1}{2\pi \mathrm{i} } \int^{y-\mathrm{i} x}_{-y-\mathrm{i} x}\! 
e^{\mathrm{i} \lambda t}\langle \mathcal{R}_\lambda \phi \,|\, \psi \rangle d\lambda .
\end{equation}
Since $\langle \mathcal{R}_\lambda \phi \,|\, \psi \rangle$ is a meromorphic function whose poles are 
generalized eigenvalues $\{ \lambda _n\}^\infty_{n=0}$, we can deform the contour from the lower half plane 
to the upper half plane.
With the aid of the residue theorem, we can prove the next theorems.
\\[0.2cm]
\textbf{Theorem \thedef \,(Spectral decomposition).} 
\\
For any $\phi, \psi \in \mathrm{Exp}_+$, there exists $t_0 > 0$ such that the equality
\begin{equation}
\langle e^{T_1t} \phi\,|\, \psi \rangle = \sum^\infty_{n=0} D_n e^{\mathrm{i} \lambda _nt}
\langle \mu_n \,|\, \phi \rangle \cdot \langle \mu_n \,|\, \psi \rangle 
\label{4-26}
\end{equation}
holds for $t>t_0$.
Similarly, the dual semigroup $(e^{T_1t})^\times$ satisfies
\begin{equation}
(e^{T_1t})^\times \phi = \sum^\infty_{n=0} D_n e^{\mathrm{i} \lambda _nt}
\langle \mu_n \,|\, \phi \rangle \cdot \mu_n
\label{4-27}
\end{equation}
for $\phi \in \mathrm{Exp}_+$ and $t>t_0$, where the right hand side converges with respect to the 
strong dual topology on $\mathrm{Exp}_+'$.
\\[0.2cm]
\textbf{Theorem \thedef \,(Completeness).}
\\
(i) A system of generalized eigenfunctions $\{  \mu_n \}^\infty_{n=0}$
is complete in the sense that if $\langle \mu_n \,|\, \psi \rangle = 0$ for $n= 0,1,\cdots $,
then $\psi =0$.
\\
(ii) $\mu_0, \mu_1, \cdots$ are linearly independent of each other:
if $\sum^\infty_{n=0} a_n \mu_n = 0$ with $a_n\in \mathbf{C}$, then $a_n = 0$ for every $n$.
\\
(iii) The decomposition of $(e^{T_1t})^\times$ using $\{  \mu_n \}^\infty_{n=0}$ is uniquely
expressed as (\ref{4-27}).
\\[0.2cm]
\textbf{Corollary \thedef \,(Linear stability).}
\\
When $0<k<k_c$, the order parameter $\eta (t) = (Z_1, P_0)$ for the linearized system (\ref{4-7})
decays exponentially to zero as $t\to \infty$ if the initial condition is an element of $\mathrm{Exp}_+$.
\\[0.2cm]
\textbf{Proof.} We start with the proof of Corollary 4.7.
When an initial condition of the system (\ref{4-7}) is given by $\phi \in \mathrm{Exp}_+$,
the order parameter is given by $\eta (t) = (Z_1, P_0) = (e^{T_1t}\phi, P_0)$.
If $0<k<k_c$, all generalized eigenvalues lie on the upper half plane, so that
$\mathrm{Re}[\mathrm{i} \lambda _n] < 0$ for $n=0,1,\cdots $.
Then the corollary follows from Eq.(\ref{4-26}).

Next, let us prove Thm.4.6.

(i) If $\langle \mu_n  \,|\, \psi \rangle = 0$ for all $n$, Eq.(\ref{4-26}) provides
$(e^{T_1t}\phi , \psi) = (\phi , (e^{T_1t})^* \psi) =0$ for any $\phi \in \mathrm{Exp}_+$.
Since $\mathrm{Exp}_+$ is dense in $L^2(\mathbf{R}, g(\omega )d\omega )$, we obtain $(e^{T_1t})^*\psi = 0$
for any $t > t_0$, which proves $\psi = 0$. 

(ii) Suppose that $\sum^\infty_{n=0} a_n \mu_n = 0$.
Prop.3.23 provides
\begin{eqnarray*}
0 = (e^{T_1t})^\times \sum^\infty_{n=0} a_n \mu_n 
 = \sum^\infty_{n=0} a_n (e^{T_1t})^\times  \mu_n 
 = \sum^\infty_{n=0} a_n e^{\mathrm{i} \lambda _n t} \mu_n .
\end{eqnarray*}
Changing the label if necessary, we can assume that
\begin{eqnarray*}
\mathrm{Re}[\mathrm{i}  \lambda _0] \geq \mathrm{Re}[\mathrm{i}  \lambda _1] 
 \geq \mathrm{Re}[\mathrm{i}  \lambda _2]  \geq \cdots,
\end{eqnarray*}
without loss of generality.
Suppose that $\mathrm{Re}[\mathrm{i}  \lambda _0] = \cdots = \mathrm{Re}[\mathrm{i}  \lambda _k]$
and $\mathrm{Re}[\mathrm{i}  \lambda _k] > \mathrm{Re}[\mathrm{i}  \lambda _{k+1}]$.
Then, the above equality provides
\begin{eqnarray*}
0 = \sum^k_{n=0} a_n e^{\mathrm{i} \mathrm{Im}[\mathrm{i} \lambda _n] t}  \mu_n 
 + \sum^\infty_{n=k+1} a_n e^{(\mathrm{i} \lambda _n - \mathrm{Re}[\mathrm{i}  \lambda _0] )t} \mu_n .
\end{eqnarray*} 
Taking the limit $t\to \infty$ yields
\begin{eqnarray*}
0 = \lim_{t\to \infty}\sum^k_{n=0} a_n e^{\mathrm{i} \mathrm{Im}[\mathrm{i} \lambda _n] t}  \mu_n .
\end{eqnarray*}
Since the finite set $\mu_0, \cdots , \mu_k$ of eigenvectors are linearly independent
as in a finite-dimensional case, we obtain $a_n = 0$ for $n = 0, \cdots , k$.
The same procedure is repeated to prove $a_n = 0$ for every $n$.

(iii) This immediately follows from Part (ii) of the theorem.

Finally, let us prove Thm.4.5.
Recall that generalized eigenvalues are roots of the equation (\ref{4-18}).
Hence, there exist positive numbers $B$ and $\{ r_j\}^\infty_{j=1}$ such that
\begin{equation}
\left| 1 - \frac{k}{2\mathrm{i} }\langle A(\lambda )P_0 \,|\, P_0 \rangle \right| \geq B
\label{4-28}
\end{equation} 
holds for $\lambda = r_j e^{\mathrm{i} \theta }\,\, (0 < \theta < \pi)$.
Take a positive number $d$ so that $\mathrm{Im} (\lambda _n) > -d$ for all $n = 0,1,\cdots $.
Fix a small positive number $\delta $ and define a closed curve $C(j) = C_1 + \cdots + C_6$ by
\begin{eqnarray*}
& & C_1 = \{ x - \mathrm{i} d \, | \, -r_j \leq x \leq r_j\} \\
& & C_2 = \{ r_j - \mathrm{i} y \, | \, 0 \leq y\leq d \} \\
& & C_3 = \{ r_j e^{\mathrm{i} \theta } \, | \, 0 \leq \theta \leq \delta \} \\
& & C_4 = \{ r_j e^{\mathrm{i} \theta } \, | \, \delta  \leq \theta \leq \pi - \delta \},
\end{eqnarray*}
and $C_5$ and $C_6$ are defined in a similar manner to $C_3$ and $C_2$, respectively, see Fig.\ref{fig5} (b).

Let $\lambda _0, \lambda _1, \cdots , \lambda _{N(j)}$ be generalized eigenvalues inside the closed curve $C(j)$.
Due to Eq.(\ref{4-21}), we have
\begin{eqnarray*}
\frac{1}{2\pi \mathrm{i} }\int_{C(j)} e^{\mathrm{i} \lambda t} 
         \langle \mathcal{R}_\lambda \phi \,|\, \psi \rangle d\lambda 
= \sum^{N(j)}_{n=1} D_ne^{\mathrm{i}  \lambda _nt} \langle \mu_n \,|\, \phi \rangle
                       \langle \mu_n \,|\, \psi \rangle.
\end{eqnarray*}
Taking the limit $j\to \infty\,\, (r_j \to \infty)$ provides
\begin{eqnarray*}
\langle e^{T_1t} \phi\,|\, \psi \rangle 
     &=& \lim_{j\to \infty} \frac{1}{2\pi \mathrm{i} }\int_{C_2+\cdots +C_6} e^{\mathrm{i} \lambda t} 
         \langle \mathcal{R}_\lambda \phi \,|\, \psi \rangle d\lambda  \\
& & \quad \quad + \lim_{j\to \infty} \sum^{N(j)}_{n=0} D_n e^{\mathrm{i} \lambda _nt}
\langle \mu_n \,|\, \phi \rangle \cdot \langle \mu_n \,|\, \psi \rangle .
\end{eqnarray*}
We can prove by the standard way that the integrals along $C_2, C_3, C_5$ and $C_6$ tend to zero as $j\to \infty$.
The integral along $C_4$ is estimated as
\begin{eqnarray*}
\left| \int_{C_4} e^{\mathrm{i} \lambda t} \langle \mathcal{R}_\lambda \phi \,|\, \psi \rangle d\lambda  \right| 
&\leq & \max_{\lambda \in C_4} |\langle \mathcal{R}_\lambda \phi \,|\, \psi \rangle | \cdot
       \int^{\pi / 2}_{\delta }\! 2r_j e^{-r_jt \sin \theta } d\theta \\
 &\leq & \max_{\lambda \in C_4} |\langle \mathcal{R}_\lambda \phi \,|\, \psi \rangle | \cdot
       \int^{\pi / 2}_{\delta }\! 2r_j e^{-2r_jt\theta /\pi} d\theta \\
 &\leq & \max_{\lambda \in C_4} |\langle \mathcal{R}_\lambda \phi \,|\, \psi \rangle | \cdot
       \frac{\pi}{t}\left( e^{-2r_j \delta t / \pi} - e^{-r_jt}\right) .\\
\end{eqnarray*} 
It follows from (\ref{4-20}) that
\begin{eqnarray*}
& &  \langle \mathcal{R}_\lambda \phi \,|\, \psi \rangle 
= \frac{1}{2\mathrm{i} /k - \langle A(\lambda )P_0 \,|\, P_0 \rangle } \times \\
& &\quad \left( 
\frac{2\mathrm{i} }{k} \langle A(\lambda )\phi \,|\, \psi \rangle - \langle A(\lambda ) P_0\,|\, P_0 \rangle
\langle A(\lambda )\phi \,|\, \psi \rangle + \langle A(\lambda ) P_0\,|\, \phi \rangle
\langle A(\lambda ) P_0 \,|\, \psi \rangle \right).
\end{eqnarray*}
Since $\phi, \psi \in \mathrm{Exp}_+$, there exist positive constants $C_1,C_2, \beta_1, \beta_2$ such that
\begin{eqnarray*}
|\phi (\lambda )| \leq C_1 e^{\beta_1 |\lambda |} ,\quad |\psi (\lambda )| \leq C_2 e^{\beta_2 |\lambda |}.
\end{eqnarray*}
Using the definition (\ref{4-16}) of $A(\lambda )$, we can show that there exist positive constants
$D_0, \cdots ,D_4$ such that
\begin{eqnarray}
|\langle \mathcal{R}_\lambda \phi \,|\, \psi \rangle | \leq
\frac{D_0 + (D_1 + D_2C_1 e^{\beta _1 |\lambda  |} + D_3C_2 e^{\beta _2 |\lambda  |} + 
       D_4C_1C_2 e^{(\beta _1 + \beta_2 )|\lambda |}) \cdot |g(\lambda )|}
{|2\mathrm{i} /k - \langle A(\lambda )P_0 \,|\, P_0 \rangle|}.
\label{4-29}
\end{eqnarray}
When $|g(\lambda )| \to \infty$ as $|\lambda | \to \infty$, this yields
\begin{eqnarray*}
|\langle \mathcal{R}_\lambda \phi \,|\, \psi \rangle | \leq
\frac{D_0 + (D_1 + D_2C_1 e^{\beta _1 |\lambda  |} + D_3C_2 e^{\beta _2 |\lambda  |} + 
       D_4C_1C_2 e^{(\beta _1 + \beta_2 )|\lambda |}) }
{2\pi} + o(|\lambda |).
\end{eqnarray*}
When $|g(\lambda )|$ is bounded as $|\lambda | \to \infty$, Eq.(\ref{4-28}) is used to estimate (\ref{4-29}).
For both cases, we can show that there exists $D_5 > 0$ such that
\begin{eqnarray*}
|\langle \mathcal{R}_\lambda \phi \,|\, \psi \rangle | \leq D_5e^{(\beta_1 + \beta_2) r_j},
\quad (\lambda  = r_j e^{\mathrm{i} \theta }).
\end{eqnarray*}
Therefore, we obtain
\begin{eqnarray*}
\left| \int_{C_4} e^{\mathrm{i} \lambda t} \langle \mathcal{R}_\lambda \phi \,|\, \psi \rangle d\lambda  \right| 
\leq \frac{\pi D_5}{t} \left(e^{(\beta_1 + \beta_2 - 2\delta t/\pi)r_j} - e^{(\beta_1 + \beta_2 - t) r_j} \right) .
\end{eqnarray*}
Thus if $t > t_0 := \pi (\beta_1 + \beta_2)/(2\delta )$, this integral tends to zero as $j\to \infty$,
which proves Eq.(\ref{4-26}).
Since Eq.(\ref{4-26}) holds for each $\psi \in \mathrm{Exp}_+$, the right hand side of Eq.(\ref{4-27})
converges with respect to the weak dual topology on $\mathrm{Exp}_+'$.
Since $\mathrm{Exp}_+$ is a Montel space, a weakly convergent series also converges
with respect to the strong dual topology. \hfill $\blacksquare$


\appendix

\section{Pettis integrals and vector valued holomorphic functions on the dual space}

The purpose in this Appendix is to give the definition and the existence theorem of Pettis integrals.
After that, a few results on vector-valued holomorphic functions are given.
For the existence of Pettis integrals, the following property 
\\
\textbf{(CE)} for any compact set $K$, the closed convex hull of $K$ is compact,
\\
which is sometimes called the convex envelope property, is essentially used.
For the convenience of the reader, sufficient conditions for the property are listed below.
We also give conditions for $X$ to be barreled because it is assumed in (X3).
Let $X$ be a locally convex Hausdorff vector space, and $X'$ its dual space.
\begin{itemize}
\item The closed convex hull $\overline{co}(K)$ of a compact set $K$ in $X$ is compact if and only if 
$\overline{co}(K)$ is complete in the Mackey topology on $X$ (Krein's theorem, see K\"{o}the \cite{Kot}, \S 24.5).

\item $X$ has the convex envelope property if $X$ is quasi-complete.

\item If $X$ is bornological, the strong dual $X'$ is complete.
In particular, the strong dual of a metrizable space is complete.

\item If $X$ is barreled, the strong dual $X'$ is quasi-complete.
In particular, $X'$ has the convex envelope property.

\item Montel spaces, Fr\'{e}chet spaces, Banach spaces and Hilbert spaces are barreled.

\item The product, quotient, direct sum, (strict) inductive limit, completion of barreled spaces are barreled.
\end{itemize}
See Tr\'{e}ves \cite{Tre} for the proofs.

Let $X$ be a topological vector space over $\mathbf{C}$ and $(S, \mu)$ a measure space.
Let $f : S\to X$ be a measurable $X$-valued function.
If there exists a unique $I_f \in X$ such that $\langle \xi \,|\, I_f \rangle
 = \int_{S}\! \langle \xi \,|\, f \rangle d\mu$ for any $\xi \in X'$, $I_f$ is called the 
\textit{Pettis integral} of $f$.
It is known that if $X$ is a locally convex Hausdorff vector space with the convex envelope property,
$S$ is a compact Hausdorff space with a finite Borel measure $\mu$, and if $f:S\to X$ is continuous,
then the Pettis integral of $f$ exists (see Rudin \cite{Rud}).
In Sec.3.5, we have defined the integral of the form $\int_{\gamma }\! \mathcal{R}_\lambda \phi d\lambda $,
where $\mathcal{R}_\lambda \phi $ is an element of the dual $X(\Omega )'$.
Thus our purpose here is to define a ``dual version" of Pettis integrals.

In what follows, let $X$ be a locally convex Hausdorff vector space over $\mathbf{C}$,
$X'$ a strong dual with the convex envelope property,
and let $S$ be a compact Hausdorff space with a finite Borel measure $\mu$.
For our purpose in Sec.3.5, $S$ is always a closed path on the complex plane.
Let $f: S\to X'$ be a continuous function with respect to the strong dual topology on $X'$.
\\[0.2cm]
\textbf{Theorem \thedef.}
(i) Under the assumptions above, there exists a unique $I(f) \in X'$ such that
\begin{equation}
\langle I(f) \,|\, x \rangle = \int_{S}\! \langle f \,|\, x \rangle d\mu 
\label{A-1}
\end{equation}
for any $x\in X$. $I(f)$ is denoted by $I(f) = \int_{S}\! fd\mu$ and called the Pettis integral of $f$.

(ii) The mapping $f\mapsto I(f)$ is continuous in the following sense;
for any neighborhood $U$ of zero in $X'$ equipped with the weak dual topology, there exists a neighborhood $V$
of zero in $X'$ such that if $f(s) \in V$ for any $s\in S$, then $I(f)\in U$.

(iii) Furthermore, suppose that $X$ is a barreled space.
Let $T$ be a linear operator densely defined on $X$ and $T'$ its dual operator with the domain 
$\mathsf{D}(T') \subset X'$.
If $f(S) \subset \mathsf{D}(T') $ and the set $\{ \langle T' f(s) \,|\, x \rangle \}_{s\in S}$ is bounded for each $x\in X$,
then, $I(f) \in \mathsf{D}(T')$ and  $T' I(f) = I(T'f)$ holds; that is,
\begin{equation}
T' \int_{S}\! fd\mu = \int_{S}\! T' f d\mu
\label{A-2}
\end{equation}
holds.

The proof of (i) is done in a similar manner to that of the existence of Pettis integrals on $X$ \cite{Rud}.
Note that $T$ is not assumed to be continuous for the part (iii).
When $T$ is continuous, the set $\{ \langle T' f(s) \,|\, x \rangle \}_{s\in S}$ is bounded because $T'$ and $f$ are continuous.
\\[0.2cm]
\textbf{Proof.} At first, note that the mapping $\langle \cdot \,|\, x \rangle : X' \to \mathbf{C}$
is continuous because $X$ can be canonically embedded into the dual of the strong dual $X'$.
Thus $\langle f(\cdot) \,|\, x \rangle : S \to \mathbf{C}$ is continuous and it is integrable
on the compact set $S$ with respect to the Borel measure.

Let us show the uniqueness. If there are two elements $I_1(f), I_2(f) \in X'$ satisfying Eq.(\ref{A-1}),
we have $\langle I_1(f) \,|\, x \rangle = \langle I_2(f) \,|\, x \rangle$ for any $x\in X$.
By the definition of $X'$, it follows $I_1(f) = I_2(f)$.

Let us show the existence.
We can assume without loss of generality that $X$ is a vector space over $\mathbf{R}$ and $\mu$ is a probability measure.
Let $L \subset X$ be a finite set and put
\begin{equation}
V_L(f) = V_L := \{ x'\in X' \, | \, \langle x' \,|\, x \rangle = \int_{S}\! \langle f \,|\, x \rangle d\mu, \,\,  \forall x\in L\}.
\end{equation}
Since $\langle \cdot \,|\, x \rangle$ is a continuous mapping, $V_L$ is closed.
Since $f$ is continuous, $f(S)$ is compact in $X'$. Due to the convex envelope property, the closed convex hull
$\overline{co}(f(S))$ is compact.
Hence, $W_L: = V_L \cap \overline{co}(f(S))$ is also compact.
By the definition, it is obvious that $W_{L_1} \cap W_{L_2} = W_{L_1 \cup L_2}$.
Thus if we can prove that $W_L$ is not empty for any finite set $L$, a family $\{ W_L\}_{L\in \{ \text{finite set}\}}$
has the finite intersection property.
Then, $\bigcap_L W_L$ is not empty because $\overline{co}(f(S))$ is compact.
This implies that there exists $I(f) \in \bigcap_L W_L$ such that 
$\langle I(f) \,|\, x \rangle = \int_{S}\! \langle f \,|\, x \rangle d\mu$ for any $x\in X$.

Let us prove that $W_L$ is not empty for any finite set $L = \{ x_1, \cdots ,x_n\} \subset X$.
Define the mapping $\mathcal{L} : X' \to \mathbf{R}^n$ to be
\begin{eqnarray*}
\mathcal{L}(x') = \left( \langle x' \,|\, x_1 \rangle\, , \cdots , \langle x' \,|\, x_n \rangle \right).
\end{eqnarray*}
This is continuous and $\mathcal{L}(f(S))$ is compact in $\mathbf{R}^n$.
Let us show that the element
\begin{equation}
y : = \left( \int_{S}\! \langle f \,|\, x_1 \rangle d\mu\, , 
\cdots  , \int_{S}\! \langle f \,|\, x_n \rangle d\mu \right)
\end{equation}
is included in the convex hull $co (\mathcal{L}(f(S)))$ of $\mathcal{L}(f(S))$.
If otherwise, there exist real numbers $c_1, \cdots  ,c_n$ such that for any 
$(z_1, \cdots ,z_n) \in co(\mathcal{L}(f(S)))$, the inequality
\begin{eqnarray*}
\sum^n_{i=1} c_i z_i < \sum^n_{i=1} c_i y_i, \quad y = (y_1, \cdots ,y_n)
\end{eqnarray*} 
holds (this is a consequence of Hahn-Banach theorem for $\mathbf{R}^n$).
In particular, since $\mathcal{L}(f(S)) \subset co (\mathcal{L}(f(S)))$,
\begin{eqnarray*}
\sum^n_{i=1} c_i \langle f \,|\, x_i \rangle < \sum^n_{i=1} c_i y_i.
\end{eqnarray*}
Integrating both sides (in the usual sense) yields $\sum^n_{i=1} c_i y_i < \sum^n_{i=1} c_i y_i$.
This is a contradiction, and therefore $y\in co (\mathcal{L}(f(S)))$.
Since $\mathcal{L}$ is linear, there exists $v\in co (f(S))$ such that $y = \mathcal{L}(v)$.
This implies that $v\in V_L \cap co (f(S))$, and thus $W_L$ is not empty. 
By the uniqueness, $\bigcap_L W_L = \{ I(f)\}$.
Part (ii) of the theorem immediately follows from Eq.(\ref{A-1}) and properties of the usual integral.

Next, let us show Eq.(\ref{A-2}).
When $X$ is a barreled space, $I(f)$ is included in $\mathsf{D}(T') $ so that $T'I(f)$ is well defined.
To prove this, it is sufficient to show that the mapping
\begin{eqnarray*}
x\mapsto \langle I(f) \,|\,Tx  \rangle 
= \int_{S}\! \langle f \,|\, Tx \rangle d\mu = \int_{S}\! \langle T'f \,|\, x \rangle d\mu  
\end{eqnarray*}
from $\mathsf{D}(T) \subset X$ into $\mathbf{C}$ is continuous.
By the assumption, the set $\{ \langle T'f(s) \,|\, x \rangle\}_{s\in S}$
is bounded for each $x\in X$.
Then, Banach-Steinhaus theorem implies that the family $\{ T'f(s)\}_{s\in S}$ of continuous linear functionals
are equicontinuous.
Hence, for any $\varepsilon >0$, there exists a neighborhood $U$ of zero in $X$ such that
$|\langle T'f(s) \,|\, x \rangle| < \varepsilon $ for any $s\in S$ and $x\in U$.
This proves that the above mapping is continuous, so that $I(f)\in \mathsf{D}(T') $
and $T'I(f) = T' \bigcap_L W_L$.

For a finite set $L \subset X$, put
\begin{eqnarray*}
& & V_L(T'f) = \{ x' \in X' \, | \langle x' \,|\, x \rangle = \int_{S}\! \langle T'f \,|\, x \rangle d\mu, 
\,\, \forall x\in L \}, \\
& & T'V_{TL}(f) = \{ T'x' \in X' \, | \, x'\in \mathsf{D}(T'),\,\, 
\langle x' \,|\, x \rangle = \int_{S}\! \langle f \,|\, x \rangle d\mu, 
\,\, \forall x\in TL \}.
\end{eqnarray*}
Put $W_L (f) = V_L(f) \cap \overline{co}(f(S))$ as before.
It is obvious that $\bigcap _L W_L (f) \subset \bigcap _L W_{TL} (f)$.
Therefore, 
\begin{eqnarray*}
\{ T' I(f)\} &=& T' \bigcap _L W_{L} (f) \subset T'\bigcap _L W_{TL} (f)\cap \mathsf{D}(T')  \\
&\subset & T' \bigcap_{L} \left( V_{TL} (f) \cap \overline{co}(f(S)) \cap \mathsf{D}(T') \right) \\
&\subset & \bigcap_{L} \left(  T'V_{TL} (f) \cap  T'\overline{co}(f(S)) \cap \mathsf{R}(T') \right).
\end{eqnarray*}
On the other hand, if $y'\in T'V_{TL}(f)$, there exists $x' \in X'$ such that $y' = T'x'$ and 
$\langle x' \,|\, x \rangle = \int_{S}\! \langle f \,|\, x \rangle d\mu$ for any $x\in TL$.
Then, for any $x\in L\cap \mathsf{D}(T)$,
\begin{eqnarray*}
\langle y' \,|\, x \rangle = \langle T'x'  \,|\, x \rangle
 = \langle x' \,|\, Tx \rangle = \int_{S}\! \langle f \,|\, Tx \rangle d\mu
 = \int_{S}\! \langle T'f \,|\, x \rangle d\mu.
\end{eqnarray*}
This implies that $y' \in V_{L\cap \mathsf{D}(T)} (T'f)$, and thus $T'V_{TL}(f) \subset V_{L\cap \mathsf{D}(T)}(T'f)$.
Hence, we obtain
\begin{eqnarray*}
\{ T' I(f)\} &\subset &
\bigcap_{L}  V_{L\cap \mathsf{D}(T)}(T'f) \cap  \overline{co}(T'f(S))
 = \bigcap _L W_{L\cap \mathsf{D}(T)} (T'f). 
\end{eqnarray*}
If $\langle x' \,|\, x \rangle = \int_{S}\! \langle f \,|\, x \rangle d\mu$ for dense subset of $X$,
then it holds for any $x\in X $.
Hence, we have
\begin{equation}
\{ I(T'f)\} = \bigcap_L W_L (T'f) = \bigcap_L W_{L\cap \mathsf{D}(T)} (T'f) \supset \{ T' I(f)\}.
\end{equation}
which proves $T' I(f) = I(T'f)$. \hfill $\blacksquare$
\\[-0.2cm]

Now that we can define the Pettis integral on the dual space, we can develop the ``dual version" of the theory
of holomorphic functions. Let $X$ and $X'$ be as in Thm.A.1.
Let $f : D\to X'$ be an $X'$-valued function on an open set $D\subset \mathbf{C}$.
\\[0.2cm]
\textbf{Definition \thedef.}
(i) $f$ is called weakly holomorphic if $\langle f \,|\, x \rangle$ is holomorphic on $D$ in the 
classical sense for any $x\in X$ (more exactly, it should be called weak-dual-holomorphic).
\\
(ii) $f$ is called strongly holomorphic if
\begin{equation}
\lim_{z_0 \to z} \frac{1}{z_0 - z} \left( f(z_0)-f(z)\right),\quad (\text{the strong dual limit})
\end{equation}
exists in $X'$ for any $z\in D$ (more exactly, it should be called strong-dual-holomorphic).
\\[0.2cm]
\textbf{Theorem \thedef.} Suppose that the strong dual $X'$ satisfies the convex envelope property
and $f : D\to X'$ is weakly holomorphic.
\\[0.2cm]
(i) If $f$ is strongly continuous, Cauchy integral formula and Cauchy integral theorem hold:
\begin{eqnarray*}
& &  f(z) = \frac{1}{2\pi \mathrm{i} }\int_{\gamma }\! \frac{f(z_0)}{z_0 - z}dz_0, \quad
 \int_{\gamma }\! f(z_0)dz_0 = 0 ,
\end{eqnarray*}
where $\gamma \subset D$ is a closed curve enclosing $z\in D$.
\\
(ii) If $f$ is strongly continuous and if $X'$ is quasi-complete, 
$f$ is strongly holomorphic and is of $C^\infty$ class.
\\
(iii) If $X$ is barreled, the weak holomorphy implies the strong continuity.
Thus (i) and (ii) above hold; $f$ is strongly holomorphic and is expanded in a Taylor series as
\begin{equation}
f(z) = \sum^\infty_{n=0} \frac{f^{(n)}(a)}{n!}(z-a)^n, \quad (\text{strong dual convergence}),
\label{A-6}
\end{equation}
near $a\in D$.
Similarly, a Laurent expansion and the residue theorem hold if $f$ has an isolated singularity.
\\[0.2cm]
\textbf{Proof.}
(i) Since $f$ is continuous with respect to the strong dual topology, the Pettis integral
\begin{eqnarray*}
I(z) = \frac{1}{2\pi \mathrm{i} }\int_{\gamma }\! \frac{f(z_0)}{z_0 - z}dz_0 
\end{eqnarray*}
exists. By the definition of the integral,
\begin{eqnarray*}
\langle I(z) \,|\, x \rangle
 = \frac{1}{2\pi \mathrm{i} } \int_{\gamma }\! \frac{\langle f(z_0) \,|\, x \rangle}{z_0 - z} dz_0 
\end{eqnarray*}
for any $x\in X$.
Since $\langle f(z) \,|\, x \rangle$ is holomorphic in the usual sense, the right hand side above is equal to 
$\langle f(z) \,|\, x \rangle$. Thus we obtain $I(z) = f(z)$,
which gives the Cauchy formula. The Cauchy theorem also follows from the classical one.

(ii) Let us prove that $f$ is strongly holomorphic at $z_0$.
Suppose that $z_0 = 0$ and $f(z_0) = 0$ for simplicity.
By the same way as above, we can verify that
\begin{eqnarray*}
\frac{f(z)}{z} &=& \frac{1}{2\pi \mathrm{i} } \int_{\gamma }\! \frac{f(z_0)}{z_0(z_0 - z)} dz_0 \\
&=& \frac{1}{2\pi \mathrm{i} } \int_{\gamma }\! \frac{f(z_0)}{z_0^2} dz_0
     + \frac{z}{2\pi \mathrm{i} } \int_{\gamma }\! \frac{f(z_0)}{z_0^2(z_0 - z)} dz_0.
\end{eqnarray*}
Since $X'$ is quasi-complete, the above converges as $z\to 0$ to yield
\begin{eqnarray*}
f'(0) := \lim_{z\to 0} \frac{f(z)}{z} = \frac{1}{2\pi \mathrm{i} } \int_{\gamma }\! \frac{f(z_0)}{z_0^2} dz_0.
\end{eqnarray*}
In a similar manner, we can verify that
\begin{equation}
f^{(n)}(z):=\frac{d^n}{dz^n}f(z) = \frac{n!}{2\pi \mathrm{i} } \int_{\gamma }\! \frac{f(z_0)}{(z_0 - z)^{n+1}} dz_0
\label{A-7}
\end{equation}
exists for any $n= 0,1,2,\cdots $.

(iii) If $X$ is barreled, weakly bounded sets in $X'$ are strongly bounded (see Thm.33.2 of Tr\'{e}ves \cite{Tre}).
By using it, let us prove that a weakly holomorphic $f$ is strongly continuous.
Suppose that $f(0) = 0$ for simplicity.
Since $\langle f(z) \,|\, x \rangle$ is holomorphic in the usual sense, Cauchy formula provides
\begin{eqnarray*}
\frac{\langle f(z) \,|\, x \rangle}{z} 
= \frac{1}{2\pi \mathrm{i} } \int_{\gamma }\! \frac{1}{z_0 - z} \frac{\langle f(z_0) \,|\, x \rangle}{z_0} dz_0.
\end{eqnarray*}
Suppose that $|z| < \delta $ and $\gamma $ is a circle of radius $2\delta $ centered at the origin.
Since $\langle f(\cdot) \,|\, x \rangle $ is holomorphic, there exists a positive number $M$
such that $|\langle f(z_0) \,|\, x \rangle| < M$ for any $z_0\in \gamma $.
Then,
\begin{eqnarray*}
\left| \frac{\langle f(z)  \,|\, x \rangle}{z}  \right| 
\leq \frac{1}{2\pi}\cdot \frac{1}{\delta }\cdot \frac{M}{2\delta } \cdot 4\pi \delta = \frac{M}{\delta }.
\end{eqnarray*}
This shows that the set $B:= \{ f(z)/z \,\, | \,\, |z| < \delta \}$ is weakly bounded in $X'$.
Since $X$ is barreled, $B$ is strongly bounded.
By the definition of bounded sets, for any convex balanced neighborhood $U$ of zero in $X'$ equipped with the strong dual,
there is a number $t > 0$ such that $tB \subset U$.
This proves that
\begin{eqnarray*}
f(z) - f(0) = f(z) \in \frac{z}{t} U \subset \frac{\delta }{t}U
\end{eqnarray*}
for $|z-0|< \delta $, which implies the continuity of $f$ with resect to the strong dual topology.

If $X$ is barreled, $X'$ is quasi-complete and has the convex envelope property.
Thus the results in (i) and (ii) hold.

Finally, let us show that $f(z)$ is expanded in a Taylor series around $a\in D$.
Suppose $a = 0$ for simplicity.
Let us prove that 
\begin{eqnarray*}
S_m = \sum^m_{n=0} \frac{1}{n!} \frac{d^n f}{dz^n}(0) z^n
\end{eqnarray*}
forms a Cauchy sequence with respect to the strong dual topology.
It follows from (\ref{A-7}) that
\begin{eqnarray*}
\frac{1}{n!} \langle f^{(n)}(0)  \,|\, x \rangle
 = \frac{1}{2\pi \mathrm{i} } \int_{\gamma }\! \frac{\langle  f(z_0) \,|\, x \rangle }{z_0^{n+1}} dz_0
\end{eqnarray*}
for any $x\in X$.
Suppose that $\gamma $ is a circle of radius $2\delta $ centered at the origin.
There exists a constant $M_x > 0$ such that $|\langle f(z_0) \,|\, x \rangle| < M_x$ for any $z_0\in \gamma $,
which implies that the set $\{ f(z_0) \, | \, z_0\in \gamma \} $ is weakly bounded.
Because $X$ is barreled, it is strongly bounded.
Therefore, for any bounded set $B\subset X$, there is a positive number $M_B$ such that 
$|\langle  f(z_0)  \,|\, x \rangle | < M_B$ for $x\in B$ and $z_0\in \gamma $.
Then, we obtain
\begin{eqnarray*}
\left| \frac{1}{n!} \langle f^{(n)}(0) \,|\, x \rangle \right|
\leq \frac{1}{2\pi} \cdot \frac{M_B}{(2\delta )^{n+1}} \cdot 4\pi \delta  = \frac{M_B}{(2\delta )^n}.
\end{eqnarray*}
By using this, it is easy to verify that $\{ \langle S_m \,|\, x \rangle\} ^\infty_{m=0}$
is a Cauchy sequence uniformly in $x\in B$ when $|z| < \delta $.
Since $X'$ is quasi-complete, $S_m$ converges as $m\to \infty$ in the strong dual topology.
By the Taylor expansion in the classical sense, we obtain
\begin{eqnarray*}
\langle f(z)\,|\, x \rangle =\sum^\infty_{n=0}\frac{1}{n!}\frac{d^n}{dz_0^n}\Bigl|_{z_0=0} \langle f(z_0) \,|\, x \rangle z^n
=\sum^\infty_{n=0} \frac{1}{n!} \langle  f^{(n)}(0) \,|\, x\rangle z^n.
\end{eqnarray*}
Since $\lim_{m\to \infty}S_m$ exists and $\langle \, \cdot \, \,|\, x \rangle : X' \to \mathbf{C}$ is continuous, we have
\begin{eqnarray*}
\langle f(z)  \,|\, x\rangle = \langle \sum^\infty_{n=0} \frac{1}{n!} f^{(n)}(0)z^n \,|\, x \rangle,
\end{eqnarray*}
for any $x\in X$.
This proves Eq.(\ref{A-6}) for $a = 0$. 
The proof of a Laurent expansion, when $f$ has an isolated singularity, is done in the same way.
Then, the proof of the residue theorem immediately follows from the classical one.
\hfill $\blacksquare$
\\[0.2cm]
\textbf{Remark.}
In a well known theory of Pettis integrals on a space $X$ \cite{Rud}, not a dual $X'$,
we need not assume that $X$ is barreled because every locally convex space $X$ has the property that
any weakly bounded set is bounded with respect to the original topology.
Since the dual $X'$ does not have this property, we have to assume that $X$ is barreled so that 
any weakly bounded set in $X'$ is strongly bounded.


\vspace*{0.5cm}
\textbf{Acknowledgements.}

This work was supported by Grant-in-Aid for Young Scientists (B), No.22740069 from MEXT Japan.


\end{document}